\documentclass[leqno,12pt]{article}
\usepackage{amsmath,amssymb,wasysym,rotating,a4wide,color}

\usepackage[all,cmtip]{xy}

\usepackage{verbatim}

\newdir^{ (}{{}*!/-3pt/\dir^{(}}
\newdir_{ (}{{}*!/-3pt/\dir_{(}}

\entrymodifiers={+!!<0pt,\fontdimen22\textfont2>}

\def\Langle{\langle\!\langle}
\def\Rangle{\rangle\!\rangle}
\def\sgn{\mathop{\rm sgn}\nolimits}
\def\bigtimes{\mathop{\raise-2pt\hbox{\huge$\times$}}}

\newbox\circbulletbox
\setbox\circbulletbox=\hbox{\,{\large$\circledcirc$}\kern-8.7pt\raise .6pt\hbox{$\bullet$}\,}
\def\circbullet{\copy\circbulletbox}

\let\le\leqslant
\let\ge\geqslant

\def\mycirc{{\kern1pt\circ\kern2pt}}

\def\cycl{\mathop{\rm cycl}\nolimits}

\def\Aut{\mathop{\rm Aut}\nolimits}
\def\Frob{\mathop{\rm Frob}\nolimits}

\def\Gal{\mathop{\rm Gal}\nolimits}

\def\Spec{\mathop{\rm Spec}\nolimits}

\def\mod{\mathop{\rm mod}\nolimits}
\def\Ker{\mathop{\rm Ker}\nolimits}

\def\Mat{\mathop{\rm Mat}\nolimits}

\def\Cent{\mathop{\rm Cent}\nolimits}
\def\Norm{\mathop{\rm Norm}\nolimits}
\def\Conjug{\mathop{\rm Conj}\nolimits}

\def\diag{\mathop{\rm diag}\nolimits}
\def\anti{\mathop{\rm anti}\nolimits}
\def\ord{\mathop{\rm ord}\nolimits}

\def\proj{\mathop{\rm pr}\nolimits}
\def\GL{\mathop{\rm GL}\nolimits}

\def\ab{{\rm ab}}

\def\et{{\rm \acute{e}t}}

\def\geom{{\rm geom}}
\def\arith{{\rm arith}}

\def\id{{\rm id}}

\let\phi\varphi
\let\theta\vartheta
\let\epsilon\varepsilon
\let\setminus\smallsetminus

\newtheorem{Thm}{Theorem}[section]
\newtheorem{Prop}[Thm]{Proposition}
\newtheorem{Lem}[Thm]{Lemma}

\newtheorem{Cor}[Thm]{Corollary}

\newtheorem{Rem}[Thm]{Remark}
\newtheorem{Ex}[Thm]{Example}

\newtheorem{Class}[Thm]{Classification}

\numberwithin{Thm}{section}
\def\UseTheoremCounterForNextEquation{\setcounter{equation}{\value{Thm}}\addtocounter{Thm}{1}}
\renewcommand{\theequation}
             {\arabic{section}.\arabic{Thm}}

\def\qed{{\hskip0pt\unskip\unskip\nobreak\hfil\penalty50
          \hskip1em\hbox{}\nobreak\hfil
           {$\square$}
          \parfillskip=0pt\finalhyphendemerits=0
          \par}\medskip}

\newenvironment{Proof}
               {\noindent{\bf Proof.}\ }
               {\qed}

               {\noindent{\bf Proof of #1.}\ }
               {\qed}

\newcommand{\BA}{{\mathbb{A}}}

\newcommand{\BC}{{\mathbb{C}}}

\newcommand{\BF}{{\mathbb{F}}}

\newcommand{\BP}{{\mathbb{P}}}
\newcommand{\BQ}{{\mathbb{Q}}}

\newcommand{\BZ}{{\mathbb{Z}}}

\newbox\mybox
\def\arrover#1{\mathrel{
       \setbox\mybox=\hbox spread 1.4em
              {\hfil$\scriptstyle#1$\hfil}
       \vbox{\offinterlineskip\copy\mybox
             \hbox to\wd\mybox{\rightarrowfill}}}}

\def\larrover#1{\mathrel{
       \setbox\mybox=\hbox spread 1.4em
              {\hfil$\scriptstyle#1\vphantom{g}$\hfil}
       \vbox{\offinterlineskip\copy\mybox
             \hbox to\wd\mybox{\leftarrowfill}}}}

\def\ontoover#1{\mathrel{
       \setbox\mybox=\hbox spread 1.4em
              {\hfil$\scriptstyle#1\vphantom{g}$\hfil}
       \vbox{\offinterlineskip\copy\mybox
             \hbox to\wd\mybox{\rightarrowfill\hskip-2.8mm
                               $\rightarrow$}}}}
\def\leftontoover#1{\mathrel{
       \setbox\mybox=\hbox spread 1.4em
              {\hfil$\scriptstyle#1\vphantom{g}$\hfil}
       \vbox{\offinterlineskip\copy\mybox
             \hbox to\wd\mybox{$\leftarrow$\hskip-2.8mm
                               \leftarrowfill}}}}
\let\longto\longrightarrow
\let\into\hookrightarrow
\let\onto\twoheadrightarrow
\def\longonto{\ontoover{\ }}

\def\invlim{\mathop{\vtop{\hbox{\rm lim}\vskip-8pt
        \hbox{\hskip1pt$\scriptstyle\longleftarrow$}\vskip-1pt}}}

\begin{document}

\title{Profinite iterated monodromy groups\\
arising from quadratic polynomials}

\author{Richard Pink\\[12pt]
\small Department of Mathematics \\[-3pt]
\small ETH Z\"urich\\[-3pt]
\small 8092 Z\"urich\\[-3pt]
\small Switzerland \\[-3pt]
\small pink@math.ethz.ch\\[12pt]}


\date{September 23, 2013}

\maketitle

\begin{abstract}
We study in detail the profinite group $G$ arising as geometric \'etale iterated monodromy group of an arbitrary quadratic polynomial over a field of characteristic different from two. This is a self-similar closed subgroup of the group of automorphisms of a regular rooted binary tree. 
(When the base field is $\BC$ it is the closure of the finitely generated iterated monodromy group for the usual topology which is also often studied.)
Among other things we prove that the conjugacy class and hence the isomorphism class of $G$ depends only on the combinatorial type of the postcritical orbit of the polynomial. 

We represent a chosen instance of $G$ by explicit recursively defined generators. The uniqueness up to conjugacy depends on a certain semirigidity property, which ensures that arbitrary conjugates of these generators under the automorphism group of the tree always generate a subgroup that is conjugate to~$G$. We determine the Hausdorff dimension, the maximal abelian factor group, and the normalizer of~$G$ using further explicit generators. The description of the normalizer is then used to describe the arithmetic \'etale iterated monodromy group of the quadratic polynomial. 

The methods used are purely group theoretical and do not involve fundamental groups over $\BC$ at all.
\end{abstract}

{\renewcommand{\thefootnote}{}
}

\newpage
\tableofcontents


\newpage
\addtocounter{section}{-1}
\section{Introduction}
\label{Intro}

Let $f$ be a rational function of degree two in one variable, with coefficients in a field $k$ of characteristic different from~$2$. For every integer $n\ge1$ let $f^n$ denote the $n^{\rm th}$ \emph{iterate of~$f$}, obtained by substituting $n$ copies of $f$ into each other. The result is a rational function of degree~$2^n$, which we view as a finite morphism from the projective line $\BP^1_k$ over $k$ to itself. 

Let $\bar k$ be a separable closure of~$k$. Then $f$ has precisely two critical points in $\BP^1(\bar k)$; let $C$ denote the set of these. Let $P := \bigcup_{n\ge1}f^n(C)$ denote the (strictly) \emph{postcritical orbit} of~$f$.
Then each $f^n$ is a connected ramified covering of $\BP^1_k$ which is unramified over $\BP^1_k\setminus P$. As such it is determined up to isomorphism by the action of the \'etale fundamental group $\pi_1^\et(\BP^1_k \setminus P)$ on a certain set with $2^n$ elements. Varying~$n$, these coverings form an infinite tower $\ldots \smash{\stackrel{f}{\longto}}\BP^1_k \smash{\stackrel{f}{\longto}}\BP^1_k \smash{\stackrel{f}{\longto}}\BP^1_k$ whose structure up to isomorphism is determined by the action of the \'etale fundamental group on a regular rooted binary tree~$T$. The action corresponds to a homomorphism $\rho\colon \pi_1^\et(\BP^1_k \setminus P) \to\Aut(T)$, whose image is the \emph{iterated monodromy group} associated to $f$ which we study in this article.

By construction it is a profinite group. If $k$ is the field of complex numbers one can also define a discrete iterated monodromy group using the usual fundamental group of $\BP^1(\BC) \setminus P$, whose closure in the automorphism group of the tree is the profinite iterated monodromy group defined above. In the present article we disregard that point of view entirely and work purely algebraically and with profinite groups only.

In particular we do not require $k$ to be separably closed. In fact, 
it is natural to compare the image of $\pi_1^\et(\BP^1_k \setminus P)$ with the image of $\pi_1^\et(\BP^1_{\bar k} \setminus P)$. 
We call the former the \emph{arithmetic} iterated monodromy group $G^\arith$ and the latter the \emph{geometric} iterated monodromy group $G^\geom$ associated to~$f$. Then $G^\geom$ is a closed normal subgroup of $G^\arith$, and the monodromy action induces a natural surjective homomorphism $\bar\rho\colon \Gal(\bar k/k)\onto G^\arith/G^\geom$. Determining them is a fundamental problem in the area.

Our ultimate motivation for this, though not developed at all in the present article, is to understand the arithmetic properties of the tower of coverings which manifest themselves in the iterated monodromy groups. We are intrigued by the analogy between the Galois representation on the tree associated to~$f$, and the linear $\ell$-adic Galois representation on the Tate module of an abelian variety or a Drinfeld module or on the $\ell$-adic cohomology of an algebraic variety. There are many deep results for the latter, and it would be delightful if counterparts for iterated monodromy groups could be developed as well. Natural candidates for study are for instance the images of Frobenius elements in the case that $k$ is a number field or a finite field.%
\footnote[1]{For example: It is not hard to show that the set of conjugacy classes in the automorphism group of the tree has the cardinality of the continuum. On the other hand, the conjugacy class of the image of a Frobenius element depends only on a quadratic polynomial over a finite field and an element of that finite field. Thus the conjugacy classes arising from Frobenius elements form a countable set only. The images of Frobenius elements must therefore satisfy some strong special properties. The question is: which? 

Note that for a non-trivial linear $\ell$-adic Galois representation the set of conjugacy classes in $\GL_n(\BQ_\ell)$ likewise has the same cardinality as~$\BQ_\ell$, but the characteristic polynomials of Frobenius elements usually have coefficients in a number field, which also leaves only countably many conjugacy classes for them.}

\medskip
To find generators of $G^\geom$ we look at inertia groups. Note that the monodromy action is ramified at any point $p\in P$. Since $G^\geom$ is a pro-$2$-group and the base field does not have characteristic~$2$, the image of an inertia group at $p$ is a quotient of the tame inertia group and therefore topologically generated by a single element. Let $b_p\in G^\geom$ be a generator. Then $G^\geom$ is topologically generated by the elements $b_p$ for all $p\in P$. In particular, if the postcritical orbit $P$ is finite, then $G^\geom$ is topologically finitely generated. Also, the $b_p$ can be chosen such that their product in some given order is~$1$.

If we remove the bottom level from the tower of coverings $\ldots \smash{\stackrel{f}{\longto}}\BP^1_k \smash{\stackrel{f}{\longto}}\BP^1_k \smash{\stackrel{f}{\longto}}\BP^1_k$ we obtain the same tower again. From this one deduces that the monodromy representation and the iterated monodromy groups are \emph{self-similar} in the sense that a copy arises as a proper subquotient. This implies a weak form of recursion relations for the generators~$b_p$. Namely, on deleting the root, the tree $T$ decomposes into two regular rooted binary trees, each related to the tower of coverings with the bottom level removed. Choosing an isomorphism from $T$ to each of these subtrees, it turns out that each $b_p$ is conjugate in $\Aut(T)$ to an explicit element given in terms of certain other generators acting on the two half trees. Specifically, let $\sigma\in \Aut(T)$ be an element of order $2$ that interchanges the two half trees. For any two elements $u,v\in\Aut(T)$ let $(u,v)$ denote the element of $\Aut(T)$ which acts by $u$ on the first half tree and by $v$ on the second. Then $b_p$ is conjugate under $\Aut(T)$ to 
\UseTheoremCounterForNextEquation
\begin{equation}\label{0GenRecConj1}
\qquad\left\{\begin{array}{ll}
\ \ \ \sigma & \hbox{if $p=f(c)$ for some $c\in C\setminus P$,}\\[3pt]
(b_c,1)\,\sigma & \hbox{if $p=f(c)$ for some $c\in C\cap P$,}\\[3pt]
(b_q,1) & \hbox{if $p=f(q)$ for a unique $q\in P\setminus C$,} \\[3pt]
(b_q,b_{q'}) & \hbox{if $p=f(q)=f(q')$ for distinct $q,q'\in P\setminus C$,}
\end{array}\right\}
\end{equation}
see Proposition \ref{17GenRec1}.
All the results about $G^\geom$ in this article are purely algebraic consequences of these facts.

\medskip
{}From now on assume that $f$ is a quadratic polynomial. Then one of the critical points is $\infty$ and satisfies $f(\infty)=\infty$. The other critical point is then necessarily a point $p_0\in k$. Abbreviate $p_i := f^i(p_0)$ for all $i\ge1$, so that $P=\{\infty\}\sqcup P'$ for the  strict forward orbit $P' := \{p_i\mid i\ge1\}$ of $p_0$.
If $P$ is infinite, the points $p_i$ for $i\ge0$ are all distinct, and one can easily deduce from (\ref{0GenRecConj1}) that $G^\geom=G^\arith=\Aut(T)$: see Theorem \ref{18GgeomarithThm}.

So assume henceforth that $P$ and hence $P'$ is finite. Let $r$ denote the cardinality of~$P'$. Then the points $p_i$ for $1\le i\le r$ are all distinct, and $p_{r+1} = p_{s+1}$ for some unique $0\le s<r$. If $s=0$, we also have $p_r=p_0$ and the critical point $p_0$ is \emph{periodic}; otherwise it is \emph{strictly pre-periodic}: see Classification \ref{17Cases}. 

As some product of the $b_p$ for all $p\in P$ is~$1$, the group $G^\geom$ is already topologically generated by the elements $b_{p_i}$ for $1\le i\le r$. Moreover, the recursion relations (\ref{0GenRecConj1}) for the generators $b_{p_i}$ involve only generators of the form $b_{p_j}$. On $P'$ they thus reduce to
\UseTheoremCounterForNextEquation
\begin{equation}\label{0GenRecConj2}
\hbox{$b_{p_i}$ is conjugate under $\Aut(T)$ to} \ \ 
\left\{\begin{array}{ll}
\ \ \ \sigma & \hbox{if $i=1$ and $s>0$,}\\[3pt]
(b_{p_r},1)\,\sigma & \hbox{if $i=1$ and $s=0$,}\\[3pt]
(b_{p_{i-1}},1) & \hbox{if $1<i\le r$ and $i\not=s+1$,} \\[3pt]
(b_{p_s},b_{p_r}) & \hbox{if $1<i\le r$ and $i=s+1$,} 
\end{array}\right\}
\end{equation}
see Proposition \ref{17GenRec}. One of our main results is:

\begin{Thm}\label{0Semirigid}
Any two closed subgroups of $\Aut(T)$ which are generated by elements $b_{p_i}$ for $1\le i\le r$ which satisfy the relations (\ref{0GenRecConj2}) 
are conjugate. 
\end{Thm}

Thus $G^\geom$ up to conjugacy depends only on the combinatorics of the postcritical orbit of~$f$. This may be somewhat surprising, because the relations seem so much weaker than the recursion relations one usually has for the generators of discrete iterated monodromy groups over~$\BC$, which are definite equalities not just up to conjugacy. 

By Theorem \ref{0Semirigid} the analysis of $G^\geom$ reduces to the study of the closed subgroup $G\subset\Aut(T)$ that is topologically generated by any single choice of generators for all $1\le i\le r$ which satisfy (\ref{0GenRecConj2}). As generators we choose the elements $a_i\in\Aut(T)$ that are uniquely characterized by the recursion relations
\UseTheoremCounterForNextEquation
\begin{equation}\label{0GenRecEq}
a_i\ =\ \left\{\begin{array}{ll}
\ \ \ \sigma & \hbox{if $i=1$ and $s>0$,}\\[3pt]
(a_r,1)\,\sigma & \hbox{if $i=1$ and $s=0$,}\\[3pt]
(a_{i-1},1) & \hbox{if $1<i\le r$ and $i\not=s+1$,} \\[3pt]
(a_s,a_r) & \hbox{if $1<i\le r$ and $i=s+1$.}
\end{array}\right\}
\kern-20pt
\end{equation}
Understanding $G$ is then a purely group theoretical problem that is interesting in its own right. Most of the article is actually devoted to it. The methods used are standard finite and profinite group theory, and time and again the exploitation of the self-similarity properties of $G$ resulting from the recursion relations (\ref{0GenRecEq}). 
As there is a marked difference in structural details between the periodic and the strictly pre-periodic case, the main body of the paper is divided accordingly. Among other things:

\begin{enumerate}
\item[$\bullet$] We determine the maximal abelian factor group of~$G$: see Theorems \ref{22Gab} and \ref{32Gab}.

\item[$\bullet$] We calculate the Hausdorff dimension of~$G$: see Theorems \ref{23Hausdorff} and \ref{33Hausdorff}.

\item[$\bullet$] We prove Theorem \ref{0Semirigid}. Actually we show that for any elements $b_{p_i}\in\Aut(T)$ satis\-fy\-ing the weak recursion relations (\ref{0GenRecConj2}), 
on conjugating them by the same element of $\Aut(T)$ we can achieve that all $b_{p_i}$ lie in~$G$, are conjugate to the respective $a_i$ under~$G$, and topologically generate~$G$: see Theorem \ref{24SemiRigid} or \ref{34SemiRigid}. We view this result as a kind of \emph{semirigidity} property of the generators. 

\item[$\bullet$] In one interesting case of a group with two generators, namely the closure of the `Basilica group', we prove a stronger \emph{rigidity} property that on conjugating the $b_{p_i}$ by the same element of $\Aut(T)$ we can make them equal to~$a_i$\,: see Theorem \ref{25SmallRigid}. This rigidity property, and its failure whenever the number of generators is greater than two, was in fact the motivation for the author to search for a weaker version that would still give $G$ up to conjugacy in the general case, and which he found in the semirigidity property described above.

\item[$\bullet$] We determine the normalizer $N\subset\Aut(T)$ of $G$ and describe it using further explicit generators. This result requires significantly more effort than the others. In the periodic case we establish an explicit isomorphism $N/G \cong (\BZ_2^\times)^r$, where $\BZ_2$ denotes the ring of $2$-adic integers: see Theorem \ref{25NormThm}. 
In the strictly pre-periodic case with $r=2$ the group $G$ is infinite pro-dihedral isomorphic to $\BZ_2\rtimes\{\pm1\}$ and its normalizer is isomorphic to $\BZ_2\rtimes\BZ_2^\times$: see Subsection \ref{37Small}.
In the strictly pre-periodic case with $r\ge3$ we construct an isomorphism between $N/G$ and a countably infinite product $\BF_2^{\raise2pt\hbox{$\kern1pt\scriptstyle\infty$}}$  of copies of the cyclic group of order~$2$. For technical reasons this requires a further division into subcases: see Theorems \ref{35NormThm}, \ref{36NormThm}, and \ref{37NormThm}. 

\item[$\bullet$] We study a particularly useful class of elements of $G$ called \emph{odometers}. We show that their proportion in $G$ is $2^{-r}>0$ and that they are all conjugate under~$N$: see Proposition \ref{28ManyOdos} or \ref{38ManyOdos} and Theorem \ref{28OdosConj} or \ref{38OdosConj}. We describe how the normalizer of any odometer in $G$ sits inside $N$ and how it maps to $N/G$: see Proposition \ref{28OdoNorm} or \ref{38OdoNorm}. In the periodic case we can also strengthen the semirigidity theorem mentioned above to include odometers in Theorem \ref{28OdoSemiRigid}, but we do not know how to achieve anything similar in the strictly pre-periodic case.
\end{enumerate}

Let us now return to the iterated monodromy groups associated to a quadratic polynomial over~$k$. By Theorem \ref{0Semirigid}, after conjugating everything by an automorphism of $T$ we may assume that $G^\geom=G$. Then $G^\arith$ is contained in the normalizer~$N$, and describing it is equivalent to describing the subgroup $G^\arith/G$ of $N/G$. In fact, we give the composite homomorphism $\bar\rho\colon \Gal(\bar k/k) \onto G^\arith/G \into N/G$ explicitly in terms of the description of $N/G$ mentioned above.

For this we look at the restriction of $\rho$ to the decomposition group at a point in~$P$. The image of the corresponding inertia group is always cyclic or pro-cyclic, and the decomposition group acts on it through the cyclotomic character. Using this we can deduce that $\bar\rho$ always factors through the cyclotomic character. 
Moreover, the generator $b_\infty$ of the image of the inertia group at $\infty$ is an odometer, so our previously assembled results on odometers in $G$ easily determine $\bar\rho$ completely.

Namely, in the periodic case the homomorphism $\bar\rho\colon\Gal(\bar k/k) \to N/G \cong (\BZ_2^\times)^r$ is the cyclotomic character $\Gal(\bar k/k)\to \BZ_2^\times$ followed by the diagonal embedding: see Theorem \ref{28GarithThm}.
In the strictly pre-periodic case with $r=2$ the homomorphism $\bar\rho\colon\Gal(\bar k/k) \to N/G \cong \BZ_2^\times/\{\pm1\}$ is the cyclotomic character followed by the projection $\BZ_2^\times\onto\BZ_2^\times/\{\pm1\}$.
In these two cases, therefore, the index $[G^\arith:G^\geom]$ is infinite if the ground field $k$ is finitely generated over its prime field. 

In the strictly pre-periodic case with $r\ge3$ the homomorphism $\bar\rho\colon\Gal(\bar k/k) \to N/G \cong \BF_2^{\raise2pt\hbox{$\kern1pt\scriptstyle\infty$}}$  is the cyclotomic character followed by the projection $\BZ_2^\times \onto (\BZ/8\BZ)^\times$ and a homomorphism $(\BZ/8\BZ)^\times \to \BF_2^{\raise2pt\hbox{$\kern1pt\scriptstyle\infty$}}$  that is given precisely in Theorem \ref{38GarithThm1}. In this case the group $G^\arith/G$ is elementary abelian of order dividing $4$, and in fact dividing $2$ except in one particular case: see Corollary \ref{38GarithThm2}. Thus $G^\arith$ depends only very mildly on the field~$k$. Moreover, in this case we have $G^\arith=G$ whenever $k$ contains the eighth roots of unity, in particular whenever $k$ is finite and its order is a square. It is conceivable that this state of affairs has some influence on the properties of the images of Frobenius elements.

The results on $G^\arith$ also yield upper bounds for the images of Galois on the set of preimages $\coprod_{n\ge0}f^{-n}(x')$ for any given point $x'\not\in P$, using Proposition \ref{18SpecDecomp}. In combination with the results on $G^\arith$ this implies a special case of a conjecture of Rafe Jones \cite[Conj. 3.2]{JonesSurvey2013}.

\medskip
Let us now point out the relation with existing results in the literature, without any claim of completeness. Although most results concern discrete iterated monodromy groups, some of them have a direct relation with the profinite iterated monodromy groups arising as closures of the discrete ones.
First, for generalities on self-similar groups acting on regular rooted trees see for instance Nekrashevych \cite{Nekrashevych-2003} or Grigorchuk-Savchuk-\v{S}uni\'{c} \cite{GrigorchukSavchukSunic2007}.

Next, much is already known in the following small cases: In the periodic case with $r=1$ the group $G$ is pro-cyclic. 
In the periodic case with $r=2$ the group $G$ is the closure of the so-called `Basilica' group, studied for instance in Grigorchuk-$\dot{\rm Z}$uk \cite{GrigorchukZuk2002}. 
In the periodic case with $r=3$ corresponding to both the `Douady rabbit' and the `airplane' some of our results follow from Nekrashevych \cite[\S8]{Nekrashevych-2007}.
In the strictly pre-periodic case with $s=1$ and $r=2$ the group $G$ is infinite pro-dihedral. 
In the strictly pre-periodic case with $s=1$ and $r=3$ the group $G$ is the closure of the Grigorchuk group, introduced in Grigorchuk \cite{Grigorchuk1980} and studied extensively, among others in Grigorchuk \cite[\S\S12-15]{Grigorchuk2000b}, Bartholdi-Grigorchuk \cite{Bartholdi-Grigorchuk-2002}, Grigorchuk-Sidki \cite{GrigorchukSidki2004}, Nekrashevych \cite{Nekrashevych-2007}.

An important result related to ours is contained in the paper \cite{Bartholdi-Nekrashevych-2008} by Bartholdi and Nekrashevych. There the authors describe the discrete iterated monodromy group of an arbitrary quadratic polynomial over $\BC$ in terms of generators satisfying explicit recursion relations similar to (\ref{0GenRecEq}). The only difference is that some of the terms $(a_{i-1},1)$ may be replaced by $(1,a_{i-1})$, and $(a_s,a_r)$ may be replaced by $(a_r,a_s)$. This goes a large part of the way towards our result that the closure $G^\geom$ is conjugate to our group $G$ in all cases.

Some arithmetic aspects of iterated monodromy groups associated to quadratic polynomials are studied in Boston-Jones \cite{BostonJones2007}, \cite{BostonJones2009}, \cite{JonesBoston2012}. Also, for a recent survey on arboreal Galois representations of global fields see Jones \cite{JonesSurvey2013}. The author hopes that the results of the present article might shed more light in these directions. 

\medskip
Here are several open questions that one might pursue next. 
Within the scope of the present paper, can some of the many case distinctions be avoided using a different setup?
It would also be interesting to obtain more results about conjugacy of generators and odometers in the strictly pre-periodic case, like an improvement of Proposition \ref{33WConjInGIsGConj} and analogues of Theorem \ref{28OdosConj} (b) and Theorem \ref{28OdoSemiRigid}. Again, such results might show the way towards a cleaner overall structure of the material.

Beyond that, a very natural problem is to generalize the results of this article to iterated monodromy groups of rational quadratic morphisms instead of quadratic polynomials.
The case of quadratic morphisms with infinite postcritical orbits is treated in a sequel to the present article by the same author \cite{Pink2013c}, using the same methods. 
By contrast the classification of finite postcritical orbits of arbitrary quadratic morphisms involves tuples of four integers, as opposed to two integers for quadratic polynomials; hence one can expect the complexity for those to increase roughly by a factor of $2$ on a logarithmic scale. 
\smiley

Finally there is the wide open field of applying the present results to the arithmetic properties of iterated quadratic polynomials, especially over a number field or a finite field. A starting point for this might be a study of Frobenius elements in $G^\arith$.

\medskip
\begin{samepage}
At last, the author would like to thank Volodia Nekrashevych for introducing him to important techniques in the area.

\end{samepage}

\numberwithin{Thm}{subsection}
\def\UseTheoremCounterForNextEquation{\setcounter{equation}{\value{Thm}}\addtocounter{Thm}{1}}
\renewcommand{\theequation}
             {\arabic{section}.\arabic{subsection}.\arabic{Thm}}


%
%

\newpage
\section{Generalities}
\label{1Generalities}


\subsection{Basics}
\label{11Basics}

Let $T$ be the infinite tree whose vertices are the finite words over the alphabet $\{0,1\}$ and where each vertex $t$ is connnected by an edge to the vertices $t0$ and $t1$. The empty word is called the \emph{root of~$T$}, making $T$ an infinite regular binary tree. For any integer $n\ge0$ we let $T_n$ denote the finite rooted subtree whose vertices are all words of length $\le n$. The words of length $n$ are precisely the vertices at distance $n$ from the root; we call the set of these the \emph{level $n$ of~$T$}. Thus the level $n$ consists of $2^n$ points.

\medskip
Unfortunately there do not seem to be established abbreviations for the automorphism groups of $T_n$ and~$T$. In this article we will write $W_n := \Aut(T_n)$ and $W := \Aut(T)$, based on the vague justification that $\Aut(T_n)$ is an iterated $\underline{\rm w}$reath product.

Namely, for any elements $u,v\in W$ let $(u,v)$ denote the element of $W$ defined by $t0\mapsto v(t)0$ and $t1\mapsto v(t)1$ for any word~$t$. This defines an embedding of $W\times W$ into $W$ which we identify with its image. 
Let $\sigma\in W$ denote the element of order $2$ defined by $t0\mapsto t1\mapsto t0$ for any word~$t$, and let $\langle\sigma\rangle$ be the subgroup of $W$ generated by it.
Then we can write $W$ as the semidirect product
\UseTheoremCounterForNextEquation
\begin{equation}\label{11WWreath}
W = (W\times W)\rtimes\langle\sigma\rangle.
\end{equation}
For any integer $n\ge1$ the same definitions define a semidirect product decomposition
\UseTheoremCounterForNextEquation
\begin{equation}\label{11WnWreath}
W_n = (W_{n-1}\times W_{n-1})\rtimes\langle\sigma\rangle.
\end{equation}
Since $W_0=1$, this describes $W_n$ as an iterated wreath product.
Calculating with this decomposition requires only the basic relations
\UseTheoremCounterForNextEquation
\begin{equation}\label{11BasicRels}
\left\{\!\begin{array}{rcl}
(u,v)\,(u',v') &\!\!=\!\!& (uu',vv') \\[3pt]
\sigma\,(u,v) &\!\!=\!\!& (v,u)\,\sigma \\[3pt]
\sigma^2 &\!\!=\!\!& 1
\end{array}\right\}
\end{equation}
for all $u,v,u',v'\in W$. 

\medskip
Every automorphism of $T$ fixes the root and thus stabilizes $T_n$ for every $n\ge0$. The restriction of automorphisms therefore induces a natural homomorphism $W\to W_n$, ${w\mapsto w|_{T_n}}$. This homomorphism is surjective and induces an isomorphism $W\cong\invlim_n W_n$.
Since $W_n$ is finite, this description realizes $W$ as a profinite topological group.

\medskip
All subgroups $G\subset W$ that we study in this article are closed and hence again profinite. Throughout we let $G_n$ denote the image of $G$ in~$W_n$. Since $G$ is closed, we obtain a natural isomorphism $G\cong\invlim_n G_n$.

\medskip
For any subset $S$ of a group we write $\langle S\rangle$ for the subgroup generated by~$S$ and abbreviate $\langle w_1,\ldots,w_r\rangle := \langle\{w_1,\ldots,w_r\}\rangle$. If the ambient group is a topological group, the closures of these subgroups are denoted $\Langle S\Rangle$ and $\Langle w_1,\ldots,w_r\Rangle$.


\subsection{Size}
\label{12Size}

For any integer $n\ge1$ the decomposition (\ref{11WnWreath}) implies that $|W_n| = 2\cdot|W_{n-1}|^2$. Since $|W_0|=1$, by induction on $n$ this easily implies that
\UseTheoremCounterForNextEquation
\begin{equation}\label{12WnOrder}
\qquad\qquad |W_n|\ =\ 2^{2^n-1} \qquad\qquad \hbox{for all $n\ge0$.}
\end{equation}
In particular $W_n$ is a finite $2$-group, and so $W=\invlim_n W_n$ is a pro-$2$-group. 

\medskip
The size of a closed subgroup $G\subset W$ is measured by its \emph{Hausdorff dimension}, which is defined as
\UseTheoremCounterForNextEquation
\begin{equation}\label{12Hausdorff}
\lim_{n\to\infty} \;\frac{\log_2|G_n|}{\log_2|W_n|}
\ =\ \lim_{n\to\infty} \frac{\log_2|G_n|}{2^n-1}
\end{equation}
if this limit exists. Its calculation usually requires some kind of recursive description of~$G$.


\subsection{Conjugacy}
\label{13Conjugacy}

For elements $w$, $w'\in W$ we write $w\sim w'$ if and only if $w$ and $w'$ are conjugate under~$W$.

\begin{Lem}\label{13ConjEquiv}
\begin{enumerate}
\item[(a)] For any elements $w$, $w'\in W$ we have $w\sim w'$ if and only if
$$\left\{\begin{array}{lclcccl}
w=(u,v) & \!\!\hbox{and}\! & 
w'=(u',v') & \!\!\hbox{and}\! & 
u\sim u' & \!\!\!\hbox{and}\! & v\sim v', \hbox{\ or} \\[3pt]
w=(u,v) & \!\!\hbox{and}\! & 
w'=(u',v') & \!\!\hbox{and}\! & 
u\sim v' & \!\!\!\hbox{and}\! & v\sim u', \hbox{\ or} \\[3pt]
w=(u,v)\,\sigma & \!\!\hbox{and}\! & 
w'=(u',v')\,\sigma & \!\!\hbox{and}\! & 
uv\sim u'v'.\!\!\! && 
\end{array}\right\}$$
\item[(b)] In particular, for any $u,v\in W$ we have
$$(u,v)\,\sigma 
\sim (uv,1)\,\sigma 
\sim (vu,1)\,\sigma 
\sim (1,uv)\,\sigma 
\sim (1,vu)\,\sigma
\sim (v,u)\,\sigma.$$
\item[(c)] 
The same assertions hold for any $n\ge1$ and $w,w'\in W_n$ and $u,v,u',v,\in W_{n-1}$.
\end{enumerate}
\end{Lem}

\begin{Proof}
If $w\sim w'$, either both elements act trivially or both act non-trivially on~$T_1$. 

In the first case we have $w=(u,v)$ and $w'=(u',v')$ for some $u,v,u',v'\in W$. Then for any $x,y\in W$ we have $(x,y)\,(u,v)\,(x,y)^{-1} = (xux^{-1},yvy^{-1})$, and so $w$ and $w'$ are conjugate under the subgroup $W\times W$ if and only if $u\sim u'$ and $v\sim v'$. This and the fact that $\sigma\,(u,v)\,\sigma^{-1} = (v,u)$ implies that $w$ and $w'$ are conjugate under the coset $(W\times W)\,\sigma$ if and only if $u\sim v'$ and $v\sim u'$. Together this yields the first two lines in (a).

In the second case we have $w=(u,v)\,\sigma$ and $w'=(u',v')\,\sigma$ for some $u,v,u',v'\in W$. If $w$ and $w'$ are conjugate under~$W$, then so are $w^2 = (u,v)\,\sigma\,(u,v)\,\sigma = (uv,vu)$ and $w^{\prime2}=(u'v',v'u')$. By the first case this implies that $uv$ is conjugate to $u'v'$ or $v'u'$, but since anyway $u'v'\sim v'u'$, it follows that $uv\sim u'v'$. Conversely assume that $uv\sim u'v'$. Then $u'v'=xuvx^{-1}$ for some $x\in W$, and so 
\begin{eqnarray*}
(u,v)\,\sigma 
&\sim& (x,u^{\prime-1}xu)\,(u,v)\,\sigma\,(x,u^{\prime-1}xu)^{-1} \\
&=& (x,u^{\prime-1}xu)\,(u,v)\,(u^{-1}x^{-1}u',x^{-1})\,\sigma \\
&=& (u',u^{\prime-1}xuvx^{-1})\,\sigma \\
&=& (u',v')\,\sigma,
\end{eqnarray*}
hence $w\sim w'$. Together this yields the third line in (c). 
That in turn directly implies (b).
Also, the same arguments apply with $w,w'\in W_n$ and $u,v,u',v,\in W_{n-1}$, proving (c).
\end{Proof}

\medskip
For later use we include some general facts about conjugacy.

\begin{Lem}\label{13ConjGen}
Let $G$ be a pro-$p$ group for a prime $p$ and let $a_1,\ldots,a_r$ be elements such that $G = \Langle a_1,\ldots,a_r\Rangle$. For every $i$ let $b_i\in G$ be a conjugate of~$a_i$. Then $G = \Langle b_1,\ldots,b_r\Rangle$.
\end{Lem}

\begin{Proof}
Taking an inverse limit reduces this to the case that $G$ is a finite $p$-group. 
We must then show that the subgroup $H := \langle b_1,\ldots,b_r\rangle$ is equal to $G=\langle a_1,\ldots,a_r\rangle$. But for each $i$ the elements $a_i$ and $b_i$ have the same image in the maximal abelian factor group $G_\ab$ of~$G$; hence the composite homomorphism $H\into G\onto G_\ab$ is surjective. Since $G$ is nilpotent, it follows that $H=G$, as desired.
\end{Proof}

\begin{Lem}\label{13ConjLimit}
Let $G$ be a profinite group given as the inverse limit of a filtered system of finite groups $G_i$ such that the projection maps $G\to G_i$ are surjective. Then two elements $g,g'\in G$ are conjugate if and only if their images in $G_i$ are conjugate for all~$i$.
\end{Lem}

\begin{Proof}
The `only if' part is obvious. For the `if' part assume that the images of $g,g'$ in $G_i$ are conjugate for all~$i$. For each $i$ let $U_i$ denote the set of elements $u\in G$ for which the images of $ugu^{-1}$ and $g'$ in $G_i$ coincide. Since $U_i$ is a union of cosets under the kernel of $G\onto G_i$, which is an open subgroup of~$G$, it is a closed subset of~$G$. By assumption it is also non-empty. Moreover, for any transition morphism $G_j\onto G_i$ in the system we have $U_j\subset U_i$. By the compactness of $G$ it thus follows that the intersection $U$ of all $U_i$ is non-empty. Any element $u\in U$ then satisfies $ugu^{-1}=g'$, as desired.
\end{Proof}

\medskip
The following result (compare \cite[Thm.3.1]{GNS2001}) uses only the basic relations (\ref{11BasicRels}) and is a good warm up exercise for our later calculations.

\begin{Prop}\label{13ConjPowers}
For any $w\in W$ and any $k\in\BZ_2^\times$ the element $w^k$ is conjugate to~$w$.
\end{Prop}

\begin{Proof}
(Compare Remark \ref{16OdoNormExpl}.)
Fix $k\in\BZ_2^\times$. By Lemma \ref{13ConjLimit} it suffices to prove that for every $n\ge0$ we have $w^k|_{T_n}\sim w|T_n$ for all $w\in W$. This is trivial for $n=0$, so assume that $n>0$ and the assertion holds for $n-1$ and all $w\in W$. 

Consider an element $w\in W$. If $w=(u,v)$ for $u,v\in W$, by the induction hypothesis we have $u^k|_{T_{n-1}}\sim u|_{T_{n-1}}$ and $v^k|_{T_{n-1}}\sim v|_{T_{n-1}}$. By Lemma \ref{13ConjEquiv} this implies that 
$$w^k|_{T_n} = 
(u^k|_{T_{n-1}},v^k|_{T_{n-1}}) 
\sim (u|_{T_{n-1}},v|_{T_{n-1}}) 
= w|_{T_n},$$
as desired.
Otherwise we have $w=(u,v)\,\sigma$ for $u,v\in W_{n-1}$. Then $w\sim(uv,1)\,\sigma$ by Lemma \ref{13ConjEquiv}.
We calculate $((uv,1)\,\sigma)^2 = (uv,1)\,\sigma\,(uv,1)\,\sigma = (uv,uv)$, and writing $k=2\ell+1$ with $\ell\in\BZ_2$ we deduce that
$$((uv,1)\,\sigma)^k = 
\bigl(((uv,1)\,\sigma)^2\bigr)^\ell\,(uv,1)\,\sigma
= \bigl((uv,uv)\bigr)^\ell\,(uv,1)\,\sigma
= \bigl((uv)^{\ell+1},(uv)^\ell\bigr)\,\sigma.$$
By Lemma \ref{13ConjEquiv} this is conjugate to $((uv)^k,1)\,\sigma$;
hence $w^k \sim ((uv)^k,1)\,\sigma$.
But by the induction hypothesis we have $(uv)^k|_{T_{n-1}}\sim uv|_{T_{n-1}}$. Thus with Lemma \ref{13ConjEquiv} we find that
$$w^k|_{T_n} \sim ((uv)^k,1)\,\sigma|_{T_n}
=  ((uv)^k|_{T_{n-1}},1)\,\sigma
\sim (uv|_{T_{n-1}},1)\,\sigma
= (uv,1)\,\sigma|_{T_n}
\sim w|_{T_n},$$
as desired.
\end{Proof}


\subsection{Recursive construction of elements}
\label{15Recursive}

\begin{Prop}\label{15RecRelsProp}
For all indices $i$ in a set $I$ consider distinct new symbols~$a_i$, finite words $f_i$ and $g_i$ over the alphabet $\{a_j,a_j^{-1}\mid j\in I\} \sqcup W$, and numbers $\nu_i\in\{0,1\}$.  
Then there exist unique elements $a_i\in W$ for all $i\in I$ which satisfy the recursion relations 
$$a_i = (f_i,g_i)\,\sigma^{\nu_i}.$$
Here the right hand side is interpreted as an element of $W$ by substituting the actual elements $a_i$ for the associated symbols. Also, the empty word is permitted for $f_i$ and $g_i$ and represents the identity element $1\in W$, as usual. 
\end{Prop}

\begin{Proof}
The restrictions $a_i|_{T_0}$ are necessarily trivial, and if for some $n>0$ the restrictions $a_i|_{T_{n-1}}$ are already known for all~$i$, then the recursion relations uniquely describe the restrictions $a_i|_{T_n}$ for all~$i$. By induction we therefore obtain unique elements $a_i|_{T_n}$ for all $n\ge0$, which combine to give the desired elements $a_i\in W$.
\end{Proof}

\medskip
Often $f_i$ and $g_i$ are words over the alphabet $\{a_j,a_j^{-1}\mid j\in I\}$ only, in which case the recursion relations do not involve any previously known elements of~$W$. Nevertheless one can construct very interesting elements in this way, as in (\ref{16OdoRec}), (\ref{2RecRels}), and (\ref{3RecRels}). In general one can construct elements with special properties with respect to given elements of~$W$, as in Remark \ref{16OdoNormExpl}, Proposition \ref{25Explicit}, and in (\ref{37W0Def}).

\medskip
One should be aware that different recursion relations can determine the same elements of~$W$. In particular:

\begin{Prop}\label{15RecTriv}
In Proposition \ref{15RecRelsProp} assume that for each $i\in I$ the expression $(f_i,g_i)\,\sigma^{\nu_i}$ has the form
$(h\kern1pt f' h^{-1},k\kern1pt g' \kern1pt k^{-1})$
for words $f'$, $g'$ over the alphabet $\{a_j,a_j^{-1}\mid j\in I\}$ only and arbitrary words $h$, $k$, which may all depend on~$i$.
Then $a_i=1$ for all $i\in I$. 
\end{Prop}

\begin{Proof}
Trivially $a_i|_{T_0}=1$ for all $i$, and if $a_i|_{T_{n-1}}=1$ for all~$i$, the recursion relations imply that 
$$a_i|_{T_n} 
= (h\kern1pt f' h^{-1},k\kern1pt g' \kern1pt k^{-1})|_{T_n}
\sim (f',g')|_{T_n}
= (f'|_{T_{n-1}},g'|_{T_{n-1}})
\sim (1,1)
= 1$$
and hence $a_i|_{T_n}=1$ for all $i$ as well. Induction on $n$ thus shows that $a_i=1$ on the whole tree~$T$.
\end{Proof}

\medskip
Also, elements satisfying certain kinds of recursion relations only up to conjugacy are conjugate to the elements satisfying the same recursion relations precisely:

\begin{Prop}\label{15RecConj=Conj}
Consider symbols $a_i$, words $f_i$, $g_i$, and numbers $\nu_i$ for all $i\in I$ as in Proposition \ref{15RecRelsProp}. Assume that for each $i\in I$ the expression $(f_i,g_i)\,\sigma^{\nu_i}$ has the form
$$\biggl\{\begin{array}{l}
(h\kern1pt a_j^\lambda h^{-1},k\kern1pt a_\ell^\mu \kern1pt k^{-1}) \quad \hbox{or} \\[3pt]
(h\kern1pt a_j^\lambda k^{-1},k\kern1pt a_j^\mu \kern1pt h^{-1})\,\sigma 
\end{array}$$
for words $h$, $k$, indices $j$, $\ell\in I$, and integers~$\lambda$,~$\mu$ which may depend on~$i$.
Consider other distinct new symbols $b_i$ and let $f'_i$ and $g'_i$ be the words over the alphabet ${\{b_j,b_j^{-1}\mid j\in I\}\sqcup W}$ obtained from $f_i$ and $g_i$ by replacing each occurrence of $a_j$ by~$b_j$. 
Let $a_i\in W$ be the elements satisfying $a_i = (f_i,g_i)\,\sigma^{\nu_i}$ which are furnished by Proposition \ref{15RecRelsProp}. Then for any elements $b_i\in W$ the following statements are equivalent:
\begin{enumerate}
\item[(a)] For every $i\in I$ the element $b_i$ is conjugate to $(f'_i,g'_i)\,\sigma^{\nu_i}$ under~$W$.
\item[(b)] For every $i\in I$ the element $b_i$ is conjugate to $a_i$ under~$W$.
\end{enumerate}
\end{Prop}

\begin{Proof}
Assume (b) and consider any $i\in I$. Suppose first that $(f_i,g_i)\,\sigma^{\nu_i}$ has the form $(h\kern1pt a_j^\lambda h^{-1},k\kern1pt a_\ell^\mu \kern1pt k^{-1})$. Then $(f'_i,g'_i)\,\sigma^{\nu_i}$ has the form $(h'\kern1pt b_j^\lambda h^{\prime-1},k'\kern1pt b_\ell^\mu \kern1pt k^{\prime-1})$ for some words $h'$,~$k'$. 
Moreover the assumption (b) implies that $b_i\sim a_i$ and $b_j^\lambda \sim a_j^\lambda $ and $b_\ell^\mu \sim a_\ell^\mu$, and hence
$$\begin{array}{rcccc}
b_i \ \sim\ a_i 
&\!\!\!\!=\!\!\!\!& (h\kern1pt a_j^\lambda h^{-1},k\kern1pt a_\ell^\mu \kern1pt k^{-1})
&\!\!\!\!\sim\!\!\!\!& (a_j^\lambda ,a_\ell^\mu) \\[2pt]
&&&&\wr \\[4pt]
(f'_i,g'_i)\,\sigma^{\nu_i} 
&\!\!\!\!=\!\!\!\!& (h'\kern1pt b_j^\lambda h^{\prime-1},k'\kern1pt b_\ell^\mu \kern1pt k^{\prime-1})
&\!\!\!\!\sim\!\!\!\!& (b_j^\lambda ,b_\ell^\mu) \rlap{.}
\end{array}$$
Now suppose that $(f_i,g_i)\,\sigma^{\nu_i}$ has the form $(h\kern1pt a_j^\lambda k^{-1},k\kern1pt a_j^\mu \kern1pt h^{-1})\,\sigma$. Then $(f'_i,g'_i)\,\sigma^{\nu_i}$ has the form $(h'\kern1pt b_j^\lambda k^{\prime-1},k'\kern1pt b_j^\mu \kern1pt h^{\prime-1})\,\sigma$ for some words $h'$,~$k'$. Moreover the assumption (b) implies that $b_i\sim a_i$ and $b_j^{\lambda+\mu}\sim a_j^{\lambda+\mu}$, hence using Lemma \ref{13ConjEquiv} we deduce that
$$\begin{array}{rcccccc}
b_i \ \sim\ a_i 
&\!\!\!\!=\!\!\!\!& (h\kern1pt a_j^\lambda k^{-1},k\kern1pt a_j^\mu \kern1pt h^{-1})\,\sigma
&\!\!\!\!\sim\!\!\!\!& (h\kern1pt a_j^{\lambda+\mu}\kern1pt h^{-1},1)\,\sigma
&\!\!\!\sim\!\!\!& (a_j^{\lambda+\mu},1)\,\sigma \\[2pt]
&&&&&&\wr \\[4pt]
(f'_i,g'_i)\,\sigma^{\nu_i} 
&\!\!\!\!=\!\!\!\!& (h'\kern1pt b_j^\lambda k^{\prime-1},k'\kern1pt b_j^\mu \kern1pt h^{\prime-1})\,\sigma
&\!\!\!\!\sim\!\!\!\!& (h'\kern1pt b_j^{\lambda+\mu}h^{\prime-1},1)\,\sigma
&\!\!\!\sim\!\!\!& (b_j^{\lambda+\mu},1)\,\sigma \rlap{.}
\end{array}$$
Together this shows that (b) implies (a).

Conversely assume (a). To deduce (b), by Lemma \ref{13ConjLimit} it suffices to show that $b_i|_{T_n}$ is conjugate to $a_i|_{T_n}$ in~$W_n$ for all $i$ and~$n$.
This is trivial for $n=0$, so assume that $n>0$ and that it holds for $n-1$ and all~$i$. Let $\sim$ denote conjugacy in~$W_n$, respectively in~$W_{n-1}$, so that the statements of Lemma \ref{13ConjEquiv} also hold in~$W_n$. 

Consider any $i\in I$. Suppose first that $(f_i,g_i)\,\sigma^{\nu_i}$ has the form $(h\kern1pt a_j^\lambda h^{-1},k\kern1pt a_\ell^\mu \kern1pt k^{-1})$. Then
$(f'_i,g'_i)\,\sigma^{\nu_i}$ has the form $(h'\kern1pt b_j^\lambda h^{\prime-1},k'\kern1pt b_\ell^\mu \kern1pt k^{\prime-1})$ for some words $h'$,~$k'$. 
Moreover the induction hypothesis implies that $b_j^\lambda|_{T_{n-1}}\sim a_j^\lambda|_{T_{n-1}}$ and $b_\ell^\mu |_{T_{n-1}}\sim a_\ell^\mu |_{T_{n-1}}$. Using the assumption (a) we therefore find that
$$\begin{array}{rcccccc}
b_i|_{T_n} \ \sim\ (f'_i,g'_i)\,\sigma^{\nu_i} |_{T_n}
&\!\!\!\!=\!\!\!\!& (h'\kern1pt b_j^\lambda h^{\prime-1},k'\kern1pt b_\ell^\mu \kern1pt k^{\prime-1}) |_{T_n}
&\!\!\!\!\sim\!\!\!\!& (b_j^\lambda ,b_\ell^\mu)|_{T_n}
&\!\!\!\!=\!\!\!\!& (b_j^\lambda |_{T_{n-1}},b_\ell^\mu|_{T_{n-1}}) \\[4pt]
&&&&&&\wr \\[2pt]
a_i|_{T_n}
&\!\!\!\!=\!\!\!\!& (h\kern1pt a_j^\lambda h^{-1},k\kern1pt a_\ell^\mu \kern1pt k^{-1})|_{T_n}
&\!\!\!\!\sim\!\!\!\!& (a_j^\lambda ,a_\ell^\mu)|_{T_n}
&\!\!\!\!=\!\!\!\!& (a_j^\lambda |_{T_{n-1}},a_\ell^\mu|_{T_{n-1}}) \rlap{.}
\end{array}$$
Now suppose that $(f_i,g_i)\,\sigma^{\nu_i}$ has the form $(h\kern1pt a_j^\lambda k^{-1},k\kern1pt a_j^\mu \kern1pt h^{-1})\,\sigma$. Then $(f'_i,g'_i)\,\sigma^{\nu_i}$ has the form $(h'\kern1pt b_j^\lambda k^{\prime-1},k'\kern1pt b_j^\mu \kern1pt h^{\prime-1})\,\sigma$ for some words $h'$,~$k'$. Moreover the induction hypothesis implies that $b_j^{\lambda+\mu}|_{T_{n-1}}\sim a_j^{\lambda+\mu}|_{T_{n-1}}$. Using the assumption (a) and Lemma \ref{13ConjEquiv} we thus find that
$$\begin{array}{rcccccc}
b_i|_{T_n} \ \sim\ (f'_i,g'_i)\,\sigma^{\nu_i} |_{T_n}
&\!\!\!\!=\!\!\!\!& (h'\kern1pt b_j^\lambda k^{\prime-1},k'\kern1pt b_j^\mu \kern1pt h^{\prime-1})\,\sigma |_{T_n}
&\!\!\!\!\sim\!\!\!\!& (b_j^{\lambda+\mu},1)\,\sigma |_{T_n}
&\!\!\!\!=\!\!\!\!& (b_j^{\lambda+\mu}|_{T_{n-1}},1)\,\sigma \\[4pt]
&&&&&&\wr \\[2pt]
a_i|_{T_n}
&\!\!\!\!=\!\!\!\!& (h\kern1pt a_j^\lambda k^{-1},k\kern1pt a_j^\mu \kern1pt h^{-1})\,\sigma |_{T_n}
&\!\!\!\!\sim\!\!\!\!& (a_j^{\lambda+\mu},1)\,\sigma |_{T_n}
&\!\!\!\!=\!\!\!\!& (a_j^{\lambda+\mu}|_{T_{n-1}},1)\,\sigma \rlap{.}
\end{array}$$
Together this shows the desired assertion for $n$ and all~$i$; and so it follows for all $n\ge0$ by induction. This proves that (a) implies (b).
\end{Proof}


\subsection{Signs}
\label{15Signs}

For any integer $n\ge0$, any element $w\in W$ fixes the root and thus permutes the level $n$ of~$T$. We let $\sgn_n(w)$ denote the sign of the induced permutation of the level~$n$, which defines a homomorphism 
\UseTheoremCounterForNextEquation
\begin{equation}\label{15SignDef}
\sgn_n\colon W \longto \{\pm1\}.
\end{equation}
For any $m\ge n$ this homomorphism factors through a homomorphism $W_m\to\{\pm1\}$, which we again denote by~$\sgn_n$.
For $n=0$ the homomorphism is trivial. For $n\ge1$ the definition of the embedding $W\times W\into W$ implies that
\UseTheoremCounterForNextEquation
\begin{equation}\label{15SignRec}
\sgn_n((w,w')) = \sgn_{n-1}(w)\cdot\sgn_{n-1}(w')
\end{equation}
for all $w,w'\in W$. Also, for $n\ge1$ the element $\sigma$ of order $2$ has no fixed points on level~$n$ and hence precisely $2^{n-1}$ orbits of length $2$. Therefore
\UseTheoremCounterForNextEquation
\begin{equation}\label{15SignSigma}
\sgn_n(\sigma) = 
\biggl\{\!\begin{array}{rl}
-1 & \hbox{if $n=1$,}\\[3pt]
1 & \hbox{if $n\not=1$.}
\end{array}
\end{equation}
Together these formulas can be used to calculate the signs of any recursively described elements of~$W$.

\begin{Ex}\label{15Ex}
\rm Consider the elements $b_1,b_2,\ldots\in W$ defined by the recursion relations
$$\biggl\{\begin{array}{ll}
b_1 = \sigma & \\[3pt]
b_i = (b_{i-1},1) & \hbox{for all $i>1$.}\\
\end{array}$$
Then the formula (\ref{15SignRec}) implies that $\sgn_n(b_i) = \sgn_{n-1}(b_{i-1})$ for all $n\ge1$ and $i>1$. By (\ref{15SignSigma}) and induction it follows that for all $n,i\ge1$ we have
$$\sgn_n(b_i) \ =\ 
\biggl\{\!\begin{array}{rl}
-1 & \hbox{if $n    =i$,}\\[3pt]
 1 & \hbox{if $n\not=i$.}
\end{array}$$
\end{Ex}

\begin{Prop}\label{15SignWnSurj}
For any $n\ge0$ the homomorphism
$$W_n\to \{\pm1\}^n,\ w\mapsto (\sgn_m(w))_{m=1}^n$$ 
induces an isomorphism from the maximal abelian factor group $W_{n,\ab}$ to $\{\pm1\}^n$.
\end{Prop}

\begin{Proof}
Example \ref{15Ex} implies that the homomorphism is surjective; hence it induces a surjective homomorphism $W_{n,\ab}\onto\{\pm1\}^n$. It remains to show that this is an isomorphism. For $n=0$ that is obvious, so assume that $n>0$ and that it is an isomorphism for $W_{n-1}$. 
Then the decomposition $W_n = (W_{n-1}\times W_{n-1})\rtimes\langle\sigma\rangle$ implies that $W_{n,\ab}$ is also the maximal abelian factor group of $(W_{n-1,\ab}\times W_{n-1,\ab})\rtimes\langle\sigma\rangle$. By the induction hypothesis this is isomorphic to the maximal abelian factor group of $(\{\pm1\}^{n-1}\times\{\pm1\}^{n-1})\rtimes\langle\sigma\rangle$, where $\sigma$ acts by interchanging the two factors $\{\pm1\}^{n-1}$. The commutators of $\sigma$ with $\{\pm1\}^{n-1}\times\{\pm1\}^{n-1}$ form the diagonally embedded subgroup $\diag(\{\pm1\}^{n-1})$; hence the maximal abelian factor group of $W_n$ is isomorphic to $\{\pm1\}^{n-1}\times\langle\sigma\rangle$ and hence to $\{\pm1\}^n$. For reasons of cardinality the surjective homomorphism $W_{n,\ab}\onto\{\pm1\}^n$ is therefore an isomorphism. By induction this now follows for all $n\ge0$.
\end{Proof}

\begin{Prop}\label{15SignGnSurj}
\begin{enumerate}
\item[(a)] For any subgroup $G_n\subset W_n$ we have $G_n=W_n$ if and only if the combined homomorphism $G_n\to \{\pm1\}^n$, $g\mapsto (\sgn_m(g))_{m=1}^n$ is surjective.
\item[(b)] For any closed subgroup $G\subset W$ we have $G=W$ if and only if the combined homomorphism $G\to \{\pm1\}^\infty$, $g\mapsto (\sgn_m(g))_{m=1}^\infty$ is surjective.
\end{enumerate}
\end{Prop}

\begin{Proof}
(This is essentially in Stoll \cite[p.~241 Thm.]{Stoll1992}.)
Since $W_n$ is a finite $2$-group, it is nilpotent, and hence a subgroup $G_n\subset W_n$ is equal to $W_n$ if and only if it surjects to the maximal abelian factor group of~$W_n$. Thus part (a) follows from Proposition \ref{15SignWnSurj}. Part (b) follows from (a) by taking inverse limits.
\end{Proof}


\subsection{Odometers}
\label{16Odometer}

As another useful preparation for the main part of the article (and for fun) consider the element $a\in W$ defined by the recursion relation
\UseTheoremCounterForNextEquation
\begin{equation}\label{16OdoRec}
a = (a,1)\,\sigma.
\end{equation}
It is called the \emph{standard odometer} \cite[Ex.1]{GrigorchukSunic2007} or \emph{adding machine} \cite[2.5.1]{Nekrashevych-2009}. The reason for this is that if we view each vertex on some level $n$ of~$T$, that is, each word of length $n$ over the alphabet $\{0,1\}$, as an integer in the interval $[0,2^n-1]$ written in binary notation, then $a$ acts on these by adding $1$ modulo $2^n$. This is precisely what odometers of real life vehicles do with integers to base $10$ in place of~$2$.

\medskip
More generally we will call any element of $W$ that is conjugate to $a$ an \emph{odometer}.
For a part of the following result compare \cite[Prop.3.3]{GrigorchukSunic2007}:

\begin{Prop}\label{16OdoEquiv}
For any $w\in W$ the following are equivalent:
\begin{enumerate}
\item[(a)] $w$ is an odometer, i.e., it is conjugate to $a$.
\item[(b)] $w$ is conjugate to $(w,1)\,\sigma$.
\item[(c)] $w$ acts transitively on level $n$ for all $n\ge1$.
\item[(d)] $\sgn_n(w)=-1$ for all $n\ge1$.
\end{enumerate}
\end{Prop}

\begin{Proof}
If $w\sim a$, with Lemma \ref{13ConjEquiv} we deduce that $w\sim a = (a,1)\,\sigma \sim (w,1)\,\sigma$, proving the implication (a)$\Rightarrow$(b).

Next assume that $w$ is conjugate to $(w,1)\,\sigma$. Then it clearly acts non-trivially and hence transitively on level~$1$. Consider $n\ge2$ such that $w$ acts transitively on level $n-1$. 
Since $(w,1)\,\sigma$ interchanges the two embedded half trees in~$T$, it interchanges 
their respective levels $n-1$. By assumption the square $((w,1)\,\sigma)^2 = (w,w)$ acts transitively on each of these levels $n-1$. Thus $(w,1)\,\sigma$ itself acts transitively on the their union, which is the level $n$ of~$T$. The same then follows for $w\sim (w,1)\,\sigma$, as desired. By induction on $n$ this proves the implication (b)$\Rightarrow$(c).


Next, if $w$ acts transitively on a level $n\ge1$, it acts by a cycle of length $2^n$ which has sign $-1$. This directly yields the implication (c)$\Rightarrow$(d).

To establish the remaining implication (d)$\Rightarrow$(a), by Lemma \ref{13ConjLimit} it suffices to prove: For any $n\ge0$ and any $w$ as in (d) we have $w|{T_n} \sim a|{T_n}$. This is trivial for $n=0$, so assume that $n>0$ and that it holds universally for $n-1$.
Consider any $w\in W$ satisfying (d). Then $\sgn_1(w)=-1$ shows that $w$ acts nontrivially on $T_1$; hence $w=(u,v)\,\sigma$ for some $u,v\in W$.
For all $n\ge1$ the formulas (\ref{15SignRec}) and (\ref{15SignSigma}) then imply that 
$$\sgn_n(uv)\ =\ \sgn_n(u)\cdot\sgn_n(v)\cdot\sgn_{n+1}(\sigma)\ =\ \sgn_{n+1}(w)\ =\ -1,$$
thus $uv$ also satisfies the condition in (d). Using Lemma \ref{13ConjEquiv} and the induction hypothesis for $uv$ it follows that
$$\begin{array}{rcccccc}
w|_{T_n} &\!\!\!\!=\!\!\!\!& (u,v)\,\sigma|_{T_n}
&\!\!\!\!\sim\!\!\!\!& (uv,1)\,\sigma|_{T_n}
&\!\!\!\!=\!\!\!\!& (uv|_{T_{n-1}},1)\sigma \\[4pt]
&&&&&&\wr \\[2pt]
&& a|_{T_n}
&\!\!\!\!=\!\!\!\!& (a,1)\,\sigma|_{T_n}
&\!\!\!\!=\!\!\!\!& (a|_{T_{n-1}},1)\sigma
\end{array}$$
and so the desired assertion holds for $n$. By induction it follows for all $n\ge0$, proving the implication (d)$\Rightarrow$(a).
\end{Proof}

%
\begin{Prop}\label{16OdoProp}
For any odometer $w\in W$ we have:
\begin{enumerate}
\item[(a)] The subgroup $\Langle w\Rangle$ of $W$ is isomorphic to~$\BZ_2$.
\item[(b)] It is its own centralizer.
\item[(c)] Its normalizer is isomorphic to $\BZ_2\rtimes\BZ_2^\times$.
\end{enumerate}
\end{Prop}

\begin{Proof}
The implication (a)$\Rightarrow$(c) in Proposition \ref{16OdoEquiv} shows that $w$ acts transitively on each level $n$ through a cycle of length $2^n$. For $n\to\infty$ this shows that $w$ has infinite order, which implies (a).
Assertion (b) results from the fact that in any symmetric group on $m$ letters, any $m$-cycle generates its own centralizer. For (c) observe that $\Aut(\BZ_2) = \BZ_2^\times$. Moreover, for any $k\in \BZ_2^\times$ the element $w^k\in W$ still has the same properties (c) or (d) from Proposition \ref{16OdoEquiv}, and is therefore conjugate to~$w$. This shows that the normalizer of $\Langle w\Rangle$ in $W$ surjects to $\Aut(\Langle w\Rangle) = \BZ_2^\times$. Since the centralizer is $\Langle w\Rangle\cong\BZ_2$, this implies (c).
\end{Proof}

\medskip
For the standard odometer we can make the isomorphism in Proposition \ref{16OdoProp} explicit as follows:

\begin{Prop}\label{16OdoNormExpl}
For any $k\in\BZ_2^\times$ consider the element $z_k\in W$ defined using Proposition \ref{15RecRelsProp} by the recursion relation
$$z_k = (z_k,\smash{a^{\frac{k-1}{2}}}z_k).$$
\begin{enumerate}
\item[(a)] For any $k\in\BZ_2^\times$ we have $z_kaz_k^{-1} = a^k$.
\item[(b)] For any $k,k'\in\BZ_2^\times$ we have $z_kz_{k'}=z_{kk'}$.
\item[(c)] There is an isomorphism $\BZ_2\rtimes\BZ_2^\times \stackrel{\sim}{\longto} \Norm_W(\Langle a\Rangle)$, $(i,k)\mapsto a^iz_k$.
\end{enumerate}
\end{Prop}

\begin{Proof}
To show (a) we abbreviate $\ell := \frac{k-1}{2} \in \BZ_2$ and $z := z_k = (z,a^\ell z)$. Since $a^2= (a,1)\,\sigma\,(a,1)\,\sigma = (a,a)$, we have $a^{-k} = a^{-1-2\ell} = a^{-1}\,(a^{-\ell},a^{-\ell})$. Therefore
\begin{eqnarray*}
t \ :=\ z a z^{-1} a^{-k}
&\!\!=\!\!& (z,a^\ell z)\, (a,1)\,\sigma\, (z^{-1},z^{-1}a^{-\ell})\, 
\sigma^{-1}\, (a^{-1},1)\, (a^{-\ell},a^{-\ell}) \\
&\!\!=\!\!& \bigl(zaz^{-1}a^{-\ell}a^{-1}a^{-\ell}\,,\,a^\ell zz^{-1}a^{-\ell}\bigr) \\
&\!\!=\!\!& (zaz^{-1}a^{-k},1) \\
&\!\!=\!\!& (t,1).
\end{eqnarray*}
By Proposition \ref{15RecTriv} this implies that $t=1$, proving (a).
To show (b) we abbreviate in addition $\ell' := \frac{k'-1}{2}$ and $z' := z_{k'} = (z',a^{\ell'} z')$ as well as $k'' := kk'$ and $\ell'' := \frac{k''-1}{2} = \ell+k\ell'$ and $z'' := z_{k''} = (z'',a^{\ell''} z'')$. Then
\begin{eqnarray*}
u \ :=\ z z' z^{\prime\prime-1}
&=& (z,a^\ell z)\, (z',a^{\ell'} z')\, (z^{\prime\prime-1},z^{\prime\prime-1}a^{-\ell''}) \\
&=& (zz'z^{\prime\prime-1} , a^\ell z a^{\ell'} z' z^{\prime\prime-1}a^{-\ell''}) \\
&=& (u , a^\ell z a^{\ell'} z^{-1} u a^{-\ell''}) \\
&\stackrel{\smash{(a)}}{=}& (u , a^\ell a^{k\ell'} u a^{-\ell''}) \\
&=& (u , a^{\ell''} u a^{-\ell''}).
\end{eqnarray*}
Again by Proposition \ref{15RecTriv} this implies that $u=1$, proving (b).
Finally, (c) is a direct consequence of (a) and (b).
\end{Proof}


\subsection{Iterated monodromy groups}
\label{17Monodromy}

Let $k$ be a field of characteristic different from~$2$, and let $\bar k$ be a separable closure of~$k$. Let $\BP^1_k$ denote the projective line over~$k$. Let $f\colon \BP^1_k\to\BP^1_k$ be a morphism over~$k$ of degree~$2$, called a \emph{quadratic morphism} for short. 
Let $C\subset\BP^1(\bar k)$ denote the set of critical points of~$f$, which by Hurwitz
has cardinality~$2$. For any integer $n\ge0$ let $f^n\colon\BP^1_k\to\BP^1_k$ denote the \emph{$n^{th}$ iterate} of~$f$, defined by setting $f^0:=\id$ and $f^{n+1} := f\circ f^n$.
Then the \emph{(strictly) postcritical orbit of $f$} is the set
\UseTheoremCounterForNextEquation
\begin{equation}\label{17PostCritDef}
\textstyle P\ :=\ \bigcup_{n\ge1}f^n(C) \ \subset\ \BP^1(\bar k).
\end{equation}
This is precisely the set of points in $\BP^1(\bar k)$ over which some iterate $f^n$ is ramified. Note that $C$ and $P$ may or may not be contained in $\BP^1(k)$.
We can view $P$ as the set of vertices of a directed graph with an edge from $p$ to $f(p)$ for every $p\in P$, and with a subset $f(C)$ of two distinct specially marked `entry points'. For a classification of the possibilities for this graph see \cite[Class.\ 2.3]{Pink2013}, respectively Cases \ref{17Cases} below in the polynomial case.

\medskip
Fix a separably closed overfield $L$ of $\bar k$ and a point $x_0\in\BP^1(L)\setminus P$. For instance, any point $x_0\in \BP^1(L)$ that is not algebraic over $k$ is in order.
For any $n\ge0$ let $f^{-n}(x_0)$ denote the set of $2^n$ points $x_n\in\BP^1(L)$ with $f^n(x_n)=x_0$. Let $T_{x_0}$ denote the infinite 
graph whose set of vertices is the exterior disjoint union $\coprod_{n\ge0}f^{-n}(x_0)$, where any vertex $x_n\in f^{-n}(x_0)$ for $n>0$ is connected by an edge towards $f(x_n)\in f^{-(n-1)}(x_0)$. By construction this is a regular rooted binary tree with root $x_0\in f^{-0}(x_0)$ and level $n$ set $f^{-n}(x_0)$ for every~$n$.

For any $n\ge0$ the morphism $f^n\colon\BP^1_k\to\BP^1_k$ defines a connected unramified covering of $\BP^1_k \setminus P_n$ for the finite set $P_n := \bigcup_{m=1}^nf^m(C)$.
Up to isomorphism this covering is determined by the associated monodromy action on $f^{-n}(x_0)$ of the \'etale fundamental group $\pi_1^\et(\BP^1_k \setminus P_n,x_0)$. Varying $n$, these actions combine to a natural continuous homomorphism
\UseTheoremCounterForNextEquation
\begin{equation}\label{17MonoRep}
\rho\colon\ \pi_1^\et(\BP^1_k \setminus P,x_0)
:= \smash{\invlim_n}\ \pi_1^\et(\BP^1_k \setminus P_n,x_0)
 \longto \Aut(T_{x_0}),
\end{equation}
which encodes the Galois theoretic properties of the whole tower of ramified coverings~$f^n$.
Recall the natural short exact sequence
\UseTheoremCounterForNextEquation
\begin{equation}\label{17ShortExact}
1 \longto \pi_1^\et(\BP^1_{\bar k} \setminus P,x_0)
\longto \pi_1^\et(\BP^1_k \setminus P,x_0)
\longto \Gal(\bar k/k) \longto 1.
\end{equation}
We are interested in the images 
\UseTheoremCounterForNextEquation
\begin{eqnarray}
G^\geom  &:=& \rho\bigl(\pi_1^\et(\BP^1_{\bar k} \setminus P,x_0)\bigr) 
\ \subset\Aut(T_{x_0}), \\
\UseTheoremCounterForNextEquation
G^\arith &:=& \rho\bigl(\pi_1^\et(\BP^1_{     k} \setminus P,x_0)\bigr) 
\ \subset\Aut(T_{x_0}), 
\end{eqnarray}
which are called the \emph{geometric}, respectively \emph{arithmetic, iterated monodromy group} associated to~$f$. By construction $G^\geom$ is a closed normal subgroup of $G^\arith$, and the short exact sequence (\ref{17ShortExact}) induces a surjective homomorphism
\UseTheoremCounterForNextEquation
\begin{equation}\label{17FactorMonoRep}
\Gal(\bar k/k) \onto G^\arith/G^\geom.
\end{equation}

These groups are independent of the choice of $L$ and $x_0$ in the sense that they do not change when $L$ is replaced by any separably closed overfield, and for a different choice of $x_0$ there exists an isomorphism between the two resulting trees which is equivariant for the monodromy action. Thus if we identify $T_{x_0}$ with any other regular rooted binary tree, the resulting monodromy groups will be unique up to conjugacy there. Throughout the following we will identify $T_{x_0}$ with the standard tree $T$ from Subsection \ref{11Basics}.

The group $G^\geom$ also does not change under extending~$k$. Thus in principle it can be calculated over any field of definition of the quadratic polynomial~$f$. In particular, if $k$ has characteristic zero, then $G^\geom$ can be calculated over~$\BC$. It is then the closure of the image of the topological fundamental group $\pi_1(\BP^1(\BC) \setminus P,x_0)$. The latter can be described by explicit generators as a self-similar subgroup of $\Aut(T)$ as, say, in Nekrashevych \cite{Nekrashevych-2003}, \cite[Ch.5]{Nekrashevych-2005}, \cite{Nekrashevych-2009}, Bartholdi-Nekrashevych \cite{Bartholdi-Nekrashevych-2008}, Dau \cite{Dau}, Grigorchuk et al.~\cite{GrigorchukSavchukSunic2007}. Moreover, in \cite{Pink2013} we showed how this in turn can be used to describe $G^\geom$ in positive characteristic as well by lifting the given quadratic polynomial to characteristic zero.

\medskip
In the present article, however, we are only interested in the profinite monodromy groups, and we will determine them using only the self-similarity of the monodromy representation and the combinatorics of the postcritical orbit, as follows.

\medskip
Let $x_1$ and $x_1'$ denote the two points in $f^{-1}(x_0)$. Then $f^{-(n+1)}(x_0) = f^{-n}(x_1) \sqcup f^{-n}(x_1')$ for all $n\ge0$, so after deleting the root $x_0$ the tree $T_{x_0}$ decomposes into the two regular rooted binary trees $T_{x_1}$ and $T_{x_1'}$ with roots $x_1$ and~$x_1'$, respectively. Abbreviate $\tilde P := f^{-1}(P) \allowbreak \subset \BP^1(\bar k)$. Then by functoriality $f\colon\BP^1_k\to\BP^1_k$ induces a homomorphism 
\UseTheoremCounterForNextEquation
\begin{equation}\label{17FStar}
f_*\colon\ \pi_1^\et\bigl({\BP^1_k \setminus \tilde P},x_1\bigr) \ \longto\ 
\pi_1^\et({\BP^1_k \setminus P},x_0).
\end{equation}
Its image is the subgroup of index $2$ of $\pi_1^\et({\BP^1_k \setminus P},x_0)$ which fixes $x_1\in f^{-1}(x_0)$. It therefore fixes both $x_1$ and $x_1'$ and 
acts on each of the subtrees $T_{x_1}$ and $T_{x_1'}$. 

To describe the action on these subtrees note that the inclusion $f(P)\subset P$ implies that $P\subset \tilde P$ and hence $\BP^1_k \setminus \tilde P \subset \BP^1_k \setminus P$. Thus there is a natural surjective homomorphism 
\UseTheoremCounterForNextEquation
\begin{equation}\label{17IdStar}
\id_*\colon\ \pi_1^\et\bigl({\BP^1_k \setminus \tilde P},x_1\bigr) \longonto \pi_1^\et({\BP^1_k \setminus P},x_1).
\end{equation}
The action of $\pi_1^\et\bigl({\BP^1_k \setminus \tilde P},x_1\bigr)$ on $T_{x_1}$ through $f_*$ is the same as that obtained by composing (\ref{17IdStar}) with the natural action of $\pi_1^\et({\BP^1_k \setminus P},x_1)$ on~$T_{x_1}$. 
On the other hand let $\tau$ denote the non-trivial covering automorphism of $f\colon\BP^1_k\to\BP^1_k$, so that $x_1'=\tau(x_1)$. Then $\tau$ induces a natural isomorphism
\UseTheoremCounterForNextEquation
\begin{equation}\label{17TauStar}
\tau_*\colon\ \pi_1^\et({\BP^1_k \setminus \tilde P},x_1) 
\ \stackrel{\sim}{\longto}\ \pi_1^\et({\BP^1_k \setminus \tilde P},x_1').
\end{equation}
The action of $\pi_1^\et\bigl({\BP^1_k \setminus \tilde P},x_1\bigr)$ on $T_{x_1'}$ through $f_*$ is the same as that obtained by composing (\ref{17TauStar}) with the surjection (\ref{17IdStar}) for $x_1'$ in place of~$x_1$ and the natural action of $\pi_1^\et({\BP^1_k \setminus P},x_1')$ on~$T_{x_1'}$. 

Now recall that except for the change of base point, the actions of $\pi_1^\et({\BP^1_k \setminus P},x)$ on~$T_x$ for $x=x_0,x_1,x_1'$ all describe the same inverse system of ramified coverings $f^n\colon\BP^1_k\to \nobreak \BP^1_k$. Choose isomorphisms $T_{x_1} \cong T_{x_0}\cong T_{x_1'}$ and corresponding identifications of the respective fundamental groups (which over $\BC$ usually involves choosing paths from $x_0$ to $x_1$ and~$x_1'$). Then the above statements say that the monodromy actions $\rho\circ f_*$ and $({\rho\circ\id_*},{\rho\circ\id_*\circ\tau_*})$ of $\pi_1^\et({\BP^1_k \setminus \tilde P},x_0)$ on $T_{x_0}\setminus\{x_0\} = T_{x_1}\sqcup T_{x_1'}$ coincide. In particular, with the corresponding identifications $T_{x_0}\cong T$ and $\Aut(T_{x_0})\cong W$ we obtain a commutative diagram:
\UseTheoremCounterForNextEquation
\begin{equation}\label{17RecRepDiag}
\vcenter{\xymatrix@R+10pt{
\pi_1^\et({\BP^1_k \setminus P},x_0) \ar[r]^-{\rho} &
\Aut(T_{x_0}) \ar@{=}[r]^-\sim & 
W \\
\pi_1^\et({\BP^1_k \setminus \tilde P},x_1) 
\ar[u]^-{f_*} \ar[r] \ar[dr] \ar@{=}[d]_-{\wr} &
\Aut(T_{x_1}) \times \Aut(T_{x_1'}) \ar@{^{ (}->}[u] \ar@{->>}[d]_-{\proj_1} &
W\times W \ar@{^{ (}->}[u] \ar@{->>}[d]_-{\proj_1} \ar@{=}[l]_-{\sim} \\
\pi_1^\et({\BP^1_k \setminus \tilde P},x_0) \ar@{->>}[d]_-{\id_*} \ar[dr] &
\Aut(T_{x_1}) \ar@{=}[d]_-{\wr} & 
W \ar@{=}[l]_-{\sim} \\
\pi_1^\et({\BP^1_k \setminus P},x_0) \ar[r]^{\rho} &
\Aut(T_{x_0}) & 
W \ar@{=}[l]_-{\sim} \ar@{=}[u] \\}}
\end{equation}
This shows that the monodromy representation $\rho$ is \emph{fractal} or \emph{self-similar} in that a copy of $\rho$ occurs in a proper subquotient. The analysis of $\rho$ and of its geometric and arithmetic monodromy groups therefore lends itself to recursive arguments.

By the construction of $G^\arith$ the diagram (\ref{17RecRepDiag}) induces a commutative diagram:
\UseTheoremCounterForNextEquation
\begin{equation}\label{17RecArithDiag}
\vcenter{\xymatrix@R+10pt{
\pi_1^\et({\BP^1_k \setminus P},x_0) \ar@{->>}[r]^-{\rho} &
G^\arith \ar@{^{ (}->}[r] & 
W \\
\pi_1^\et({\BP^1_k \setminus \tilde P},x_1) 
\ar[u]^-{f_*} \ar@{->>}[r] \ar@{->>}[d]_-{\id_*} &
G^\arith\cap(W\times W) \ar@{^{ (}->}[u] \ar@{->>}[d]_-{\proj_1} \ar@{^{ (}->}[r] &
W\times W \ar@{^{ (}->}[u] \ar@{->>}[d]_-{\proj_1} \\
\pi_1^\et({\BP^1_k \setminus P},x_0) \ar@{->>}[r]^-{\rho} &
G^\arith \ar@{^{ (}->}[r] & 
W \\}}
\end{equation}
Here the lower middle vertical arrow is surjective, because the others in the lower left square are surjective. Thus $G^\arith$ is \emph{self-similar} in that it is isomorphic to a proper subquotient of itself. With $\bar k$ in place of~$k$, the same remarks apply to $G^\geom$ and the commutative diagram:
\UseTheoremCounterForNextEquation
\begin{equation}\label{17RecGeomDiag}
\vcenter{\xymatrix@R+10pt{
\pi_1^\et({\BP^1_{\bar k} \setminus P},x_0) \ar@{->>}[r]^-{\rho} &
G^\geom \ar@{^{ (}->}[r] & 
W \\
\pi_1^\et({\BP^1_{\bar k} \setminus \tilde P},x_1) 
\ar[u]^-{f_*} \ar@{->>}[r] \ar@{->>}[d]_-{\id_*} &
G^\geom\cap(W\times W) \ar@{^{ (}->}[u] \ar@{->>}[d]_-{\proj_1} \ar@{^{ (}->}[r] &
W\times W \ar@{^{ (}->}[u] \ar@{->>}[d]_-{\proj_1} \\
\pi_1^\et({\BP^1_{\bar k} \setminus P},x_0) \ar@{->>}[r]^-{\rho} &
G^\geom \ar@{^{ (}->}[r] & 
W \\}}
\end{equation}
Note that the groups $G^\arith$ and $G^\geom$ and several of the maps in these diagrams depend on the chosen identifications of trees and are therefore unique only up to conjugacy by~$W$. 

\medskip
Now we can describe generators of~$G^\geom$.
Suppose first that $P$ is finite, and write $P=\{p_1,\ldots,p_m\}$ for distinct elements~$p_i$. Let $\Delta$ denote the maximal pro-$2$ factor group of ${\pi_1^\et(\BP^1_k \setminus P,x_0)}$. Since the characteristic of $k$ is different from~$2$, Grothendieck's description of the tame fundamental group \cite[exp. XIII, Cor. 2.12]{SGA1} says that for every $p_i\in P$ there exists a generator $\delta_{p_i}$ of the pro-$2$ inertia group at some geometric point above~$p_i$, such that $\Delta$ is the quotient of the free pro-$2$ group on the elements $\delta_{p_1},\ldots,\delta_{p_m}$ by a single relation of the form $\delta_{p_1}\!\cdots\delta_{p_m}=1$. Setting $b_{p_i} := \rho(\delta_{p_i})\in G^\geom$, we deduce that
\UseTheoremCounterForNextEquation
\begin{equation}\label{17GgeomGensFin}
G^\geom = \Langle b_{p_1},\ldots,b_{p_m}\Rangle
\quad\hbox{with}\quad
b_{p_1}\!\cdots b_{p_m}=1.
\end{equation}

If $P$ is infinite, write $P=\{p_1,p_2,\ldots\}$ for distinct elements~$p_i$, such that for every $n\ge0$ we have $P_n=\{p_1,\ldots,p_{m_n}\}$ for some $m_n\ge0$. We can then apply the same argument to ${\pi_1^\et(\BP^1_k \setminus P_n,x_0)}$ for every $n\ge0$, which naturally acts on the finite tree $T_n$ of level $\le n$. If $b_{p_i,n} \in W_n = \Aut(T_n)$ denotes the image of a suitable generator of an inertia group above $p_i\in P_n$, we deduce that $G^\geom_n = \langle b_{p_1,n},\ldots,b_{p_{m_n},n}\rangle$ and $b_{p_1,n}\cdots b_{p_{m_n},n}=1$.
Passing to the limit over $n$ we obtain elements $b_{p_i}\in W$ representing the inertia above $p_i\in P$.
These elements satisfy $b_{p_i}|_{T_n}=1$ whenever $i>m_n$; hence the sequence $b_{p_1},b_{p_2},\ldots$ converges to the identity element $1\in W$. The equations $b_{p_1,n}\cdots b_{p_{m_n},n}=1$ for all $n$ thus imply that the infinite product $b_{p_1}\!\!\cdot\! b_{p_2}\cdots$ also converges to~$1$. Together we deduce that in this case
\UseTheoremCounterForNextEquation
\begin{equation}\label{17GgeomGensInfin}
G^\geom = \Langle b_{p_1},b_{p_2},\ldots\Rangle
\quad\hbox{with}\quad
b_{p_1}\!\cdot b_{p_2}\!\cdots\,=1.
\end{equation}

\medskip
The generators thus constructed enjoy the following recursion relations up to conjugacy:

\begin{Prop}\label{17GenRec1}
For any $p\in P$ the element $b_p$ is conjugate under $W$ to
$$\left\{\begin{array}{ll}
\ \ \ \sigma & \hbox{if $p=f(c)$ for some $c\in C\setminus P$,}\\[3pt]
(b_c,1)\,\sigma & \hbox{if $p=f(c)$ for some $c\in C\cap P$,}\\[3pt]
(b_q,1) & \hbox{if $p=f(q)$ for a unique $q\in P\setminus C$,} \\[3pt]
(b_q,b_{q'}) & \hbox{if $p=f(q)=(q')$ for distinct $q,q'\in P\setminus C$.}
\end{array}\right\}$$
\end{Prop}

\begin{Proof}
First note that if $p=f(c)$ for some critical point $c\in C$, then $c$ is unique and $p\not=f(q)$ for all $q\in P\setminus C$, because $f$ is totally ramified at~$c$. Thus the first two cases are mutually exclusive and exclusive from the other two.
If $p\not=f(c)$ for any $c\in C$, the definition (\ref{17PostCritDef}) of $P$ implies that $p=f(q)$ for some $q\in P\setminus C$. But the number of such $q$ is always $\le2$, because $f$ is a morphism of degree~$2$. Thus the four cases are mutually exclusive and cover all possibilities.

Suppose that $p=f(q)$ for some $q\in P\setminus C$. Then the morphism $f$ is unramified over~$p$, and $f^{-1}(p)\subset\BP^1(\bar k)$ consists of $p$ and precisely one other point, say~$t$. The unramifiedness means that the inertia group above $p$ acts trivially on the covering $f\colon\BP^1_k\to\BP^1_k$, and hence that $b_p$ acts trivially on level $1$ of the tree~$T$. Thus $b_p=(u,v)$ for some elements $u,v\in W$. Also, under the functoriality of the fundamental group $f_*$ sends a generator of an inertia group above $q$ to a generator of an inertia group above~$p$. In view of Proposition \ref{13ConjPowers} this means that $T$ with the action of $\delta_q$ is isomorphic to one half subtree of $T$ with the action of $\delta_p$, or equivalently, that $b_q$ is conjugate to $u$ or~$v$. 
The other of the two entries $u,v$ arises in the same way from the inertia at~$t$. For this there are two possibilities. 

If $p=f(q')$ for some $q'\in P\setminus C$ distinct from~$q$, then $q'=t$ and the same argument as above shows that the other entry among $u,v$ is conjugate to~$b_{q'}$. Thus $b_p$ is conjugate to $(b_q,b_{q'})$ or to $(b_{q'},b_q)$. As these two elements of $W$ are already conjugate to each other under $\sigma\in W$, it follows that anyway $b_p$ is conjugate to $(b_q,b_{q'})$.

If $q$ is the only element of $P\setminus C$ with $p=f(q)$, the other point $t$ lies outside~$P$. Then the whole tower of coverings is unramified above~$t$, and so the other entry among $u,v$ is~$1$. In that case we find that $b_p$ is conjugate to $(b_q,1)$ or $(1,b_q)$, and hence anyway to $(b_q,1)$.

Now suppose that $p=f(c)$ for some $c\in C$. Then the morphism $f$ is ramified over~$p$. Thus the inertia group above $p$ acts non-trivially on the covering $f\colon\BP^1_k\to\BP^1_k$, and hence $b_p$ acts non-trivially on level $1$ of the tree~$T$. Therefore $b_p=(u,v)\,\sigma$ for some elements $u,v\in W$. 
Also, since $f$ has ramification degree $2$ at~$c$, under the functoriality of the fundamental group $f_*$ sends a generator of an inertia group above $c$ to the square of a generator of an inertia group above~$p$. Since $b_p^2 = (u,v)\,\sigma\,(u,v)\,\sigma = (uv,vu)$, using Proposition \ref{13ConjPowers} we deduce that the image of a generator of an inertia group above $c$ is conjugate to $uv$ under~$W$. Again there are now two possibilities. 

If $c\in P$, this means that $b_c$ is conjugate to~$uv$. By Lemma \ref{13ConjEquiv} this implies that $b_p = (u,v)\,\sigma$ is conjugate to $(uv,1)\,\sigma$ and hence to $(b_c,1)\,\sigma$ under~$W$.

If $c\not\in P$, the whole tower of coverings is unramified above~$c$, and so $uv=1$. By Lemma \ref{13ConjEquiv} this implies that $b_p = (u,v)\,\sigma$ is conjugate to $(uv,1)\,\sigma = \sigma$ under~$W$.

We have therefore proved the proposition in all cases.
\end{Proof}


\subsection{Specialization}
\label{18aSpecialization}

We keep the notation of Subsection \ref{17Monodromy}. 
First consider an intermediate field $k\subset k'\subset\bar k$ and a point $x'\in \BP^1(k')\setminus P$.
Viewing $x'$ as a point in $\BP^1(\bar k)$, we can construct the associated regular rooted binary tree $T_{x'}$ with set of vertices $\coprod_{n\ge0}f^{-n}(x')$, as in Subsection \ref{17Monodromy} with $\bar k$ in place of~$L$. This time we are interested in the natural action of $\Gal(\bar k/k')$ on~$T_{x'}$. Identify $T_{x'}$ with the standard tree $T$ in some way and let $G_{x'} \subset W$ denote the image of the continuous homomorphism $\Gal(\bar k/k')\to W$ describing the action on~$T_{x'}$. 
The following fact gives an upper bound for $G_{x'}$ in terms of the group $G^\arith$ associated to~$f$:

\begin{Prop}\label{18SpecDecomp}
There exists $w\in W$ such that 
$$G_{x'}\subset wG^\arith w^{-1}.$$
\end{Prop}

\begin{Proof}
Recall that changing the base point $x_0$ and/or the identifications of trees in Subsection \ref{17Monodromy} amounts to conjugating everything by an element of~$W$. Thus without loss of generality we may assume that $x_0=x'$ and that the identification of trees used in Subsection \ref{17Monodromy} is the same as above. Then the action of $\Gal(\bar k/k')$ on $T_{x'}$ is simply the composite homomorphism
$$\xymatrix@C-10pt{
\Gal(\bar k/ k') \ar@{}[r]|-= & 
\pi_1^\et(\Spec k',x_0)
\ar[rr] && 
\pi_1^\et(\BP^1_k \setminus P,x_0)
\ar[rr]^-\rho && W,\\}$$
where the arrow in the middle arises from the functoriality of the \'etale fundamental group. The inclusion for the images is then obvious.
\end{Proof}

\medskip
We are also interested in the effect on $G^\geom$ and $G^\arith$ of specializing the morphism $f$ itself. For this consider a noetherian normal integral domain $R$ with residue field $k$ of characteristic $\not=2$ and quotient field~$K$. Consider a morphism $F\colon \BP^1_S\to\BP^1_S$ over $S := \Spec R$ which is fiberwise of degree~$2$. Let $F$ also denote the induced quadratic morphism over~$K$, and $f$ the induced quadratic morphism over~$k$. Let $G_F^\geom\subset G_F^\arith\subset W$ and $G_f^\geom\subset G_f^\arith\subset W$ denote their associated monodromy groups. 

\begin{Prop}\label{18SpecSpec}
There exists $w\in W$ such that 
$$G_f^\geom \subset wG_F^\geom  w^{-1} \qquad\hbox{and}\qquad
  G_f^\arith\subset wG_F^\arith w^{-1} .$$
\end{Prop}

\begin{Proof}
Let $P\subset \BP^1_S$ denote the postcritical orbit of $F$ over~$S$, so that $P_K$ and $P_k$ are the respective postcritical orbits of $F$ over~$K$ and of~$f$. Choose geometric points $x_0$ of $\BP^1_k \setminus P_k$ and $X_0$ of $\BP^1_K\setminus P_K$. Then the respective monodromy representations and the functoriality of the \'etale fundamental group yield a commutative diagram
$$\xymatrix@C-10pt@R+5pt{
\pi_1^\et(\BP^1_K {\setminus} P_K,X_0) \ar[r] \ar[dr] \ar@{->>}[dd] &
\pi_1^\et(\BP^1_S {\setminus} P,  X_0) \ar@{=}[r]^-\sim \ar[d] &
\pi_1^\et(\BP^1_S {\setminus} P,  x_0) \ar[d] &
\pi_1^\et(\BP^1_k {\setminus} P_k,x_0) \ar[l] \ar[dl] \ar@{->>}[dd] \\
\qquad\qquad\qquad\qquad & \Aut(T_{X_0}) \ar@{=}[r]^-\sim \ar@{=}[d]^\wr 
& \Aut(T_{x_0}) \ar@{=}[d]^\wr & \qquad\qquad\qquad\qquad \\
G^\arith_F \ar@{^{ (}->}[r] & W \ar@{=}[r]^-\sim & W &
G^\arith_f \ar@{_{ (}->}[l] \\ }$$
Here the automorphism of $W$ in the lower row is the inner automorphism resulting from the respective identifications $T_{X_0}\cong T\cong T_{x_0}$. Since $\BP^1_S \setminus P$ is noetherian normal integral, the homomorphism on the upper left is an isomorphism by SGA1 \cite[exp.V, Props. \hbox{8.1--2}]{SGA1}. The commutativity of the diagram thus shows the desired assertion for the arithmetic monodromy groups.
For the geometric monodromy groups it follows from the commutative diagram and the fact that both $G_F^\geom$ and $G_f^\geom$ are topologically generated by the images of the tame inertia groups along~$P$.
\end{Proof}


\subsection{Polynomial case}
\label{18Poly}

We keep the notation of Subsection \ref{17Monodromy}, but from now on we assume that $f$ is a quadratic polynomial. Then one of the critical points is the point $\infty \in \BP^1(k)$ and satisfies $f(\infty)=\infty$. The other critical point is then necessarily a point $p_0\in\BA^1(k) = k$. Abbreviate $p_i := f^i(p_0)$ for all $i\ge1$, so that $P=\{\infty\}\sqcup P'$ where $P' := \{p_i\mid i\ge1\}$ is the  strict forward orbit of~$p_0$. Then the postcritical orbit $P$ as a directed graph can be classified as follows:

\begin{Class}\label{17Cases}
\rm If $P$ is infinite, the points $p_1,p_2,\ldots$ are all distinct. Otherwise let $r$ denote the cardinality of~$P'$. Then the points $p_1,\ldots,p_r$ are all distinct, and $p_{r+1} = p_{s+1}$ for some unique $0\le s<r$. If $s=0$, the equality $p_{r+1}=p_1$ implies that $p_r=p_0$, because $p_0$ is the only point above~$p_1$. In this case the critical point $p_0$ is called \emph{periodic}. If $s>0$, the point $p_0$ is \emph{strictly pre-periodic} in the sense that it is not periodic but the element $p_{s+1}$ in its forward orbit is. All in all we thus have the following three cases:
$$\begin{array}{lll}
\hbox{(a) Infinite case} & 
\hbox{(b) Periodic case} & 
\hbox{(c) Strictly pre-periodic case} \\[5pt]
\fbox{\ $
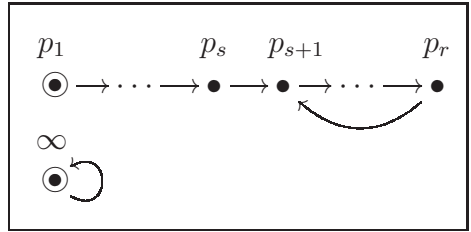
\xymatrix@R-26pt@C-12pt{
p_1 & p_2 & p_3 & \phantom{p_n} \raise3pt\hbox{$\mathstrut$}\\
\circbullet\ar[r]&\bullet\ar[r]&\bullet\ar[r]&\cdots \\
\phantom{x} &&& \\
\infty  &&& \\
\!\circbullet\! \ar@(dr,ur)[] 
&&& {\raise-5pt\hbox{\strut}} \\
}$}
\ \ &
\fbox{\ $\xymatrix@R-26pt@C-12pt{
p_1 & p_2 & \phantom{p_n} & p_r \raise3pt\hbox{$\mathstrut$}  \\ 
\circbullet\ar[r]&\bullet\ar[r]&\cdots\ar[r]&\bullet\ar@/^20pt/[lll] \\
\phantom{x} &&& \\
\infty &&& \\
\!\circbullet\! \ar@(dr,ur)[] &&& {\raise-5pt\hbox{\strut}} \\
}$}
\ \ &
\fbox{\ $\xymatrix@R-26pt@C-12pt{
p_1 & \phantom{p_n} & p_s & \kern-5pt p_{s+1}\kern-15pt & \phantom{p_n} & p_r \raise3pt\hbox{$\mathstrut$}  \\ 
\circbullet \ar[r]&\cdots\ar[r]&\bullet\ar[r]&\bullet\ar[r]&\cdots\ar[r]&\bullet\ar@/^16pt/[ll] \\
\phantom{x} &&&&& \\
\infty &&&&& \\
\!\circbullet\! \ar@(dr,ur)[] &&&&& {\raise-5pt\hbox{\strut}} \\
}$} 
\end{array}$$
\end{Class}

\medskip
Applying Proposition \ref{17GenRec1} to these cases we deduce:

\begin{Prop}\label{17GenRec}
If $f$ is a quadratic polynomial, the generators of $G^\geom$ constructed above satisfy the following conjugacy relations under~$W$, according to the cases in \ref{17Cases}:
$$\left\{\begin{array}{l}
\hbox{$b_{p_1}$ is conjugate to $\sigma$ in the cases (a) and (c).}\\[3pt]
\hbox{$b_{p_1}$ is conjugate to $(b_{p_r},1)\,\sigma$ in the case (b).}\\[3pt]
\hbox{$b_{p_i}$ is conjugate to $(b_{p_{i-1}},1)$ if $i>1$, except if $i=s+1$ in the case (c).}\\[3pt]
\hbox{$b_{p_{s+1}}$ is conjugate to $(b_{p_s},b_{p_r})$ in the case (c).}\\[3pt]
\hbox{$b_\infty$ is conjugate to $(b_\infty,1)\,\sigma$.}
\end{array}\right\}$$
\end{Prop}

\medskip
Our determination of $G^\geom$ in the quadratic polynomial case will be based solely on the weak recursion relations from Proposition \ref{17GenRec}. The results in the respective cases are Theorems \ref{18GgeomarithThm}, \ref{28GgeomThm}, and \ref{38GgeomThm}. For the description of $G^\arith$ see Theorems \ref{18GgeomarithThm}, \ref{28GarithThm}, \ref{38GarithThm1}, and Corollary \ref{38GarithThm2}.
Also note that with a suitable indexing of the postcritical orbit the product relation in (\ref{17GgeomGensFin}) or (\ref{17GgeomGensInfin}) is equivalent to
\UseTheoremCounterForNextEquation
\begin{equation}\label{17GgeomGensRelPol}
b_\infty = 
\biggl\{\begin{array}{ll}
(b_{p_1}\!\cdot b_{p_2}\!\cdots)^{-1} & \hbox{in the case (a), respectively} \\[3pt]
(b_{p_1}\!\cdots b_{p_r})^{-1} & \hbox{in the cases (b) and (c).}\\
\end{array}
\end{equation}
We can thus drop the generator~$b_\infty$, obtaining
\UseTheoremCounterForNextEquation
\begin{equation}\label{17GgeomGensPol}
G^\geom = 
\biggl\{\begin{array}{ll}
\Langle b_{p_1},b_{p_2},\ldots\Rangle & \hbox{in the case (a), respectively} \\[3pt]
\Langle b_{p_1},\ldots,b_{p_r}\Rangle & \hbox{in the cases (b) and (c).}\\
\end{array}
\end{equation}
On the other hand, by Proposition \ref{16OdoEquiv} the last line of Proposition \ref{17GenRec} means that $b_\infty$ is an odometer. We will exploit this fact in Sections \ref{28Monodromy} and \ref{38Monodromy}.


\subsection{Infinite polynomial case}
\label{19Infinite}

This is the case where $f$ is a quadratic polynomial with an infinite postcritical orbit. Thus $G^\geom = \Langle b_{p_1},b_{p_2},\ldots\Rangle$, where by Proposition \ref{17GenRec} the generators satisfy:
\UseTheoremCounterForNextEquation
\begin{equation}\label{18RecRels}
\left\{\begin{array}{l}
\hbox{$b_{p_1}$ is conjugate to $\sigma$ under~$W$, and}\\[3pt]
\hbox{$b_{p_i}$ is conjugate to $(b_{p_{i-1}},1)$ under~$W$ for any $i>1$.}
\end{array}\right\}
\end{equation}

\begin{Thm}\label{18GgeomarithThm}
In this case $G^\geom=G^\arith=W$.
\end{Thm}

\begin{Proof}
The same calculation as in Example \ref{15Ex} shows that for all $n,i\ge1$ we have
$$\sgn_n(b_{p_i}) \ =\ 
\biggl\{\!\begin{array}{rl}
-1 & \hbox{if $n    =i$,}\\[3pt]
 1 & \hbox{if $n\not=i$}
\end{array}$$
Thus Proposition \ref{15SignGnSurj} (b) implies that $G^\geom=W$. Since ${G^\geom\subset G^\arith\subset W}$, both inclusions must be equalities.
\end{Proof}


%
%
%
%
%
%
%

%
%

\newpage
\section{Periodic polynomial case}
\label{2Periodic}

Throughout this section we fix an integer $r\ge1$. Using Proposition \ref{15RecRelsProp} we define elements $a_1,\ldots,a_r\in W$ by the recursion relations
\UseTheoremCounterForNextEquation
\begin{equation}\label{2RecRels}
\biggl\{\begin{array}{ll}
a_1 = (a_r,1)\,\sigma, & \\[3pt]
a_i = (a_{i-1},1) & \hbox{for $2\le i\le r$.}\\
\end{array}\biggr\}
\end{equation}
We will study the subgroup $G:= \Langle a_1,\ldots,a_r\Rangle\subset W$ and its images $G_n\subset W_n$ for all~$n$.
The recursion relations are only one choice encoding the qualitative recursive behavior of the generators described by the diagram
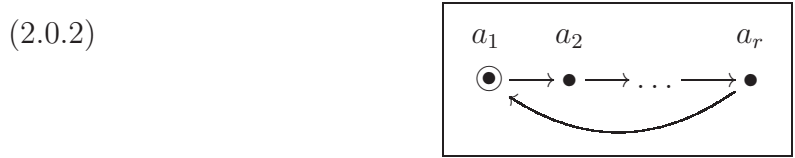
\UseTheoremCounterForNextEquation
\begin{equation}\label{2Diagram}
\fbox{\ $\xymatrix@R-24pt@C-10pt{
a_1 & a_2 && \ a_r\ \mathstrut \\ 
\circbullet\ar[r]&\bullet\ar[r]&\ldots\ar[r]&\bullet\ar@/^20pt/[lll] \\
{\vphantom{X}} \mathstrut &&&\ \\
}$}\qquad\\
\end{equation}
where the vertices are now marked by the generators and the superfluous generator corresponding to the critical point $\infty$ is dropped: see Subsection \ref{17Monodromy}. In Subsection \ref{24Conjugacy} we will show that many other choices yield the same group $G$ up to conjugacy in~$W$.


\subsection{Basic properties}
\label{21Basic}

\begin{Prop}\label{21GenSigns}
For all $n\ge1$ and all $1\le i\le r$ we have
$$\sgn_n(a_i) \ =\ 
\biggl\{\!\begin{array}{rl}
-1 & \hbox{if $n    \equiv i\mod r$,}\\[3pt]
 1 & \hbox{if $n\not\equiv i\mod r$.}
\end{array}$$
\end{Prop}

\begin{Proof}
The definition of $a_i$ implies that $\sgn_1(a_1)=-1$ and $\sgn_1(a_i) = 1$ for all $2\le i\le r$; hence the proposition holds for $n=1$. For $n>1$ the formulas (\ref{15SignRec}) and (\ref{15SignSigma}) imply that $\sgn_n(a_1) = \sgn_{n-1}(a_r)$ and $\sgn_n(a_i) = \sgn_{n-1}(a_{i-1})$ for all $2\le i\le r$. Thus the proposition for $n-1$ implies the proposition for~$n$, and so we are done by induction.
\end{Proof}

\begin{Prop}\label{21GnAll}
For any $n\ge0$ we have $G_n=W_n$ if and only if $n\le r$.
\end{Prop}

\begin{Proof}
Direct consequence of Propositions \ref{15SignGnSurj} (a) and \ref{21GenSigns}.
\end{Proof}

\begin{Prop}\label{21LevelTrans}
The group $G$ acts transitively on the level $n$ of $T$ for every $n\ge0$.
\end{Prop}

\begin{Proof}
Consider the element $w := a_1\cdots a_r\in G$. Proposition \ref{21GenSigns} implies that $\sgn_n(w)=-1$ for all $n\ge1$; hence $w$ already acts transitively on every level by Proposition \ref{16OdoEquiv}.
\end{Proof}

\begin{Prop}\label{21GenOrder}
Each generator $a_i$ has infinite order.
\end{Prop}

\begin{Proof}
If not, there exist $1\le i\le r$ and $n\ge0$ with $a_i^{2^n}=1$. Among all such pairs $(i,n)$ select one for which $i+rn$ is minimal. If $i>1$, the recursion relation implies that $(a_{i-1}^{2^n},1) = a_i^{2^n} = 1$ and hence $a_{i-1}^{2^n}=1$, contradicting the minimality of $i+rn$. Thus we must have $i=1$. Since $a_1\not=1$, we then have $n\ge1$. We can therefore calculate 
$$1 = a_1^{2^n} = (a_1^2)^{2^{n-1}} \!
= \bigl( (a_r,1)\,\sigma\,(a_r,1)\,\sigma \bigr)^{2^{n-1}} \!
= \bigl( (a_r,a_r)\bigr)^{2^{n-1}} \!
= (a_r^{2^{n-1}},a_r^{2^{n-1}})$$
and deduce that $a_r^{2^{n-1}}=1$.
Since $r+r(n-1) = rn < 1+rn$, this again contradicts the minimality of~$i+rn$. Thus we have a contradiction in all cases, as desired.
\end{Proof}

\begin{Lem}\label{21GenRes}
For any $1\le i\le r$ and $n\ge0$ we have $a_i|_{T_n}=1$ if and only if $n<i$.
\end{Lem}

\begin{Proof}
Direct consequence of the recursion relations (\ref{2RecRels}).
\end{Proof}

\medskip
More precise statements about the generators are given in Lemma \ref{22HinSemiDirect} (a) below.


\subsection{Useful subgroups}
\label{22Subgroups}

First we consider the subgroup of index $2$ which acts trivially on level $1$:
\UseTheoremCounterForNextEquation
\begin{equation}\label{22G1Def}
G^1\ :=\ G\cap(W\times W).
\end{equation}
Let $\proj_1$ and $\proj_2\colon W\times W\to W$ denote the two projections. We have the following self-similarity properties:

\begin{Prop}\label{22G1Prop}
\begin{enumerate}
\item[(a)] $G\subset (G\times G)\rtimes \langle\sigma\rangle$.
\item[(b)] $G^1\subset G\times G$.
\item[(c)] $\proj_1(G^1)=\proj_2(G^1)=G$.
\end{enumerate}
\end{Prop}

\begin{Proof}
By the recursion relations (\ref{2RecRels}), the generators of $G$ lie in $(G\times G)\rtimes \langle\sigma\rangle$, hence so does~$G$, proving (a). Also, (a) directly implies (b). 
In particular $G^1$ contains the elements $a_1^2 = (a_r,1)\,\sigma\,(a_r,1)\,\sigma = (a_r,a_r)$ as well as $a_i = (a_{i-1},1)$ for all $2\le i\le r$. Thus $\proj_1(G^1)$ contains the elements $a_r$ as well as $a_{i-1}$ for all $2\le i\le r$, that is, all the generators of~$G$, proving that $\proj_1(G^1)=G$. Finally, conjugation by the element $a_1 = (a_r,1)\,\sigma \in G$ interchanges the two factors of $W\times W$ and normalizes~$G^1$; hence also $\proj_2(G^1)=G$, and we are done.
\end{Proof}

%

\medskip
Next, for every $1\le i\le r$ we consider the subgroup
\UseTheoremCounterForNextEquation
\begin{equation}\label{22HiDef}
H^{(i)}\ :=\ 
\biggl\{\begin{array}{l}
\hbox{closure of the subgroup of $G$ generated} \\[3pt]
\hbox{by all $G$-conjugates of $a_j$ for all $j\not=i$.}
\end{array}\biggr\}_.
\end{equation}
By construction this is a normal subgroup of~$G$, and $G$ is the product $H^{(i)}\cdot\Langle a_i\Rangle$.

\begin{Lem}\label{22HiProp}
\begin{enumerate}
\item[(a)] $H^{(1)} = H^{(r)}\times H^{(r)}$.
\item[(b)] For any $i>1$ we have $H^{(i)}\cap (W\times W) \subset (H^{(i-1)}\times H^{(i-1)})\cdot\Langle(a_{i-1},a_{i-1}^{-1})\Rangle$.
\end{enumerate}
\end{Lem}

\begin{Proof}
By definition $H^{(1)}$ is the closure of the subgroup generated by the conjugates $za_jz^{-1}$ for all $2\le j\le r$ and $z\in G$. In the case $z=(x,y) \in G^1$ these conjugates are 
$$za_jz^{-1} = (x,y)\,(a_{j-1},1)\,(x,y)^{-1} = (xa_{j-1}x^{-1},1).$$
Here Proposition \ref{22G1Prop} (c) shows that $x$ runs through all of~$G$, so the closure of the subgroup generated by the $G^1$-conjugates only is $H^{(r)}\times1$. Conjugating this by $a_1 = (a_r,1)\,\sigma$ yields the subgroup $1\times H^{(r)}$. Thus all $G$-conjugates together yield the subgroup $H^{(r)}\times H^{(r)}$, proving (a).

Now consider the case $i>1$. Then $H^{(i-1)}\times H^{(i-1)}$ is a closed normal subgroup of $(G\times G)\rtimes\langle\sigma\rangle$, and the group $(\Langle a_{i-1}\Rangle\times\Langle a_{i-1}\Rangle)\rtimes\langle\sigma\rangle$ surjects to the associated factor group.
Also, the subgroup $\Langle(a_{i-1},a_{i-1}^{-1})\Rangle$ is normal in $(\Langle a_{i-1}\Rangle\times\Langle a_{i-1}\Rangle)\rtimes\langle\sigma\rangle$ and the factor group is abelian. Therefore $K := (H^{(i-1)}\times H^{(i-1)})\cdot\Langle(a_{i-1},a_{i-1}^{-1})\Rangle$ is a closed normal subgroup of $(G\times G)\rtimes\langle\sigma\rangle$ with an abelian factor group. 

For all $2\le j\le r$ with $j\not=i$ we have $a_{j-1}\in H^{(i-1)}$ and hence $a_j = (a_{j-1},1) \in K$. Also, we have $a_r\in H^{(i-1)}$ and hence $a_1^2 = \bigl((a_r,1)\,\sigma\bigr){}^2 = (a_r,a_r) \in K$. Moreover, since the factor group $\bigl((G\times G)\rtimes\langle\sigma\rangle\bigr)/K$ is abelian, the subgroup $\tilde K := K\cdot\langle a_1\rangle$ is again normal in $(G\times G)\rtimes\langle\sigma\rangle$.
In particular Proposition \ref{22G1Prop} (a) implies that $\tilde K$ is normalized by~$G$.
As we have just seen that $a_j\in\tilde K$ for all $1\le j\le r$ with $j\not=i$, from the definition of $H^{(i)}$ we deduce that $H^{(i)} \subset \tilde K$. Since $a_1\not\in W\times W$, but $a_1^2\in K$, we can conclude that $H^{(i)}\cap (W\times W) \subset \tilde K\cap (W\times W) = K$, as desired.
\end{Proof}

\begin{Lem}\label{22HinSemiDirect}
For any $n\ge0$ and $1\le i\le r$ let $H^{(i)}_n$ denote the image of $H^{(i)}$ in~$W_n$.
\begin{enumerate}
\item[(a)] The element $a_i|_{T_n}$ of $G_n$ has order $2^{\lfloor\frac{n+r-i}{r}\rfloor}$.
\item[(b)] We have a semidirect product $G_n= H^{(i)}_n \rtimes \langle a_i|_{T_n}\rangle$.
\end{enumerate}
\end{Lem}

\begin{Proof}
The assertions are trivially true for $n=0$, so assume that $n>0$ and that they hold for $n-1$ and all~$i$. 

Since $a_1$ acts non-trivially on level~$1$, its restriction $a_1|_{T_n}$ has order $>1$. From the induction hypothesis we can deduce that $(a_1|_{T_n})^2 = (a_r|_{T_{n-1}},a_r|_{T_{n-1}})$ has order $2^{\lfloor\frac{(n-1)+r-r}{r}\rfloor} = 2^{\lfloor\frac{n-1}{r}\rfloor}$. Together this implies that $a_1|_{T_n}$ has order $2\cdot2^{\lfloor\frac{n-1}{r}\rfloor} = 2^{\lfloor\frac{n+r-1}{r}\rfloor}$, proving (a) for $i=1$.
Also, Lemma \ref{22HiProp} (a) implies that $H^{(1)}_n = H^{(r)}_{n-1}\times H^{(r)}_{n-1}$, which in particular shows that $a_1|_{T_n} \not\in H^{(1)}_n$. Since we have $H^{(r)}_{n-1} \cap \langle a_r|_{T_{n-1}}\rangle = 1$ by the induction hypothesis, we can deduce that
$$H^{(1)}_n \cap \langle a_1|_{T_n}\rangle 
\ =\ H^{(1)}_n \cap \langle (a_1|_{T_n})^2\rangle 
\ =\ \bigl( H^{(r)}_{n-1}\times H^{(r)}_{n-1}\bigr) \cap \langle (a_r|_{T_{n-1}},a_r|_{T_{n-1}}) \rangle 
\ =\ 1.$$
Thus the product $G_n= H^{(1)}_n \cdot \langle a_1|_{T_n}\rangle$ is semidirect, proving (b) for $i=1$.

Now consider $i>1$. Then the induction hypothesis shows that $a_i|_{T_n} = (a_{i-1}|_{T_{n-1}},1)$ has order $2^{\lfloor\frac{(n-1)+r-(i-1)}{r}\rfloor} = 2^{\lfloor\frac{n+r-i}{r}\rfloor}$, as desired in (a). 
Moreover, Lemma \ref{22HiProp} (b) implies that $H^{(i)}_n \cap (W_{n-1}\times1) \subset H^{(i-1)}_{n-1}\times1$, and since $H^{(i-1)}_{n-1} \cap \langle a_{i-1}|_{T_{n-1}}\rangle = 1$ by the induction hypothesis, it follows that $H^{(i)}_n \cap \langle a_i|_{T_n}\rangle = 1$.
Thus the product $G_n= H^{(i)}_n \cdot \langle a_i|_{T_n}\rangle$ is semidirect, as desired in (b).

This proves the assertions for the given $n$ and all~$i$, and so they follow in general by induction.
\end{Proof}

\begin{Prop}\label{22HiSemiDirect}
For any $1\le i\le r$ we have a semidirect product $G= H^{(i)} \rtimes \Langle a_i\Rangle$, where $\Langle a_i\Rangle \cong \BZ_2$.
\end{Prop}

\begin{Proof}
Take the limit over $n$ in Lemma \ref{22HinSemiDirect}.
\end{Proof}

\begin{Thm}\label{22Gab}
The closure of the commutator subgroup of $G$ is $H^{(1)}\cap\ldots\cap H^{(r)}$, and the maximal profinite abelian factor group of $G$ is
$$ G_\ab\ =\ G\big/(H^{(1)}\cap\ldots\cap H^{(r)})
\ \stackrel{\sim}{\longto}\ 
G/H^{(1)} \times\ldots\times G/H^{(r)}
\ \cong\ \BZ_2^r.$$
\end{Thm}

\begin{Proof}
The isomorphisms in the middle and on the right hand side are direct consequences of Proposition \ref{22HiSemiDirect} and the fact that $a_j\in H^{(i)}$ whenever $j\not=i$. On the other hand, since $G$ is topologically generated by $a_1,\ldots,a_r$, its maximal abelian factor group is topologically generated by the images of $a_1,\ldots,a_r$. As their images form a basis of~$\BZ_2^r$, this is already the maximal profinite abelian factor group of $G$, as desired.
\end{Proof}

\medskip
In the following we abbreviate 
\UseTheoremCounterForNextEquation
\begin{equation}\label{22HDef}
H := H^{(r)}
\end{equation}
which by definition is the closure of the subgroup generated by all $G$-conjugates of $a_j$ for all $1\le j<r$. We let $H_n$ denote its image in~$W_n$, and $G^1_n$ the image of $G^1$ in~$W_n$.

\begin{Lem}\label{22HnSemiDirect}
We have the semidirect product decompositions:
\begin{enumerate}
\item[(a)] $G_n = H_n \rtimes \langle a_r|_{T_n}\rangle$ 
where $a_r|_{T_n}$ has order $2^{\lfloor\frac{n}{r}\rfloor}$.
\item[(b)] $G_n = (H_{n-1}\times H_{n-1}) \rtimes \langle a_1|_{T_n}\rangle$
where $a_1|_{T_n}$ has order $2^{\lfloor\frac{n+r-1}{r}\rfloor}$.
\item[(c)] $G^1_n = (H_{n-1}\times H_{n-1}) \rtimes \langle (a_r|_{T_{n-1}},a_r|_{T_{n-1}}) \rangle$ where $a_r|_{T_{n-1}}$ has order $2^{\lfloor\frac{n-1}{r}\rfloor}$.
\end{enumerate}
\end{Lem}

\begin{Proof}
Part (a) is the special case $i=r$ of Lemma \ref{22HinSemiDirect}, and (b) is the special case $i=1$ of Lemma \ref{22HinSemiDirect} combined with Lemma \ref{22HiProp} (a). Part (c) follows from (b) and the fact that $a_1|_{T_n} \not\in G^1_n$ and has square $(a_r|_{T_{n-1}},a_r|_{T_{n-1}})\in G^1_n$.
\end{Proof}

\medskip
In the limit over $n$ this implies the following more specific self-similarity properties of~$G$:

\begin{Prop}\label{22HSemiDirect}
We have the semidirect product decompositions:
\begin{enumerate}
\item[(a)] $G = H \rtimes \Langle a_r\Rangle$.
\item[(b)] $G = (H\times H) \rtimes \Langle a_1\Rangle$.
\item[(c)] $G^1 = (H\times H) \rtimes \Langle (a_r,a_r) \Rangle$.
\end{enumerate}
\end{Prop}



\subsection{Size}
\label{23Size}

\begin{Prop}\label{23GnOrder}
For all $n\ge0$ we have 
$$\log_2|G_n|\ =\ 2^n-1-\sum_{m=0}^{n-1} 
2^{n-1-m} \cdot {\textstyle\bigl\lfloor\frac{m}{r}\bigr\rfloor}.$$
\end{Prop}

\begin{Proof}
For $n>0$ Lemma \ref{22HnSemiDirect} (a) and (b) imply that
\begin{eqnarray*}
\log_2|G_{n-1}| &\!\!=\!\!& \textstyle \log_2|H_{n-1}| + \bigl\lfloor\frac{n-1}{r}\bigr\rfloor \quad\hbox{and}\\
\log_2|G_n| &\!\!=\!\!& \textstyle 2\cdot\log_2|H_{n-1}| + \bigl\lfloor\frac{n+r-1}{r}\bigr\rfloor.
\end{eqnarray*}
Together this yields the recursion formula
$$\textstyle \log_2|G_n|
\ =\ 2\cdot\bigl(\log_2|G_{n-1}| - \bigl\lfloor\frac{n-1}{r}\bigr\rfloor \bigr)  + \bigl\lfloor\frac{n+r-1}{r}\bigr\rfloor
\ =\ 2\cdot\log_2|G_{n-1}| + 1 - \bigl\lfloor\frac{n-1}{r}\bigr\rfloor.$$
Since $\log_2|G_0|=1$, a direct induction on $n$ yields the desired formula.
\end{Proof}

\begin{Thm}\label{23Hausdorff}
The Hausdorff dimension of $G$ exists and is $1 - \frac{1}{2^r-1}$.
\end{Thm}

\begin{Proof}
Proposition \ref{23GnOrder} implies that
$$\frac{\log_2|G_n|}{2^n-1}
\ =\ 1 - \frac{2^n}{2^n-1} \cdot \sum_{m=0}^{n-1} 2^{-1-m} \cdot {\textstyle\bigl\lfloor\frac{m}{r}\bigr\rfloor}.$$
By (\ref{12Hausdorff}) the Hausdorff dimension is its limit for $n\to\infty$, which comes out as
\begin{eqnarray*}
1 - \sum_{m=0}^{\infty} 2^{-1-m} \cdot {\textstyle\bigl\lfloor\frac{m}{r}\bigr\rfloor}
&=& 1 - \sum_{i=0}^{r-1} \sum_{j=0}^{\infty} 2^{-1-(i+rj)} \cdot {\textstyle\bigl\lfloor\frac{i+rj}{r}\bigr\rfloor} \\
&=& 1 - \sum_{i=0}^{r-1} 2^{-1-i} \cdot \sum_{j=0}^{\infty} j\cdot 2^{-rj} \\
&=& 1 - (1-2^{-r}) \cdot \frac{2^{-r}}{(1-2^{-r})^2} 
\ \ =\ \ 1 - \frac{1}{2^r-1},
\end{eqnarray*}
as desired.
\end{Proof}


\subsection{Conjugacy of generators}
\label{24Conjugacy}

The recursion relations (\ref{2RecRels}) are only one choice among many describing the same qualitative recursive behavior of the generators, and besides their simplicity there is nothing inherently better in them to set them apart from similar choices (for instance, compare \cite[\S3]{Bartholdi-Nekrashevych-2008}).
In this subsection we show in a precise sense that up to conjugacy in $W$ the group $G$ is actually the same for all recursion relations that are conjugate to the given ones. We can thus say that $G$ really depends only on the combinatorics of the diagram (\ref{2Diagram}), perhaps in the same way as a symplectic group of matrices of a given size $n$ depends up to conjugacy only on~$n$ but not on the alternating form, even though the choice of an alternating form cannot be avoided.

\begin{Thm}\label{24SemiRigid}
For any elements $b_1,\ldots,b_r\in W$ the following are equivalent:
\begin{enumerate}
\item[(a)] $\biggl\{\begin{array}{l}
\hbox{$b_1$ is conjugate to $(b_r,1)\,\sigma$ under $W$, and} \\[3pt]
\hbox{$b_i$ is conjugate to $(b_{i-1},1)$ under $W$ for any $2\le i\le r$.}\\
\end{array}\biggr\}$
\item[(b)] $b_i$ is conjugate to $a_i$ under $W$ for any $1\le i\le r$.
\item[(c)] There exist $w\in W$ and $x_i\in G$ such that $b_i=(wx_i)a_i(wx_i)^{-1}$ for any $1\le i\le r$.
\end{enumerate}
Moreover, for any $w$ as in (c) we have $\Langle b_1,\ldots,b_r\Rangle=wGw^{-1}$.
\end{Thm}

\begin{Proof}
The equivalence (a)$\Leftrightarrow$(b) is a special case of Proposition \ref{15RecConj=Conj}. The implication (c)$\Rightarrow$(b) is obvious, and the last statement follows from Lemma \ref{13ConjGen}. It remains to prove the implication (b)$\Rightarrow$(c). 
This follows by taking the limit over $n$ from the following assertion for all $n\ge0$: 
\begin{enumerate}
\item[($*_n$)] For any elements $b_i\in W_n$ conjugate to $a_i|_{T_n}$ under~$W_n$, there exist $w\in W_n$ and $x_i\in G_n$ such that $b_i=(wx_i)(a_i|_{T_n})(wx_i)^{-1}$ for each~$i$.
\end{enumerate}
This is trivial for $n=0$, so assume that $n>0$ and that ($*_{n-1}$) is true. 
Take elements $b_i\in W_n$ conjugate to $a_i|_{T_n}$ under~$W_n$.

Since $a_1 = (a_r,1)\,\sigma$ and $a_i = (a_{i-1},1)$ for all $i>1$, the $b_i$ have the form $b_1=(*,*)\,\sigma$ and $b_i=(*,1)$ or $(1,*)$ for $i>1$, where $*$ stands for some elements of~$W_{n-1}$. Set
$$b'_i\ :=\ 
\biggl\{\begin{array}{ll}
b_i & \hbox{if $i=1$ or $b_i=(*,1)$,} \\[3pt]
b_1^{-1}b_ib_1 & \hbox{if $i>1$ and $b_i=(1,*)$,}
\end{array}$$
which in either case is conjugate to $a_i|_{T_n}$ under~$W_n$. Suppose that ($*_n$) is true for $b'_1,\ldots,b'_r$, so that there exist $w\in W_n$ and $x'_1,\ldots,x'_r\in G_n$ such that $b'_i=(wx'_i)(a_i|_{T_n})(wx'_i)^{-1}$ for each~$i$. Set 
$$x_i\ :=\ 
\biggl\{\begin{array}{ll}
x'_i & \hbox{if $i=1$ or $b_i=(*,1)$,} \\[3pt]
x'_1(a_1|_{T_n})\,x_1^{\prime-1}x'_i & \hbox{if $i>1$ and $b_i=(1,*)$,}
\end{array}$$
which always lies in~$G_n$. Then if $i=1$ or $b_i=(*,1)$, we have $b_i=b'_i=(wx_i)(a_i|_{T_n})(wx_i)^{-1}$; otherwise we have
\begin{eqnarray*}
b_i &\!\!=\!\!& b_1b'_ib_1^{-1} \\
&\!\!=\!\!& 
(wx'_1)(a_1|_{T_n})(wx'_1)^{-1} \cdot 
(wx'_i)(a_i|_{T_n})(wx'_i)^{-1} \cdot 
(wx'_1)(a_1|_{T_n})(wx'_1)^{-1} \\
&\!\!=\!\!& 
wx'_1(a_1|_{T_n})\,x_1^{\prime-1}x'_i\cdot (a_i|_{T_n})\cdot 
x_i^{\prime-1}x'_1(a_1|_{T_n})\,x_1^{\prime-1}w^{-1} \\
&\!\!=\!\!& (wx_i)(a_i|_{T_n})(wx_i)^{-1}.
\end{eqnarray*}
Thus ($*_n$) holds for $b_1,\ldots,b_r$ with $x_1,\ldots,x_r$ in place of $x'_1,\ldots,x'_r$. On replacing $b_1,\ldots,b_r$ by $b'_1,\ldots,b'_r$ we have thus reduced ourselves to the case that $b_i=(*,1)$ for all $i>1$.

For each $i>0$ write $b_i=(c_{i-1},1)$ with $c_{i-1}\in W_{n-1}$. The fact that $b_i$ is conjugate to $a_i|_{T_n}=(a_{i-1}|_{T_{n-1}},1)$ under $W_n$ implies that $c_{i-1}$ is conjugate to $a_{i-1}|_{T_{n-1}}$ under~$W_{n-1}$. 
Also, choose an element $f\in W_n$ such that $b_1 = f(a_1|_{T_n})f^{-1}$. If $f$ acts non-trivially on level~$1$, we can replace it by $f(a_1|_{T_n})$ without changing the conjugate $b_1 = f(a_1|_{T_n})f^{-1}$. Afterwards $f$ has the form $f=(d,e)$ for some elements $d,e\in W_{n-1}$. We then have
$$b_1 \ =\ f(a_1|_{T_n})f^{-1} \ =\ (d,e)\,(a_r|_{T_{n-1}},1)\,\sigma\,(d,e)^{-1} 
\ =\ \bigl(d(a_r|_{T_{n-1}})e^{-1},ed^{-1}\bigr)\,\sigma.$$
Set $c_r := d(a_r|_{T_{n-1}})d^{-1} \in W_{n-1}$. Then for every $1\le i\le r$ the element $c_i$ is conjugate to $a_i|_{T_{n-1}}$ under $W_{n-1}$. By the induction hypothesis ($*_{n-1}$) there therefore exist $u\in W_{n-1}$ and $x_i\in G_{n-1}$ such that $c_i=(ux_i)(a_i|_{T_{n-1}})(ux_i)^{-1}$ for each~$i$. 
Put $w := (ux_r,ed^{-1}ux_r) \allowbreak {\in W_n}$. Then
\begin{eqnarray*}
w^{-1}b_1w
&\!\!=\!\!& \bigl(x_r^{-1}u^{-1},x_r^{-1}u^{-1}de^{-1}\bigr)\,
\bigl(d(a_r|_{T_{n-1}})e^{-1},ed^{-1}\bigr)\,\sigma
\, \bigl(ux_r,ed^{-1}ux_r\bigr)\\
&\!\!=\!\!& \bigl(x_r^{-1}u^{-1}d(a_r|_{T_{n-1}})e^{-1}ed^{-1}ux_r ,
x_r^{-1}u^{-1}de^{-1}ed^{-1}ux_r\bigr)\,\sigma \\
&\!\!=\!\!& \bigl(x_r^{-1}u^{-1}d(a_r|_{T_{n-1}})d^{-1}ux_r , 1\bigr)\,\sigma \\
&\!\!=\!\!& \bigl(x_r^{-1}u^{-1}c_rux_r , 1\bigr)\,\sigma \\
&\!\!=\!\!& \bigl(x_r^{-1}u^{-1} (ux_r)(a_r|_{T_{n-1}})(ux_r)^{-1} ux_r, 1\bigr)\,\sigma \\
&\!\!=\!\!& (a_r|_{T_{n-1}},1)\,\sigma \\
&\!\!=\!\!& a_1|_{T_n}.
\end{eqnarray*}
Also, for each $i>1$ we have 
\begin{eqnarray*}
w^{-1}b_iw
&\!\!=\!\!& \bigl(x_r^{-1}u^{-1},x_r^{-1}u^{-1}de^{-1}\bigr)\,
\bigl(c_{i-1},1\bigr)
\, \bigl(ux_r,ed^{-1}ux_r\bigr)\\
&\!\!=\!\!& \bigl(x_r^{-1}u^{-1}c_{i-1}ux_r,1\bigr) \\
&\!\!=\!\!& \bigl(x_r^{-1}u^{-1} (ux_{i-1})(a_{i-1}|_{T_{n-1}})(ux_{i-1})^{-1} ux_r,1\bigr) \\
&\!\!=\!\!& \bigl((x_r^{-1}x_{i-1})(a_{i-1}|_{T_{n-1}})(x_r^{-1}x_{i-1})^{-1},1\bigr).
\end{eqnarray*}
Since $x_r^{-1}x_{i-1}\in G_{n-1}$, by Proposition \ref{22G1Prop} (c) we can find $y_{i-1}\in W_{n-1}$ such that $z_i := (x_r^{-1}x_{i-1},y_{i-1})$ lies in~$G_n$. Then
\begin{eqnarray*}
w^{-1}b_iw
&\!\!=\!\!& 
(x_r^{-1}x_{i-1},y_{i-1})\,(a_{i-1}|_{T_{n-1}},1)\,(x_r^{-1}x_{i-1},y_{i-1})^{-1} \ \ \\ 
&\!\!=\!\!& z_i (a_i|_{T_n}) z_i^{-1}.
\end{eqnarray*}
With $z_1:=1$ we have thus shown that $w^{-1}b_iw = z_i (a_i|_{T_n}) z_i^{-1}$ or equivalently that $b_i = (wz_i) (a_i|_{T_n}) (wz_i)^{-1}$ for all $1\le i\le r$, where $w\in W_n$ and $z_i\in G_n$. This proves ($*_n$), and so we are done by induction.
\end{Proof}

\begin{Prop}\label{23WConj=GConjPower}
For any single $1\le i\le r$ and any $b_i\in G$ that is conjugate to $a_i$ under~$W$, there exists $k\in\BZ_2^\times$ such that $b_i$ is conjugate to $a_i^k$ under~$G$.
\end{Prop}


\begin{Proof}
By Lemma \ref{13ConjLimit} it suffices to show that for all $1\le i\le r$ and $n\ge0$:
\begin{enumerate}
\item[($*_{i,n}$)] For any $b_i\in G_n$ that is conjugate to $a_i|_{T_n}$ under~$W_n$, there exists $k\in\BZ_2^\times$ such that $b_i$ is conjugate to $a_i^k|_{T_n}$ under~$G_n$. 
\end{enumerate}
This is trivial for $n=0$, so consider $n>0$. By induction on $n$ it suffices to prove that ($*_{i-1,n-1}$) implies ($*_{i,n}$) if $i\ge2$, and that ($*_{r,n-1}$) implies ($*_{1,n}$). Consider an element $b_i\in G_n$, and an element $w\in W_n$ such that $b_i = w\,(a_i|_{T_n})\,w^{-1}$. 

Suppose first that $i\ge2$. If $w$ acts non-trivially on level~$1$, we replace $b_i$ by its $G_n$-conjugate $(a_1|_{T_n})\,b_i\,(a_1|_{T_n})^{-1}$ and $w$ by $(a_1|_{T_n})\,w$, which does not change the desired assertion. Afterwards $w$ has the form $w=(u,v)$ for some elements $u,v\in W_{n-1}$. Then the recursion relation yields
$$b_i = w\,(a_i|_{T_n})\,w^{-1}
= (u,v)\,(a_{i-1}|_{T_{n-1}},1)\,(u,v)^{-1}
= \bigl(u(a_{i-1}|_{T_{n-1}})u^{-1},1\bigr).$$
By assumption $b_i$ lies in~$G_n$, so Proposition \ref{22G1Prop} shows that $b_{i-1} := u(a_{i-1}|_{T_{n-1}})u^{-1}$ lies in~$G_{n-1}$. Thus $b_{i-1}$ is an element of $G_{n-1}$ that is conjugate to $a_{i-1}|_{T_{n-1}}$ under $W_{n-1}$; hence by the hypothesis ($*_{i-1,n-1}$) there exist $k\in\BZ_2^\times$ and $x\in G_{n-1}$ such that $b_{i-1} = x(a_{i-1}^k|_{T_{n-1}})x^{-1}$. By Proposition \ref{22G1Prop} (c)  we can choose an element $y\in G_{n-1}$ for which $z := (x,y)$ lies in~$G_n$. Then
$$b_i = (b_{i-1},1)
= \bigl(x(a_{i-1}^k|_{T_{n-1}})x^{-1},1\bigr)
= (x,y)\,(a_{i-1}^k|_{T_{n-1}},1)\,(x,y)^{-1}
= z\,(a_i^k|_{T_n})\,z^{-1},$$
which proves that $b_i$ is conjugate to $a_i^k|_{T_n}$ under~$G_n$, proving ($*_{i,n}$), as desired.

Now suppose that $i=1$. If $w$ acts non-trivially on level~$1$, we can replace it by $w\,(a_1|_{T_n})$ without changing the conjugate $b_1 = w\,(a_1|_{T_n})\,w^{-1}$. Afterwards $w$ has the form $w=(u,v)$ for some elements $u,v\in W_{n-1}$. Then the recursion relation yields
$$b_1^2 = w\,(a_1^2|_{T_n})\,w^{-1}
= (u,v)\,(a_r|_{T_{n-1}},a_r|_{T_{n-1}})\,(u,v)^{-1}
= \bigl(u(a_r|_{T_{n-1}})u^{-1},v(a_r|_{T_{n-1}})v^{-1}\bigr).$$
By assumption $b_1$ and hence $b_1^2$ lie in~$G_n$, so Proposition \ref{22G1Prop} shows that $b_r := u(a_r|_{T_{n-1}})u^{-1}$ lies in~$G_{n-1}$. Thus $b_r$ is an element of $G_{n-1}$ that is conjugate to $a_r|_{T_{n-1}}$ under $W_{n-1}$; hence by the hypothesis ($*_{r,n-1}$) there exist $k\in\BZ_2^\times$ and $x\in G_{n-1}$ such that $b_r = x(a_r^k|_{T_{n-1}})x^{-1}$. By Proposition \ref{22G1Prop} (c) we can choose an element $y\in G_{n-1}$ for which $z := (x,y)$ lies in~$G_n$. Then $b_1^2 = z\,(a_r|_{T_{n-1}},*)^k\,z^{-1}$, where $*$ represents an element of $W_{n-1}$ that we do not care about for the moment. Set $c_1 := z^{-1}\,b_1^{\kern1pt k^{-1}}\kern-1pt z$, so that $c_1^2 = (a_r|_{T_{n-1}},*)$. Since $b_1$ and $z$ lie in~$G_n$, so does~$c_1$. Moreover, Proposition \ref{13ConjPowers} implies that $b_1$ is conjugate to $b_1^{\kern1pt k^{-1}}$ and hence to $c_1$ under~$W_n$. Thus $c_1$ is again conjugate to $a_1|_{T_n}$ under~$W_n$. Since $b_1=zc_1^kz^{-1}$ with $z\in G_n$, the desired assertion for $b_1$ now reduces to that for~$c_1$. We can therefore replace $b_1$ by~$c_1$, after which we are reduced to the special case that $b_1^2 = (a_r|_{T_{n-1}},*)$.

Now we begin again by choosing a (new) element $w\in W_n$ such that $b_1 = w\,(a_1|_{T_n})\,w^{-1}$. If $w$ acts non-trivially on level~$1$, we replace it by $w\,(a_1|_{T_n})$ without changing the conjugate $b_1 = w\,(a_1|_{T_n})\,w^{-1}$. Afterwards $w$ has the form $w=(u,v)$ for some elements $u,v\in W_{n-1}$. Then the same calculation as before shows that
$(a_r|_{T_{n-1}},*) = b_1^2 = (u(a_r|_{T_{n-1}})u^{-1},*)$; hence $u$ commutes with $a_r|_{T_{n-1}}$. Therefore
\begin{eqnarray*}
b_1\, =\, w\,(a_1|_{T_n})\,w^{-1}
&\!\!=\!\!& (u,v)\,(a_r|_{T_{n-1}},1)\,\sigma\,(u,v)^{-1} \qquad\qquad \\
&\!\!=\!\!& \bigl(u(a_r|_{T_{n-1}})v^{-1},vu^{-1}\bigr)\,\sigma \\
&\!\!=\!\!& \bigl((a_r|_{T_{n-1}})uv^{-1},vu^{-1}\bigr)\,\sigma \\
&\!\!=\!\!& (a_r|_{T_{n-1}},1)\,\sigma\,(vu^{-1},uv^{-1})
\, =\, (a_1|_{T_n})\,(vu^{-1},uv^{-1}).
\end{eqnarray*}
Since $b_1$ lies in $G_n$ by assumption, this shows that $(vu^{-1},uv^{-1})$ lies in~$G^1_n$. By Lemma \ref{22HnSemiDirect} (c) this means that $vu^{-1} = h\,(a_r|_{T_{n-1}})^m$ for some $h\in H_{n-1}$ and $m\in\BZ$ such that $(a_r|_{T_{n-1}})^{-m} = (a_r|_{T_{n-1}})^m$. Therefore 
\begin{eqnarray*}
b_1 &\!\!=\!\!& (a_1|_{T_n})\,(vu^{-1},uv^{-1}) \\
&\!\!=\!\!& (a_1|_{T_n})\,\bigl(h\,(a_r|_{T_{n-1}})^m,(a_r|_{T_{n-1}})^{-m}h^{-1}\bigr) \\
&\!\!=\!\!& (a_r|_{T_{n-1}},1)\,\sigma\,(h,1)\,\bigl((a_r|_{T_{n-1}})^m,(a_r|_{T_{n-1}})^m\bigr)\,(1,h^{-1}) \\
&\!\!=\!\!& (a_r|_{T_{n-1}},1)\,(1,h)\,\sigma\,(a_r|_{T_{n-1}},a_r|_{T_{n-1}})^m\,(1,h)^{-1} \\
&\!\!=\!\!& (1,h)\,(a_r|_{T_{n-1}},1)\,\sigma\,(a_1|_{T_n})^{2m}\,(1,h)^{-1} \\
&\!\!=\!\!& (1,h)\,(a_1|_{T_n})^{1+2m}\,(1,h)^{-1}.
\end{eqnarray*}
Since $(1,h) \in G_n$ by Lemma \ref{22HnSemiDirect} (b), this shows that $b_1$ is conjugate to $a_1^{1+2m}|_{T_n}$ under~$G_n$, proving ($*_{r,n-1}$) and therefly finishing the proof of Proposition \ref{23WConj=GConjPower}.
\end{Proof}

%
%


\subsection{Small cases}
\label{25Small}

In the case $r=1$ the element $a_1$ is the odometer from Subsection \ref{16Odometer}, and $G = \Langle a_1\Rangle \cong \BZ_2$.
In all other cases the group $G$ has positive Hausdorff dimension and is very non-abelian.
The next smallest case $r=2$ is the closure of the so-called `Basilica group'; see Grigorchuk-$\dot{\rm Z}$uk \cite{GrigorchukZuk2002}. For the case $r=3$ corresponding to both the `Douady rabbit' and the `airplane' some of the results of the preceding subsection were obtained by Nekrashevych \cite[\S8]{Nekrashevych-2007}.

\medskip
The case $r=2$ stands out among the others by the rigidity property from Theorem \ref{25SmallRigid} below, which is stronger than the semirigidity from Theorem \ref{24SemiRigid}. The analogous statement does not hold for any $r>2$.


\medskip
For the following lemmas let $\Conjug_W(a_i)$ denote the conjugacy class of $a_i$ in~$W$ and $\Conjug_{W_n}(a_i|_{T_n})$ the conjugacy class of $a_i|_{T_n}$ in~$W_n$ for all $1\le i\le r$ and $n\ge1$.

\begin{Lem}\label{25SmallRigidLem1}
For $r=2$ and every $n\ge1$ we have 
$$\bigl| \Conjug_{W_n}(a_1|_{T_n}) \times \Conjug_{W_n}(a_2|_{T_n}) \bigr| 
\ =\ 2^{2^n-2}.$$
\end{Lem}

\begin{Proof}
Abbreviate $\gamma_{i,n} := \log_2\bigl| \Conjug_{W_n}(a_i|_{T_n}) \bigr|$. 
Since $\log_2|W_n|=2^n-1$ by (\ref{12WnOrder}), we have $\log_2\bigl| \Cent_{W_n}(a_i|_{T_n}) \bigr| = 2^n-1-\gamma_{i,n}$. For $n=1$ observe that $W_1$ is abelian; hence $\gamma_{i,1}=0$ for both $i=1,2$. 
For $n>1$ note first that $a_2|_{T_n} = (a_1|_{T_{n-1}},1)$ with $a_1|_{T_{n-1}}$ non-trivial; hence 
\UseTheoremCounterForNextEquation
\begin{equation}\label{25SRL11}
\Cent_{W_n}(a_2|_{T_n}) \ =\ \Cent_{W_{n-1}}(a_1|_{T_{n-1}}) \times W_{n-1}.
\end{equation}
Therefore 
\begin{eqnarray*}
2^n-1-\gamma_{2,n} 
&=& \log_2\bigl| \Cent_{W_n}(a_2|_{T_n}) \bigr| \\
&=& \log_2\bigl| \Cent_{W_{n-1}}(a_1|_{T_{n-1}}) \bigr| + 2^{n-1}-1 \\
&=& 2^{n-1}-1-\gamma_{1,n-1} + 2^{n-1}-1,
\end{eqnarray*}
and so $\gamma_{2,n} = \gamma_{1,n-1} + 1$.

Next the element $a_1|_{T_n}$ is contained in its own centralizer and acts non-trivially on level~$1$; hence 
$$\bigl|\Cent_{W_n}(a_1|_{T_n})\bigr| \ =\ 
2\cdot\bigl|\Cent_{W_{n-1}\times W_{n-1}}(a_1|_{T_n})\bigr|.$$
Also, an element $w=(u,v) \in W_{n-1}\times W_{n-1}$ commutes with $a_1|_{T_n}$ if and only if
\begin{eqnarray*}
(a_1|_{T_n})\,w
&=& (a_2|_{T_{n-1}},1)\,\sigma\,(u,v)
\ =\ ((a_2|_{T_{n-1}})v,u)\,\sigma \\
\Vert\qquad & \\
w\,(a_1|_{T_n})
&=& (u,v)\, (a_2|_{T_{n-1}},1)\,\sigma
\ =\ (u(a_2|_{T_{n-1}}),v)\,\sigma,
\end{eqnarray*}
that is, if and only if $u=v$ and $u$ commutes with $a_2|_{T_{n-1}}$.
Thus 
\UseTheoremCounterForNextEquation
\begin{equation}\label{25SRL12}
\Cent_{W_{n-1}\times W_{n-1}}(a_1|_{T_n}) 
\ =\ \diag\bigl(\Cent_{W_{n-1}}(a_2|_{T_{n-1}})\bigr).
\end{equation}
Therefore
\begin{eqnarray*}
2^n-1-\gamma_{1,n} 
&=& \log_2\bigl| \Cent_{W_n}(a_1|_{T_n}) \bigr| \\
&=& 1 + \log_2\bigl| \Cent_{W_{n-1}}(a_2|_{T_{n-1}}) \bigr| \\
&=& 1 + 2^{n-1}-1-\gamma_{2,n-1},
\end{eqnarray*}
and so $\gamma_{1,n} = \gamma_{2,n-1} + 2^{n-1} - 1$.

Finally, set $\gamma_n := \gamma_{1,n}+\gamma_{2,n}$. Then $\gamma_1=0$, and the above formulas imply that $\gamma_n = \gamma_{n-1} + 2^{n-1}$ for all $n>1$. By induction on $n$ it follows that $\gamma_n = 2^n-2$ for all $n\ge1$, as desired.
\end{Proof}

\begin{Lem}\label{25SmallRigidLem2}
For $r=2$ and every $n\ge1$ the centralizer of $G_n$ in $W_n$ 
has order~$2$.
\end{Lem}

\begin{Proof}
Since $W_1$ is abelian of order~$2$, the assertion holds for~$n=1$. For $n>1$ observe that since $a_1|_{T_n}$ and $a_2|_{T_n}$ generate~$G_n$, the equations (\ref{25SRL11}) and (\ref{25SRL12}) from the preceding proof imply that 
\begin{eqnarray*}
\Cent_{W_n}(G) 
&=& \Cent_{W_n}(a_1|_{T_n}) \cap \Cent_{W_n}(a_2|_{T_n}) \\
&=& \bigl( \Cent_{W_{n-1}}(a_1|_{T_{n-1}}) \times W_{n-1} \bigr)
    \cap \diag\bigl(\Cent_{W_{n-1}}(a_2|_{T_{n-1}})\bigr) \\
&=& \diag\bigl( \Cent_{W_{n-1}}(a_1|_{T_{n-1}}) \cap \Cent_{W_{n-1}}(a_2|_{T_{n-1}})\bigr) \\
&=& \diag\bigl( \Cent_{W_{n-1}}(G_{n-1})\bigr).
\end{eqnarray*}
Thus the order of the centralizer is independent of~$n$ and hence is always equal to~$2$.
\end{Proof}

\begin{Lem}\label{25SmallRigidLem3}
For $r=2$ the action of $W$ by joint conjugation on $\Conjug_W(a_1) \times \Conjug_W(a_2)$ is transitive.
\end{Lem}

\begin{Proof}
Let $W_n$ act by joint conjugation on $\Conjug_{W_n}(a_1|_{T_n}) \times \Conjug_{W_n}(a_2|_{T_n})$. Since $G_n$ is generated by $a_1|_{T_n}$ and $a_2|_{T_n}$, Lemma \ref{25SmallRigidLem2} means that the stabilizer of the pair $(a_1|_{T_n},a_2|_{T_n})$ has order~$2$. Thus the corresponding $W_n$-orbit has length $|W_n|/2 = 2^{2^n-2}$. By Lemma \ref{25SmallRigidLem1} this orbit therefore fills out all of $\Conjug_{W_n}(a_1|_{T_n}) \times \Conjug_{W_n}(a_2|_{T_n})$; in other words the action is transitive. Taking the inverse limit over $n$ implies the lemma by the same argument as in Lemma \ref{13ConjLimit}.
\end{Proof}

\begin{Thm}\label{25SmallRigid}
If $r=2$, the assertions in Theorem \ref{24SemiRigid} are also equivalent to:
\begin{enumerate}
\item[(d)] There exists $w\in W$ such that $b_1=wa_1w^{-1}$ and $b_2=wa_2w^{-1}$.
\end{enumerate}
\end{Thm}


\begin{Proof}
The equivalence of (d) with the assertion (b) in Theorem \ref{24SemiRigid} is precisely the content of Lemma \ref{25SmallRigidLem3}.
\end{Proof}


\subsection{Normalizer}
\label{25Normalizer}

In this subsection we determine the normalizer
\UseTheoremCounterForNextEquation
\begin{equation}\label{25NormDef}
N := \Norm_W(G).
\end{equation}
Recall from Theorem \ref{22Gab} that $K := H^{(1)}\cap\ldots\cap H^{(r)}$ is the closure of the commutator subgroup of $G$ and that we have a natural isomorphism $G_\ab = G/K \cong \BZ_2^r$. Let $\pi\colon G\onto G_\ab$ denote the natural projection, so that the elements $\pi(a_i)$ for $1\le i\le r$ form a basis of $G_\ab$ as a module over~$\BZ_2$.

\begin{Lem}\label{25Lem1}
For any $w\in N$ and any $1\le i\le r$ there exists $k_i\in\BZ_2^\times$ such that $\pi(wa_iw^{-1}) = \pi(a_i^{k_i})$.
\end{Lem}

\begin{Proof}
Since $wa_iw^{-1}$ lies in~$G$, by Proposition \ref{23WConj=GConjPower} there exists $k_i\in\BZ_2^\times$ such that $wa_iw^{-1}$ is conjugate to $a_i^{k_i}$ under~$G$. Thus $\pi(wa_iw^{-1}) = \pi(a_i^{k_i})$, as desired.
\end{Proof}

\begin{Lem}\label{25Lem2}
For any tuple $(k_i)\in(\BZ_2^\times)^r$ there exists $w\in N$ such that $\pi(wa_iw^{-1}) = \pi(a_i^{k_i})$ for each $1\le i\le r$.
\end{Lem}

\begin{Proof}
By Proposition \ref{13ConjPowers} the element $a_i^{k_i}$ is conjugate to $a_i$ under $W$ for each~$i$. Thus by Theorem \ref{24SemiRigid} there exist $w\in W$ and $x_1,\ldots,x_r\in G$ such that $a_i^{k_i}=(wx_i)a_i(wx_i)^{-1}$ for each~$i$ and $\Langle a_1^{k_1},\ldots,a_r^{k_r}\Rangle=wGw^{-1}$. But $\Langle a_1^{k_1},\ldots,a_r^{k_r}\Rangle = \Langle a_1,\ldots,a_r\Rangle = G$; hence $w$ lies in~$N$. Also $\pi(wa_iw^{-1}) = \pi((wx_i)a_i(wx_i)^{-1}) = \pi(a_i^{k_i})$ for each~$i$, as desired.
\end{Proof}

\medskip
We can actually construct explicit elements with the property in Lemma \ref{25Lem2}, as follows. 

\begin{Prop}\label{25Explicit}
For any tuple $(k_i)\in(\BZ_2^\times)^r$ write $k_i = 1+2\ell_i$ with $\ell_i\in\BZ_2$. Define elements $w_i\in W$ by the recursion relations $w_i = (w_{i-1},a_r^{\ell_{1-i}}w_{i-1})$
for all integers $i$ modulo~$r$. Then for all $i$ and $j$ modulo $r$ we have $w_ia_jw_i^{-1} = a_j^{k_{j-i}}$. In particular the element $w:=w_r$ has the property in Lemma \ref{25Lem2}.
\end{Prop}

\begin{Proof}
For all $i$ and $j$ set $y_{i,j} := w_ia_jw_i^{-1}a_j^{-k_{j-i}}$. Then for $j>1$ we have
\begin{eqnarray*}
y_{i,j} &=&
(w_{i-1},a_r^{\ell_{1-i}}w_{i-1})\,
(a_{j-1},1)\,
(w_{i-1},a_r^{\ell_{1-i}}w_{i-1})^{-1}
(a_{j-1},1)^{-k_{j-i}} \\
&=& (w_{i-1}a_{j-1}w_{i-1}^{-1}a_{j-1}^{-k_{j-i}},1) \\
&=& (y_{i-1,j-1},1),
\end{eqnarray*}
and for $j=1$ we have
\begin{eqnarray*}
y_{i,1} &=& w_ia_1w_i^{-1}a_1^{-1}(a_1^2)^{-\ell_{1-i}} \\
&=& (w_{i-1},a_r^{\ell_{1-i}}w_{i-1})\,
(a_r,1)\,\sigma\,
(w_{i-1},a_r^{\ell_{1-i}}w_{i-1})^{-1}
\sigma^{-1}(a_r,1)^{-1}
(a_r,a_r)^{-\ell_{1-i}} \\
&=& (w_{i-1},a_r^{\ell_{1-i}}w_{i-1})\,
(a_r,1)\,
(w_{i-1}^{-1}a_r^{-\ell_{1-i}},w_{i-1}^{-1})
(a_r^{-1},1)
(a_r^{-\ell_{1-i}},a_r^{-\ell_{1-i}}) \\
&=& (w_{i-1}a_rw_{i-1}^{-1}a_r^{-1-2\ell_{1-i}},1) \\
&=& (y_{i-1,r},1).
\end{eqnarray*}
These recursion relations for the $y_{i,j}$ satisfy the conditions of Proposition \ref{15RecTriv}; hence $y_{i,j}=1$ for all $i,j$.
\end{Proof}


\medskip
Now let $M$ denote the kernel of the natural homomorphism $N \to \Aut(G_\ab)$.

\begin{Lem}\label{25Lem3}
$M\subset (M\times M)\cdot G$.
\end{Lem}

\begin{Proof}
Take any element $w\in M$. If $w$ acts non-trivially on level $1$ of~$T$, we can replace it by $wa_1$, after which we have $w=(u,v)$ for certain $u,v\in W$. Then $(u,v)$ normalizes $G^1 = G\cap (W\times W)$; hence $u$ normalizes $\proj_1(G^1)=G$ by Proposition \ref{22G1Prop} (c), in other words we have $u\in N$. We claim that $u\in M$.

For this consider first any $1<i\le r$. Then, since $w\in M$, the commutator
$$[a_i^{-1},w] = (a_{i-1},1)^{-1} (u,v)\,(a_{i-1},1)\,(u,v)^{-1}
= (a_{i-1}^{-1}ua_{i-1}u^{-1},1) 
= ([a_{i-1}^{-1},u],1)$$
lies in~$K$. By Theorem \ref{22Gab} it therefore lies in $H^{(j)}$ for all $1\le j\le r$.
In particular it lies in $H^{(1)} = H^{(r)}\times H^{(r)}$ by Lemma \ref{22HiProp} (a); hence $[a_{i-1}^{-1},u]$ lies in $H^{(r)}$. Also, for any $1<j\le r$ the element $([a_{i-1}^{-1},u],1)$ lies in 
$H^{(j)}\cap (W\times W) \subset (H^{(j-1)}\times H^{(j-1)})\cdot\Langle(a_{j-1},a_{j-1}^{-1})\Rangle$ by Lemma \ref{22HiProp} (b), where the product is semidirect by Proposition \ref{22HiSemiDirect}. This implies that $[a_{i-1}^{-1},u]$ lies in $H^{(j-1)}$. Together this shows that $[a_{i-1}^{-1},u]$ lies in $H^{(1)}\cap\ldots\cap H^{(r)} = K$.

Next, since $w\in M$, the commutator
\begin{eqnarray*}
[a_1^{-1},w] &=& \sigma\,(a_r,1)^{-1} (u,v)\,(a_r,1)\,\sigma\,(u,v)^{-1} \\
&=& (1,a_r^{-1})\, (v,u)\,(1,a_r)\,(u^{-1},v^{-1}) \\
&=& (vu^{-1},a_r^{-1}ua_rv^{-1})
\end{eqnarray*}
also lies in~$K$, and hence in $H^{(j)}$ for all $1\le j\le r$.
In particular it lies in $H^{(1)} = H^{(r)}\times H^{(r)}$ by Lemma \ref{22HiProp} (a); hence 
\UseTheoremCounterForNextEquation
\begin{equation}\label{25Lem3Remem}
\hbox{$a_r^{-1}ua_rv^{-1}$ and $vu^{-1}$ lie in $H^{(r)}$,}
\end{equation}
and so their product $a_r^{-1}ua_ru^{-1} = [a_r^{-1},u]$ lies in $H^{(r)}$. Also, for any $1<j\le r$ the element $(vu^{-1},a_r^{-1}ua_rv^{-1})$ lies in $H^{(j)}\cap (W\times W) \subset (H^{(j-1)}\times H^{(j-1)})\cdot\Langle(a_{j-1},a_{j-1}^{-1})\Rangle$ by Lemma \ref{22HiProp} (b). Thus there exists $m\in\BZ_2$ such that 
$$\biggl\{\begin{array}{cl}
a_r^{-1}ua_rv^{-1} \!\!\!&\in H^{(j-1)}\cdot a_{j-1}^{-m} \qquad\hbox{and} \\[3pt]
vu^{-1}            \!\!\!&\in H^{(j-1)}\cdot a_{j-1}^m \,=\, a_{j-1}^m \cdot H^{(j-1)}.
\end{array}\biggr\}$$
Therefore their product $a_r^{-1}ua_ru^{-1} = [a_r^{-1},u]$ lies in $H^{(j-1)}$. Together this shows that $[a_r^{-1},u]$ lies in $H^{(1)}\cap\ldots\cap H^{(r)} = K$.

We have thus shown that the commutator $[a_i^{-1},u]$ lies in $K$ for every $1\le i\le r$. This means that $u\in M$, as claimed.

Finally, the same arguments applied to the element 
$$a_1^{-1}wa_1 = \sigma\,(a_r,1)^{-1} (u,v)\,(a_r,1)\,\sigma
= (1,a_r^{-1})\, (v,u)\,(1,a_r)
= (v,a_r^{-1}ua_r)$$
in place of~$w$ show that $v\in M$. Thus $w=(u,v)\in M\times M$, and we are done.
\end{Proof}

\begin{Lem}\label{25Lem4}
$M=G$.
\end{Lem}

\begin{Proof}
It suffices to prove that for any $n\ge0$ and any element $w\in M$ the restriction $w|_{T_n}$ lies in~$G_n$. This is trivial for $n=0$, so assume that $n>0$ and that the assertion holds universally for~$n-1$. We will then prove it for $n$ and any $w\in M$.

By Lemma \ref{25Lem3} we can write $w=(u,v)g$ with $u,v\in M$ and $g\in G$. By the induction hypothesis applied to $v$ there exists an element $y\in G$ such that $v|_{T_{n-1}} = y|_{T_{n-1}}$. By Proposition \ref{22G1Prop} (c) we can choose an element $x\in G$ such that $z := (x,y)\in G$. Then $w = (ux^{-1},vy^{-1})\,zg$ with $(vy^{-1})|_{T_{n-1}}=1$ and $zg\in G$. To prove the desired assertion we may replace $w$ by $w(zg)^{-1}$. Afterwards we have $w=(u,v)$ with new elements $u,v\in M$ such that $v|_{T_{n-1}}=1$.

Now recall from (\ref{25Lem3Remem}) that $vu^{-1}$ lies in $H^{(r)}$. Since $v|_{T_{n-1}}=1$, this implies that $u^{-1}|_{T_{n-1}}$ and hence $u|_{T_{n-1}}$ lies in $H_{n-1}^{(r)}$. Thus 
$w|_{T_n} = (u|_{T_{n-1}},v|_{T_{n-1}}) = (u|_{T_{n-1}},1)$ lies in $H_{n-1}^{(r)}\times 1$ and hence, by Lemma \ref{22HnSemiDirect} (b), in~$G_n$. This proves the desired assertion for~$n$. 
By induction it therefore follows for all $n\ge0$.
\end{Proof}

\begin{Thm}\label{25NormThm}
The natural homomorphism $N\to\Aut(G_\ab)$ induces an isomorphism
$$N/G \stackrel{\sim}{\longto} (\BZ_2^\times)^r \subset \Aut(\BZ_2^r) \cong \Aut(G_\ab).$$
\end{Thm}

\begin{Proof}
Lemma \ref{25Lem4} shows that the homomorphism factors through an injection $N/G \into\Aut(G_\ab)$, and Lemmas \ref{25Lem1} and \ref{25Lem2} describe its image.
\end{Proof}

\subsection{Odometers}
\label{28aOdometers}

In this subsection we study the set of odometers in~$G$, in particular concerning conjugacy. First we show that odometers are abundant in~$G$. By the \emph{proportion} of a measurable subset $X\subset G$ we mean the ratio $\mu(X)/\mu(G)$ for any Haar measure $\mu$ on~$G$.

\begin{Prop}\label{28ManyOdos}
\begin{enumerate}
\item[(a)] For any $W$\!-conjugates $b_i$ of $a_i$ and any permutation $\rho$ of $\{1,\ldots,r\}$ the product $b_{\rho1}\cdots b_{\rho r}$ is an odometer.
\item[(b)] The set of odometers in $G$ is open and closed in $G$ and has proportion $2^{-r}$.
\end{enumerate}
\end{Prop}

\begin{Proof}
Consider the continuous homomorphism $\psi:=(\sgn_i)_{i=1}^r\colon G\to \{\pm1\}^r$. Proposition \ref{21GenSigns} implies that $\psi$ is surjective and that 
$\psi^{-1}((-1,\ldots,-1))$ is the set of elements $x\in G$ satisfying $\sgn_i(x)=-1$ for all $i\ge1$. Thus by Proposition \ref{16OdoEquiv} it is the set of odometers in~$G$, and (b) follows. Proposition \ref{21GenSigns} also implies that this subset contains the element $b_{\rho1}\cdots b_{\rho r}$ in (a), which is therefore an odometer.
\end{Proof}

\medskip
In particular the (somewhat arbitrarily chosen) element 
\UseTheoremCounterForNextEquation
\begin{equation}\label{2OdoDef}
a_0 := a_1\cdots a_r
\end{equation}
is an odometer in~$G$. From
(\ref{2RecRels}) we deduce the recursion relation
\UseTheoremCounterForNextEquation
\begin{equation}\label{2OdoRecRel}
a_1^{-1}\,a_0\,a_1 
\ =\ a_2\cdots a_r\, a_1
\ =\ (a_1,1)\cdots(a_{r-1},1)\,(a_r,1)\,\sigma
\ =\ (a_0,1)\,\sigma.
\end{equation}

\begin{Thm}\label{28OdosConj}
\begin{enumerate}
\item[(a)] All odometers in $G$ are conjugate under~$N$.
\item[(b)] Any two odometers in $G$ with the same image in~$G_\ab$ are conjugate under~$G$.
\end{enumerate}
\end{Thm}

\begin{Proof}
Let $[[G,G]]$ denote the closure of the commutator subgroup of~$G$. Then Proposition \ref{22HiSemiDirect} and Theorem \ref{22Gab} show that $G$ is the disjoint union of the cosets $[[G,G]]\cdot a_1^{k_1}\cdots a_r^{k_r}$ for all $(k_1,\ldots,k_r) \in \BZ_2^r$. Moreover, by the same argument as in the proof of Proposition \ref{28ManyOdos} the set of odometers in $G$ is the union of those cosets with $(k_1,\ldots,k_r) \in (\BZ_2^\times)^r$. As $N$ permutes these cosets transitively by Theorem \ref{25NormThm}, assertion (a) reduces to assertion (b), and assertion (b) reduces to the special case $(k_1,\ldots,k_r)=(1,\ldots,1)$. In view of (\ref{2OdoDef}) we must therefore prove that any element of the coset $[[G,G]]\cdot a_0$ is conjugate to $a_0$ under~$G$.
As a preparation we show:

\begin{Lem}\label{28OdosConjLem}
Any element $c\in [[G,G]]\cdot a_0$ is conjugate under $G$ to an element of the form $(d,1)\,\sigma$ for some $d\in [[G,G]]\cdot a_0$.
\end{Lem}

\begin{Proof}
{}From (\ref{2OdoRecRel}) we deduce that
$$[[G,G]]\cdot a_0
= [[G,G]]\cdot a_1^{-1}\,a_0\,a_1 
= [[G,G]]\cdot (a_0,1)\,\sigma.$$
Also, by Theorem \ref{22Gab} and Lemma \ref{22HiProp} (a) we have $[[G,G]] \subset
H^{(1)} = H\times H$ with ${H=H^{(r)}}$. Taken together we therefore find that $c=(ha_0,h')\,\sigma$ for some ${h,h'\in H}$. Then $(h',1)$ lies in~$G$, so that we can replace $c$ by its conjugate 
$$(h',1)\,c\,(h',1)^{-1} = (h'ha_0,1)\,\sigma.$$
Thus we may without loss of generality assume that $c\in (W\times1)\cdot (a_0,1)\,\sigma$.
On the other hand Lemma \ref{22HiProp} (b) implies that $H^{(i)} \cap (W\times1) \subset H^{(i-1)}\times1$ for all $1<i\le r$. Thus from Theorem \ref{22Gab} we deduce that
$$[[G,G]] \cap (W\times1) 
= H^{(1)}\cap\ldots\cap H^{(r)} \cap (W\times1)
\subset \bigl(H^{(1)}\cap\ldots\cap H^{(r)}\bigr)\times1
= [[G,G]]\times 1,$$
and so $c\in \bigl([[G,G]]\times 1\bigr)\cdot (a_0,1)\,\sigma$, as desired.
\end{Proof}

\medskip
To prove Theorem \ref{28OdosConj}, by  Lemma \ref{13ConjLimit} it suffices to show for all $n\ge0$:
\begin{enumerate}
\item[($*_n$)] For any element $c\in [[G,G]]\cdot a_0$, the restriction $c|_{T_n}$ is conjugate to $a_0|_{T_n}$ under~$G_n$. 
\end{enumerate}
This is trivial for $n=0$, so assume that $n>0$ and that ($*_{n-1}$) is true. 
Consider any element $c\in [[G,G]]\cdot a_0$. Since we may replace $c$ by an arbitrary conjugate under~$G$, by Lemma \ref{28OdosConjLem} we may assume that $c=(d,1)\,\sigma$ for some $d\in [[G,G]]\cdot a_0$. By the induction hypothesis ($*_{n-1}$) there then exists an element $x\in G$ with $(xa_0x^{-1})|_{T_{n-1}} = d|_{T_{n-1}}$.
Next observe that by the definitions (\ref{22HDef}) and (\ref{2OdoDef}) of~$H$ and $a_0$ we have $Ha_r=Ha_0$. Thus Proposition \ref{22HSemiDirect} (a) implies that $G = H \cdot \Langle a_r\Rangle = H \cdot \Langle a_0\Rangle$. Therefore $x\in Ha_0^\ell$ for some $\ell\in\BZ_2$, and after replacing $x$ by $xa_0^{-\ell}$, which does not change the conjugate $xa_0x^{-1}$, we may without loss of generality assume that $x\in H$. Then $z := (x,x)\,a_1^{-1}$ is an element of~$G$, which by the recursion relation (\ref{2OdoRecRel}) satisfies
$$z\,a_0\,z^{-1}
= (x,x)\,a_1^{-1}\,a_0\,a_1\,(x,x)^{-1}
\ =\ (x,x)\,(a_0,1)\,\sigma\,(x,x)^{-1}
\ =\ (xa_0x^{-1},1)\,\sigma.$$
Thus 
$$(z\,a_0\,z^{-1})|_{T_n}
\ =\ \bigl((xa_0x^{-1})|_{T_{n-1}},1\bigr)\,\sigma
\ =\ (d|_{T_{n-1}},1)\,\sigma
\ =\ c|_{T_n},$$
proving the assertion ($*_n$). By induction it follows for all $n\ge0$.
\end{Proof}

\medskip
Now we can extend Theorem \ref{24SemiRigid} as follows:

\begin{Thm}\label{28OdoSemiRigid}
Consider any elements $b_1,\ldots,b_r\in W$ and set $b_0:=b_1\cdots b_r$. Then the following are equivalent:
\begin{enumerate}
\item[(a)] $\left\{\begin{array}{l}
\hbox{$b_0$ is conjugate to $(b_0,1)\,\sigma$ under $W$,} \\[3pt]
\hbox{$b_1$ is conjugate to $(b_r,1)\,\sigma$ under $W$, and} \\[3pt]
\hbox{$b_i$ is conjugate to $(b_{i-1},1)$ under $W$ for any $2\le i\le r$.}\\
\end{array}\right\}$
\item[(b)] $b_i$ is conjugate to $a_i$ under $W$ for any $0\le i\le r$.
\item[(c)] There exist $w\in W$ and $x_i\in G$ such that $b_i=(wx_i)a_i(wx_i)^{-1}$ for any $0\le i\le r$.
\end{enumerate}
Moreover, for any $w$ as in (c) we have $\Langle b_0,\ldots,b_r\Rangle=wGw^{-1}$, and the coset $wG$ is uniquely determined.
\end{Thm}

\begin{Proof}
Each of the assertions (a), (b), (c) in Theorem \ref{28OdoSemiRigid} is the corresponding assertion from Theorem \ref{24SemiRigid} augmented by a condition on~$b_0$. To prove their equivalence it thus suffices to show that the additional conditions on $b_0$ already follow from the equivalent assertions in Theorem \ref{24SemiRigid}. 
Assuming these, by \ref{24SemiRigid} (c) we can find $w\in W$ and $x_i\in G$ such that $b_i=(wx_i)a_i(wx_i)^{-1}$ for all $1\le i\le r$. Then $c := w^{-1}b_0w = (x_1a_1x_1^{-1})\cdots(x_ra_rx_r^{-1})$ is an element of $G$ with the same image in $G_\ab$ as $a_0=a_1\cdots a_r$ and is therefore an odometer. So by Theorem \ref{28OdosConj} there exists $x_0\in G$ with $c=x_0a_0x_0^{-1}$ and hence $b_0=(wx_0)a_0(wx_0)^{-1}$. In particular $b_0$ is conjugate to $a_0$ under~$W$. By Proposition \ref{16OdoEquiv} this also implies that $b_0$ is conjugate to $(b_0,1)\,\sigma$ under $W$. This proves all the additional conditions on $b_0$ in Theorem \ref{28OdoSemiRigid}, as desired.

Of the last sentence, the first half again follows from Lemma \ref{13ConjGen}. Finally consider two elements $w,w'\in W$ with the properties in (c). Then there exist $x_0,x_0'\in G$ with $b_0 = (wx_0)a_0(wx_0)^{-1} = (w'x_0')a_0(w'x_0')^{-1}$. The element $(wx_0)^{-1}(w'x_0')\in W$ then commutes with~$a_0$, and so by Proposition \ref{16OdoProp} (b) it lies in $\Langle a_0\Rangle$ and hence in~$G$. Therefore $w'G = w'x_0'G = wx_0G = wG$, showing that the coset $wG$ is unique. This finishes the proof of Theorem \ref{28OdoSemiRigid}.
(Note that in (c) we can always replace $w,x_i$ by $wy,y^{-1}x_i$ for an arbitrary element $y\in G$, so the uniqueness of the coset $wG$ is optimal.)
\end{Proof}

\medskip
We can also describe how the normalizer of any odometer in $G$ sits inside~$N$:

\begin{Prop}\label{28OdoNorm}
Consider any odometer $c\in G$ and any $\ell\in\BZ_2^\times$.
\begin{enumerate}
\item[(a)] There exists an element $w\in N$ with $wcw^{-1}=c^\ell$. 
\item[(b)] Every element $w\in W$ with $wcw^{-1}=c^\ell$ lies in~$N$.
\item[(c)] For any such $w$, the image of the coset $wG$ under the isomorphism $N/G \stackrel{\sim}{\longto} (\BZ_2^\times)^r$ of Theorem \ref{25NormThm} is $(\ell,\ldots,\ell)$.
\end{enumerate}
\end{Prop}

\begin{Proof}
Since $c^\ell$ is another odometer in~$G$, part (a) follows from Theorem \ref{28OdosConj} (a). Next consider two elements $w,w'\in W$ with $wcw^{-1}=w'cw^{\prime-1}=c^\ell$. Then $w^{-1}w'\in W$ commutes with~$c$. By Proposition \ref{16OdoProp} (b) it therefore lies in $\Langle c\Rangle$ and hence in~$G$. Thus $wG=w'G$. Since $w\in N$ for some choice of~$w$, it follows that $w\in N$ for every choice, proving (b).

Recall from the proof of Theorem \ref{28OdosConj} that the image of $c$ under the isomorphism $G_\ab\cong\BZ_2^r$ of Theorem \ref{22Gab} is an element of the form $(k_1,\ldots,k_r) \in (\BZ_2^\times)^r$. The image of $wcw^{-1}=c^\ell$ is then equal to $(\ell k_1,\ldots,\ell k_r)$; hence the image of $w$ under the isomorphism of Theorem \ref{25NormThm} is $(\ell,\ldots,\ell)$.
\end{Proof}


\subsection{Iterated monodromy groups}
\label{28Monodromy}

Now let $G^\geom\subset G^\arith\subset W$ be the geometric and arithmetic fundamental groups associated to a quadratic polynomial over~$k$, as described in Subsection \ref{17Monodromy}, and assume that we are in the periodic case \ref{17Cases} (b). Then by (\ref{17GgeomGensPol}) we have $G^\geom = \Langle b_{p_1},\ldots,b_{p_r}\Rangle$, where by Proposition \ref{17GenRec} the generators satisfy the recursion relations
\UseTheoremCounterForNextEquation
\begin{equation}\label{18bRecRels}
\left\{\begin{array}{l}
\hbox{$b_{p_1}$ is conjugate to $(b_{p_r},1)\,\sigma$ under~$W$, and}\\[3pt]
\hbox{$b_{p_i}$ is conjugate to $(b_{p_{i-1}},1)$ under~$W$ for all $2\le i\le r$.}
\end{array}\right\}
\end{equation}
By Theorem \ref{24SemiRigid} we deduce:

\begin{Thm}\label{28GgeomThm}
There exists $w\in W$ such that $G^\geom=wGw^{-1}$ and $b_{p_i}$ is conjugate to $wa_iw^{-1}$ under $G^\geom$ for every $1\le i\le r$.
\end{Thm}

For the following we change the identification of trees used in Subsection \ref{17Monodromy} by the automorphism~$w$, after which we have $G^\geom=G$, and $b_{p_i}$ is conjugate to $a_i$ under $G$ for every $1\le i\le r$.
Then $G^\arith$ is contained in the normalizer~$N$, and to describe it it suffices to describe the factor group $G^\arith/G \subset N/G$. More precisely we will determine the composite homomorphism
\UseTheoremCounterForNextEquation
\begin{equation}\label{28GalRep}
\bar\rho\colon \Gal(\bar k/k) \onto G^\arith/G \subset N/G \cong (\BZ_2^\times)^r
\end{equation}
obtained from (\ref{17FactorMonoRep}) and Theorem \ref{25NormThm}.

\begin{Thm}\label{28GarithThm}
\begin{enumerate}
\item[(a)] The homomorphism $\bar\rho$ is the cyclotomic character $\Gal(\bar k/k) \to \nobreak \BZ_2^\times$ followed by the diagonal embedding $\BZ_2^\times \into (\BZ_2^\times)^r$. 
\item[(b)] If $k$ is finitely generated over its prime field, then $G^\arith/G$ is a subgroup of finite index of $\diag(\BZ_2^\times) \subset (\BZ_2^\times)^r$.
\end{enumerate}
\end{Thm}

\begin{Proof}
For any $1\le i\le r$ let $I_{p_i} \subset D_{p_i} \subset \pi_1^\et(\BP^1_k \setminus P,x_0)$ denote the inertia group, respectively the decomposition group at a point above~$p_i$. Since $p_i$ is a $k$-rational point, we have a commutative diagram with exact rows
$$\xymatrix{
1 \ar[r] & I_{p_i} \ar@{^{ (}->}[d]
\ar[r] & D_{p_i} \ar@{^{ (}->}[d]
\ar[r] & \Gal(\bar k/k) \ar[r] \ar@{=}[d] &  1 \\
1 \ar[r] & \pi_1^\et(\BP^1_{\bar k} \setminus P,x_0) \ar@{->>}[d]_\rho
\ar[r] & \pi_1^\et(\BP^1_k \setminus P,x_0) \ar@{->>}[d]_\rho
\ar[r] & \Gal(\bar k/k) \ar[r] \ar@{->>}[d]_{\bar\rho} & 1 \\
1 \ar[r] & G
\ar[r] & G^\arith
\ar[r] & G^\arith/G \ar[r] &  1\rlap{.} \\}$$
Since the target of the homomorphism $\rho$ is a pro-$2$ group, the restriction $\rho|_{I_{p_i}}$ factors through the maximal pro-$2$ quotient of~$I_{p_i}$. As the characteristic of $k$ is different from~$2$, this quotient is the pro-$2$ tame inertia group above $p_i$ and hence isomorphic~$\BZ_2$. Moreover, $D_{p_i}/I_{p_i} \cong \Gal(\bar k/k)$ acts on it through the cyclotomic character. Thus $\Gal(\bar k/k)$ acts on the subgroup $\Langle\pi(b_{p_i})\Rangle \subset G_\ab$ through the cyclotomic character. But since $b_{p_i}$ is conjugate to $a_i$ under~$G$, we have $\pi(b_{p_i})=\pi(a_i)$, and by Theorem \ref{22Gab} these elements for all $i$ form a basis of $G_\ab$ as a module over~$\BZ_2$. Thus $\Gal(\bar k/k)$ acts on $G_\ab$ diagonally through the cyclotomic character, proving (a).

(An alternative argument uses only the image of the inertia group at~$\infty$, a generator of which is $b_\infty=(b_{p_1}\cdots b_{p_r})^{-1}$ by Subsection \ref{17Monodromy} and therefore an odometer. The homomorphism $D_\infty\to N/G$ can then be determined directly using Proposition \ref{28OdoNorm}. Compare Subsection \ref{38Monodromy}.)

For (b) let $\BF$ denote the prime field of~$k$ and consider the cyclotomic character $\Gal(\bar\BF/\BF) \allowbreak {\to \BZ_2^\times}$. In the case $\BF=\BQ$ this is surjective. Otherwise $\BF=\BF_\ell$ for an odd prime~$\ell$, and the cyclotomic character sends $\Frob_\ell$ to $\ell\in\BZ_2^\times$. Since $\ell\not=\pm1$, it topologically generates an open subgroup of~$\BZ_2^\times$. Thus in both cases the cyclotomic character has an open image. As $k$ is a finitely generated extension of~$\BF$, the natural homomorphism $\Gal(\bar k/k)\to\Gal(\bar\BF/\BF)$ also has an open image. Thus (b) follows from (a).
\end{Proof}


%
%
%
%
%
%
%

%
%

\newpage
\section{Strictly pre-periodic polynomial case}
\label{3Preperiodic}

Throughout this section we fix integers $r>s\ge1$. Using Proposition \ref{15RecRelsProp} we define elements $a_1,\ldots,a_r\in W$ by the recursion relations
\UseTheoremCounterForNextEquation
\begin{equation}\label{3RecRels}
\left\{\!\begin{array}{lll}
a_1      &\!\!\!=\, \sigma, & \\[3pt]
a_{s+1}  &\!\!\!=\, (a_s,a_r), & \\[3pt]
a_i      &\!\!\!=\, (a_{i-1},1) & \hbox{for $i\not=1,s+1$.}\\
\end{array}\!\right\}
\end{equation}
We will study the subgroup $G:= \Langle a_1,\ldots,a_r\Rangle\subset W$ and its images $G_n\subset W_n$ for all~$n$.
The recursion relations are only one choice encoding the qualitative recursive behavior of the generators described by the diagram
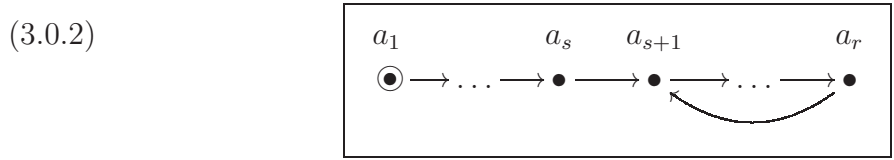
\UseTheoremCounterForNextEquation
\begin{equation}\label{3Diagram}
\fbox{\ $\xymatrix@R-24pt@C-10pt{
a_1 && a_s & a_{s+1} && \ a_r\ \mathstrut \\ 
\circbullet \ar[r]&\ldots\ar[r]&\bullet\ar[r]&\bullet\ar[r]&\ldots\ar[r]&\bullet\ar@/^16pt/[ll] \\
{\vphantom{X}} \mathstrut &&&\ \\
}$}\qquad\\
\end{equation}
where the vertices are now marked by the generators and the superfluous generator corresponding to the critical point $\infty$ is dropped: see Subsection \ref{17Monodromy}. In Subsection \ref{34Conjugacy} we will show that many other choices yield the same group $G$ up to conjugacy in~$W$.


\subsection{Basic properties}
\label{31Basic}

\begin{Prop}\label{31GenSigns}
For all $n\ge1$ and all $1\le i\le r$ we have
$$\sgn_n(a_i) \ =\ 
\left\{\!\begin{array}{rl}
-1 & \hbox{if $n=i\le s$,}\\[3pt]
-1 & \hbox{if $n\ge i>s$ and $n\equiv i\mod\,(r-s)$,}\\[3pt]
 1 & \hbox{otherwise.}
\end{array}\right.$$
\end{Prop}

\begin{Proof}
The definition of $a_i$ implies that $\sgn_1(a_1)=-1$ and $\sgn_1(a_i) = 1$ for all $2\le i\le r$; hence the formula holds for $n=1$. Assume that $n>1$ and that the formula holds for $n-1$ in place of~$n$. 
Then the formulas (\ref{15SignRec}) and (\ref{15SignSigma}) imply the recursion relations
$$\left\{\begin{array}{lll}
\sgn_n(a_1)     &\!\!\!=\, \sgn_n(\sigma) = 1 & \\[3pt]
\sgn_n(a_{s+1}) &\!\!\!=\, \sgn_{n-1}(a_s)\cdot\sgn_{n-1}(a_r) & \\[3pt]
\sgn_n(a_1)     &\!\!\!=\, \sgn_{n-1}(a_{i-1}) & \hbox{for $i\not=1,s+1$.}\\
\end{array}\!\right\}$$
Using these one directly verifies by case distinction that the formula is correct for~$n$. Thus the proposition follows by induction on~$n$.
\end{Proof}

\begin{Prop}\label{31GnAll}
For any $n\ge0$ we have $G_n=W_n$ if and only if $n\le r$.
\end{Prop}

\begin{Proof}
Direct consequence of Propositions \ref{15SignGnSurj} (a) and \ref{31GenSigns}.
\end{Proof}

\begin{Prop}\label{31LevelTrans}
The group $G$ acts transitively on the level $n$ of $T$ for every $n\ge0$.
\end{Prop}

\begin{Proof}
Consider the element $w := a_1\cdots a_r\in G$. Proposition \ref{31GenSigns} implies that $\sgn_n(w)=-1$ for all $n\ge1$; hence $w$ already acts transitively on every level by Proposition \ref{16OdoEquiv}.
\end{Proof}

\begin{Prop}\label{31GenOrder}
Each generator $a_i$ has order~$2$.
\end{Prop}

\begin{Proof}
By Proposition \ref{31GenSigns} we have $\sgn_i(a_i)=-1$ and hence $a_i\not=1$ for every $1\le i\le r$. On the other hand the elements $c_i := a_i^2$ satisfy the recursion relations
$$\left\{\!\begin{array}{lll}
c_1      &\!\!\!=\, \sigma^2 = 1, & \\[3pt]
c_{s+1}  &\!\!\!=\, (c_s,c_r), & \\[3pt]
c_i      &\!\!\!=\, (c_{i-1},1) & \hbox{for $i\not=1,s+1$.}\\
\end{array}\!\right\}$$
By Proposition \ref{15RecTriv} we therefore have $c_i=1$ for all~$i$. Thus each $a_i$ has precise order~$2$.
\end{Proof}

\begin{Lem}\label{31GenRes}
For any $1\le i\le r$ and $n\ge0$ we have $a_i|_{T_n}=1$ if and only if $n<i$.
\end{Lem}

\begin{Proof}
Direct consequence of the recursion relations (\ref{3RecRels}) and Proposition \ref{31GenSigns}.
\end{Proof}

\begin{Thm}\label{32Gab}
Let $G_\ab$ denote the maximal profinite abelian factor group of~$G$. Then the homomorphism $(\sgn_i)_{i=1}^r|_G$ induces an isomorphism $G_\ab \stackrel{\sim}{\to}\{\pm1\}^r$.
\end{Thm}

\begin{Proof}
As the elements $a_1,\ldots,a_r$ topologically generate~$G$, their images topologically generate the abelian group~$G_\ab$, and since the generators have order~$2$, we find that $G_\ab$ is a finite elementary abelian $2$-group generated by $r$ elements. On the other hand the continuous homomorphism $(\sgn_i)_{i=1}^r\colon G \to \{\pm1\}^r$ has an abelian group as target and therefore factors through~$G_\ab$. Since $\sgn_i(a_j) = (-1)^{\delta_{i,j}}$ (Kronecker delta) for all ${1\le i,j\le r}$ by Proposition \ref{31GenSigns}, the generators $a_1,\ldots,a_r$ map to a basis of $\{\pm1\}^r$. Together we deduce that $G_\ab \to\{\pm1\}^r$ is an isomorphism.
\end{Proof}

\begin{Rem}\label{31DihedralRem}
\rm Any two distinct indices $1\le i,j\le r$ yield distinct elements $a_i$ and $a_j$ of order two; hence the subgroup generated by them is dihedral (including infinite pro-dihedral) of order at least~$4$. Specifically, we have $\Langle a_i,a_j\Rangle = \Langle a_ia_j\Rangle \rtimes \langle a_i\rangle$ with $a_i(a_ia_j)a_i^{-1} = (a_ia_j)^{-1} = a_ja_i$. The order of this dihedral group is determined by the order of its cyclic part:
\end{Rem}

\begin{Lem}\label{32aiajSquares}
For any distinct $1\le i,j\le r$ we have 
$$(a_ia_j)^2 = 
\left\{\begin{array}{ll}
((a_{i-1}a_{j-1})^2,1) & \hbox{if $i,j\not=1$,} \\[3pt]
(a_{j-1},a_{j-1}) & \hbox{if $i=1$ and $j\not=s+1$,}\\[3pt]
(a_ra_s,a_sa_r) & \hbox{if $i=1$ and $j=s+1$.}
\end{array}\right.$$
\end{Lem}

\begin{Proof}
If both $i,j$ are $\not=1,s+1$, we have $(a_ia_j)^2 = \bigl( (a_{i-1},1)\,(a_{j-1},1) \bigr)^2 = ((a_{i-1}a_{j-1})^2,1)$.
If $i\not=1,s+1$ but $j=s+1$, we have $(a_ia_j)^2 = \bigl( (a_{i-1},1)\,(a_{j-1},a_r) \bigr)^2 = ((a_{i-1}a_{j-1})^2,a_r^2)$, which is again equal to $((a_{i-1}a_{j-1})^2,1)$ because $a_r^2=1$. A similar calculation applies when $i=s+1$ and $j\not=1,s+1$. This proves the first formula.

For $j\not=1,s+1$ we calculate $(a_1a_j)^2 = \sigma\,(a_{j-1},1)\,\sigma\,(a_{j-1},1) = (a_{j-1},a_{j-1})$, proving the second formula. Finally
$(a_1a_{s+1})^2 = \sigma\,(a_s,a_r)\,\sigma\,(a_s,a_r) = (a_ra_s,a_sa_r)$, proving the last.
\end{Proof}

\begin{Prop}\label{32aiajOrder}
For any $1\le i<j\le r$ we have 
$$\ord(a_ia_j) = \ord(a_ja_i) =
\left\{\begin{array}{ll}
4 & \hbox{if $j\not=i+s$,} \\[3pt]
8 & \hbox{if $j=i+s$ and $r\not=2s$,}\\[3pt]
\infty & \hbox{if $j=i+s$ and $r=2s$.}
\end{array}\right.$$
In particular we have 
$$\ord(a_sa_r) = \ord(a_ra_s) =
\left\{\begin{array}{ll}
4 & \hbox{if $r\not=2s$, \phantom{and $j=i+s$}} \\[3pt]
\infty & \hbox{if $r=2s$.}
\end{array}\right.$$
\end{Prop}

\begin{Proof}
The first formula in Lemma \ref{32aiajSquares} shows that $\ord(a_ia_j) = \ord(a_{i-1}a_{j-1})$ if $i,j>1$. By induction on~$i$ this implies that 
\UseTheoremCounterForNextEquation
\begin{equation}\label{32aiajInd}
\ord(a_ia_j) = \ord(a_1a_{j-i+1})\quad \hbox{for all $1\le i<j\le r$.}
\end{equation}
For $j\not=s+1$ the second formula in Lemma \ref{32aiajSquares} shows that 
$\ord(a_1a_j) = 2\cdot\ord(a_{j-1})=4$. This together with (\ref{32aiajInd}) implies the desired formula in the case $j\not=i+s$. 

The third formula in Lemma \ref{32aiajSquares} shows that $\ord(a_1a_{s+1}) = 2\cdot\ord(a_sa_r)$. If $r\not=2s$, by the case already proved this implies that $\ord(a_1a_{s+1}) = 2\cdot4=8$. With (\ref{32aiajInd}) we thus obtain the desired formula in the case $j=i+s$ and $r\not=2s$. On the other hand, if $r=2s$, using (\ref{32aiajInd}) we find that $\ord(a_1a_{s+1}) = 2\cdot\ord(a_sa_r) = 2\cdot\ord(a_1a_{r-s+1}) = 2\cdot\ord(a_1a_{s+1})$. Thus $\ord(a_1a_{s+1})=\infty$ in this case, and with (\ref{32aiajInd}) we deduce the desired formula in the last case.
\end{Proof}



\subsection{Useful subgroups}
\label{32Subgroups}

First we consider the subgroup of index $2$ which acts trivially on level $1$:
\UseTheoremCounterForNextEquation
\begin{equation}\label{23G1Def}
G^1\ :=\ G\cap(W\times W).
\end{equation}
Let $\proj_1$ and $\proj_2\colon W\times W\to W$ denote the two projections. We have the following self-similarity properties:

\begin{samepage}
\begin{Prop}\label{32G1Prop}
\begin{enumerate}
\item[(a)] $G\subset (G\times G)\rtimes \langle\sigma\rangle$.
\item[(b)] $G^1\subset G\times G$.
\item[(c)] $\proj_1(G^1)=\proj_2(G^1)=G$.
\end{enumerate}
\end{Prop}
\end{samepage}

\begin{Proof}
Essentially the same as that of Proposition \ref{22G1Prop}.
\end{Proof}

\medskip
Next we consider the following normal subgroup of~$G$:
\UseTheoremCounterForNextEquation
\begin{equation}\label{32HDef}
H\ :=\ 
\biggl\{\begin{array}{l}
\hbox{closure of the subgroup of $G$ generated} \\[3pt]
\hbox{by all $G$-conjugates of $a_i$ for all $i\not=s,r$.}
\end{array}\biggr\}_.
\end{equation}

\begin{Prop}\label{32HAlmostDirect}
We have:
\begin{enumerate}
\item[(a)] $G = H \cdot \Langle a_s,a_r\Rangle$.
\item[(b)] $G = (H\times H) \cdot \Langle \sigma, (a_s,a_r) \Rangle$.
\item[(c)] $G^1 = (H\times H) \cdot \Langle (a_s,a_r), (a_r,a_s) \Rangle$.
\end{enumerate}
\end{Prop}

\begin{Proof}
Assertion (a) is a direct consequence of the definition of~$H$. 

For (b) let $H'$ denote the closure of the subgroup of $G$ generated by all $G$-conjugates of $a_i$ for all $i\not=1,s+1$. Then by the same token we have $G = H' \cdot \Langle a_1,a_{s+1}\Rangle = H' \cdot \Langle \sigma, (a_s,a_r) \Rangle$. On the other hand, the generators of $H'$ are precisely the $G$-conjugates of $a_i = (a_{i-1},1)$ for all $2\le i\le r$ with $i\not=s+1$, or in other words the $G$-conjugates of $(a_j,1)$ for all $1\le j\le r$ with $j\not=s,r$. Proposition \ref{32G1Prop} implies that the conjugates of these elements under $G^1$ alone form the set of all $(x,1)$ where $x$ runs through all $G$-conjugates of $a_j$ for all $j\not=s,r$. These elements topologically generate the subgroup $H\times1\subset G$. Its conjugate by $\sigma$ is $1\times H$; hence together we find that $H'=H\times H$. This proves (b).

Finally, since $\sigma(a_s,a_r)\sigma^{-1} = (a_r,a_s)$, the intersection of $\Langle \sigma, (a_s,a_r) \Rangle$ with $W\times W$ is $\Langle (a_s,a_r), (a_r,a_s) \Rangle$; hence (b) implies (c).
\end{Proof}


\medskip
Next we consider the following closed normal subgroup of $G$:
\UseTheoremCounterForNextEquation
\begin{equation}\label{32KDef}
K\ :=\ H \cdot \Langle a_sa_r\Rangle.
\end{equation}
As the index of $\Langle a_sa_r\Rangle$ in $\Langle a_s,a_r\Rangle$ is~$2$, Proposition \ref{32HAlmostDirect} (a) implies that the index of $K$ in $G$ is at most~$2$. On the other hand, Proposition \ref{31GenSigns} shows that the signs $\sgn_s$ and $\sgn_r$ are trivial on all generators of~$H$, and the product $\sgn_s\sgn_r$ is also trivial on $a_sa_r$. Thus $K$ is contained in the subgroup $\Ker(\sgn_s\sgn_r|_G)$ of index $2$ of~$G$. Together this shows that 
\UseTheoremCounterForNextEquation
\begin{equation}\label{32KKer}
K=\Ker(\sgn_s\sgn_r|_G) \hbox{\ \ and has index $2$ in~$G$.}
\end{equation}

\begin{Lem}\label{32CommK}
The set $\{xy^{-1}\mid (x,y)\in G^1\}$ is equal to $K$.
\end{Lem}

\begin{Proof}
Using Remark \ref{31DihedralRem} we can write Proposition \ref{32HAlmostDirect} (c) in the form
$$G^1 = (H\times H) \cdot \Langle (a_sa_r,a_ra_s) \Rangle \cdot \langle (a_s,a_r) \rangle.$$
Thus an arbitrary element of $G^1$ has the form 
$$(x,y) = (h,h')\, (a_sa_r,a_ra_s)^\lambda (a_s,a_r)^\mu
= \bigl(h(a_sa_r)^\lambda a_s^\mu , h'(a_ra_s)^\lambda a_r^\mu\bigr)$$
for arbitrary $h,h'\in H$ and $\lambda\in\BZ_2$ and $\mu\in\{0,1\}$. Then
$$xy^{-1} = h(a_sa_r)^\lambda a_s^\mu a_r^{-\mu} (a_ra_s)^{-\lambda} h^{\prime-1}
= h(a_sa_r)^\lambda (a_sa_r)^\mu (a_sa_r)^\lambda h^{\prime-1}
= h(a_sa_r)^{2\lambda+\mu} h^{\prime-1}.$$
Here $2\lambda+\mu$ runs through all of~$\BZ_2$, and since $H$ is normal in $G$ and $h,h'$ are arbitrary elements of~$H$, the set of elements obtained is precisely $H \cdot \Langle a_sa_r\Rangle = K$, as desired.
\end{Proof}

\begin{Lem}\label{32SigmaConj}
For any $u\in W$ we have $(u,u^{-1}) \in G$ if and only if $u\in K$.
\end{Lem}

\begin{Proof}
If the element lies in $G$, it must already lie in~$G^1$. But by Proposition \ref{32HAlmostDirect} (c) we have
$$G^1 = (H\times H) \cdot \Langle (a_s,a_r), (a_r,a_s) \Rangle
 = (H\times H) \cdot \Langle (a_sa_r,a_ra_s)\Rangle\cdot \langle(a_s,a_r)\rangle.$$
Thus $(u,u^{-1})$ lies in $G^1$ if and only if there exist $\lambda\in\BZ_2$ and $\mu\in\{0,1\}$ such that 
\begin{eqnarray*}
(u,u^{-1}) &\in& (H\times H)\cdot(a_sa_r,a_ra_s)^\lambda\cdot(a_s,a_r)^\mu \\[5pt]
\Longleftrightarrow\quad\ \quad
\rlap{$\biggl\{\begin{array}{lll}
u &\!\!\!\in\ \;H (a_sa_r)^\lambda a_s^\mu & \hbox{and}\\[3pt]
u^{-1} &\!\!\!\in\ \;H (a_ra_s)^\lambda a_r^\mu &
\end{array}\biggr\}$} \qquad\,  && \\[5pt]
\Longleftrightarrow\qquad\quad\ 
u \ \,&\in& H (a_sa_r)^\lambda a_s^\mu \cap a_r^{-\mu}(a_ra_s)^{-\lambda}H.
\end{eqnarray*}
Since $\sgn_s$ is trivial on $H$ and~$a_r$, but not on~$a_s$, the last condition requires that $\lambda+\mu\equiv-\lambda$ modulo~$2$, in other words that $\mu=0$. Since already $(a_sa_r)^\lambda=(a_ra_s)^{-\lambda}$, the condition is therefore equivalent to $u \in H (a_sa_r)^\lambda$ for some $\lambda\in\BZ_2$, and thus to $u\in K$, as desired.
\end{Proof}

\medskip
Next we consider the following closed normal subgroup of $G$:
\UseTheoremCounterForNextEquation
\begin{equation}\label{32LDef}
L\ :=\ 
\biggl\{\begin{array}{l}
\hbox{closure of the subgroup of $G$ generated} \\[3pt]
\hbox{by all $G$-conjugates of $a_1=\sigma$.}
\end{array}\biggr\}_.
\end{equation}

\begin{Lem}\label{32LDescr}
Let $[[K,K]]$ denote the closure of the commutator subgroup of~$K$. Let $\anti\colon W\to W\times W$ denote the `antidiagonal' map $w\mapsto(w,w^{-1})$ (which is not a homomorphism!). Then
$$L =\Bigl( \bigl([[K,K]]\times[[K,K]]\bigr)\cdot\anti(K)\Bigr) \rtimes\langle\sigma\rangle.$$
\end{Lem}

\begin{Proof}
The definition of $L$ means equivalently that $L$ is the closure of the subgroup generated by $\sigma$ and the commutators $[g,\sigma]$ for all $g\in G$. Since $G=G^1\rtimes\langle\sigma\rangle$, these commutators are the same as those for all $g\in G^1$. For $(x,y)\in G^1$ we have 
$$(x,y)\,\sigma\,(x,y)^{-1}\sigma^{-1} = (xy^{-1},yx^{-1}) = \anti(xy^{-1}),$$
hence by Lemma \ref{32CommK} the set of these commutators is precisely $\anti(K)$.
For any $u,v\in K$ we calculate 
$$\anti(u)\,\anti(v)\,\anti(u^{-1}v^{-1})
= (u,u^{-1})\,(v,v^{-1})\,(u^{-1}v^{-1},vu)
= (uvu^{-1}v^{-1},1).$$
Thus the closure of the subgroup generated by $\anti(K)$ contains $[[K,K]]\times1$. Using this and its conjugate under $\sigma$ we deduce that $[[K,K]]\times[[K,K]] \subset L$.
Now the composite map
$$\anti\colon K\into K\times K \onto K/[[K,K]]\times K/ [[K,K]]$$
is a homomorphism, because the target group is abelian. Thus $\bigl([[K,K]]\times[[K,K]]\bigr)\cdot\anti(K)$ is already a subgroup, and since $K$ is closed and hence compact, it is also already closed. Thus this product is the closure of the subgroup generated by $\anti(K)$. Recombining this with $\sigma$ we deduce the desired formula in the lemma.
\end{Proof}

\begin{Lem}\label{32SquaresInL}
\begin{enumerate}
\item[(a)] For any $i,j\not=s,r$ we have $(a_ia_j)^2\in [[K,K]]$.
\item[(b)] For any $i,j\not=1,s+1$ we have $(a_ia_j)^2\in L$.
\end{enumerate}
\end{Lem}

\begin{Proof}
In (a) the assumption $i,j\not=s,r$ implies that $a_i$ and $a_j$ lie in $H$ and hence in~$K$. They also have order $2$, so that $(a_ia_j)^2 = [a_i,a_j] \in [[K,K]]$, as desired.

In (b) the assumption $i,j\not=1,s+1$ shows that $i-1,j-1\not=s,r$ and hence $(a_{i-1}a_{j-1})^2\in [[K,K]]$ by (a). Thus $(a_ia_j)^2 = \bigl((a_{i-1}a_{j-1})^2,1\bigr) \in [[K,K]] \times 1 \subset L$ by Lemma \ref{32LDescr}, proving (b).
\end{Proof}

\begin{Prop}\label{32GenCaseSquaresInH}
If $s\ge2$ and $r\ge4$, then $(a_sa_r)^2\in H$.
\end{Prop}

\begin{Proof}
If $s\ge2$ and $r>s+1$, Lemma \ref{32SquaresInL} (b) directly shows that $(a_sa_r)^2\in L$. But as $s>1$, the definitions of $L$ and $H$ imply that $L\subset H$, whence the proposition in this case.

If $r=s+1$, we first note that $(a_sa_r)^2 = ((a_{s-1}a_s)^2,1)$ by Lemma \ref{32aiajSquares}.
The assumption $r\ge4$ now implies that $s\ge3$ and therefore $s-1,s\not=1,s+1$. Thus $(a_{s-1}a_s)^2\in L$ by Lemma \ref{32SquaresInL} (b).
The fact that $s\ge3$ also implies that all $G$-conjugates of $a_2=(a_1,1)$ are among the generators of~$H$. Using Proposition \ref{32G1Prop} we find that the $G^1$-conjugates of $(a_1,1)$ are precisely the elements $(x,1)$ for all $G$-conjugates $x$ of $a_1=\sigma$. Thus by definition of $L$ implies that $L\times 1 \subset H$. Together this shows that
$(a_sa_r)^2 = \bigl((a_{s-1}a_s)^2,1\bigr) \in L\times 1 \subset H$, as desired.
\end{Proof}

\medskip
We will also need some results on finite levels. As before, for any subgroup $X\subset W$ we let $X_n$ denote its image in~$W_n$.

\begin{samepage}
\begin{Lem}\label{32IndexAtLeast4}
For any $n\ge r$ we have 
\begin{enumerate}
\item[(a)] $H_n \subset \Ker(\sgn_s|_{G_r})\cap\Ker(\sgn_r|_{G_r})$.
\item[(b)] $H_n \cap \langle a_s|_{T_n}, a_r|_{T_n}\rangle \subset \langle(a_sa_r)^2|_{T_n}\rangle$.
\item[(c)] $\bigl[ \langle a_s|_{T_n}, a_r|_{T_n} \rangle : \langle (a_sa_r)^2|_{T_n} \rangle \bigr] = 4$.
\end{enumerate}
\end{Lem}
\end{samepage}

\begin{Proof}
For $n\ge r$ both $\sgn_s$ and $\sgn_r$ factor through $W_n$ and are trivial on all generators of~$H$; hence they are trivial on~$H_n$. This proves (a). On the other hand $\sgn_s(a_s|_{T_n})=\sgn_r(a_r|_{T_n})=-1$ and $\sgn_s(a_r|_{T_n})=\sgn_r(a_s|_{T_n})=1$. Since $\langle a_s|_{T_n},a_r|_{T_n}\rangle$ is a dihedral group, this and the assertion (a) implies (b) and (c).
\end{Proof}

\begin{Prop}\label{32SmallCaseSquaresInAndNotInH}
If $s=2$ and $r=3$, then $(a_sa_r)^2|_{T_r}\in H_r$, whereas
$(a_sa_r)^2|_{T_n}\not\in H_n$ for all $n\ge r+1=4$.
\end{Prop}

\begin{Proof}
Since $s=2$ and $r=3$, Lemma \ref{32aiajSquares} implies that 
$(a_sa_r)^2 = (a_2a_3)^2 = ((a_1a_2)^2,1) = ((a_1,a_1),1) = ((\sigma,\sigma),1)$.
By the definitions of $H$ and $L$ the assumptions $s=2$ and $r=3$ also imply that $H=L$. 
Thus Lemma \ref{32LDescr} shows that 
$$H_n \cap (W_{n-1}\times W_{n-1}) 
\ =\ \bigl([K_{n-1},K_{n-1}]\times[K_{n-1},K_{n-1}]\bigr)\cdot\anti(K_{n-1})$$ 
where $[K_{n-1},K_{n-1}]$ denotes the commutator subgroup of $K_{n-1}$.
This in turn implies that 
$$\phantom{H_n \cap (W_{n-1}\times W_{n-1})} 
\llap{$H_n \cap (W_{n-1}\times 1)$} 
\ =\ \rlap{$[K_{n-1},K_{n-1}]\times 1.$}
\phantom{\bigl([K_{n-1},K_{n-1}]\times[K_{n-1},K_{n-1}]\bigr)\cdot\anti(K_{n-1})}$$ 
Therefore $(a_sa_r)^2|_{T_n}\in H_n$ if and only if $(\sigma,\sigma) \in [K_{n-1},K_{n-1}]$.

Suppose first that $n=3$. Then $n-1=2$, and $a_3|_{T_2}=1$. Thus $K_{n-1}=K_2$ contains the elements $a_1|_{T_2}=\sigma$ and $(a_2a_3)|_{T_2} = (a_1,1)|_{T_2} = (\sigma,1)$; in particular we have $K_2=W_2$.  It follows that $(\sigma,\sigma) = [\sigma,(\sigma,1)] \in [K_2,K_2]$, and so $(a_sa_r)^2|_{T_3}\in H_3$, as desired.

Suppose next that $n=4$. Then $[W_3:K_3]=2$ and a straightforward finite calculation, left to the gracious reader, shows that $[K_3,K_3]$ is the subgroup of order $8$ generated by
$$((\sigma,\sigma),1) \quad\hbox{and}\quad 
  (1,(\sigma,\sigma)) \quad\hbox{and}\quad
  ((\sigma,1)\sigma,(\sigma,1)\sigma).$$
It follows that $(\sigma,\sigma) \not\in [K_3,K_3]$, and so $(a_sa_r)^2|_{T_4}\not\in H_4$.

Finally, since $K_4$ is the image of $K_n$ in $W_4$ for any $n\ge4$, this also implies that $(a_sa_r)^2|_{T_n}\not\in H_n$ for any $n\ge4$, as desired.
\end{Proof}

\begin{Lem}\label{32Grplusone}
If $s\ge2$, then:
\begin{enumerate}
\item[(a)] $H_r = \Ker(\sgn_s|_{W_r})\cap\Ker(\sgn_r|_{W_r})$ and has index $4$ in~$W_r$.
\item[(b)] $H_r\times H_r$ is a normal subgroup of $W_{r+1} = (W_r\times W_r)\rtimes\langle\sigma\rangle$, and there is an isomorphism
$$\xymatrix@C-5pt@R-18pt{
W_{r+1}/(H_r\times H_r) \ar@{=}[r]^-\sim & 
\bigl(W_r/H_r \times W_r/H_r\bigr)\rtimes\langle\sigma\rangle
\ar[r]^-\sim & \Mat_{2\times2}(\{\pm1\})\rtimes\langle\sigma\rangle \rlap{,} \\
& (uH_r,vH_r)\cdot\sigma^\lambda \ar@{|->}[r] &
\binom{\;\sgn_s(u)\ \sgn_s(v)\;}{\sgn_r(u)\ \sgn_r(v)}\cdot \sigma^\lambda
\rlap{,} \\ }$$
where $\sigma$ acts on a $2\times2$-matrix by interchanging the columns.
\item[(c)] The image of $G_{r+1}/(H_r\times H_r)$ under this isomorphism is the subgroup
$$\textstyle \bigl\langle
\binom{-1\ \phantom{-}1\;}{\phantom{-}1\ -1\;},
\binom{\phantom{-}1\ -1\;}{-1\ \phantom{-}1\;}
\bigr\rangle
\rtimes\langle\sigma\rangle.$$
\item[(d)] The normalizer of this last subgroup in $\Mat_{2\times2}(\{\pm1\})\rtimes\langle\sigma\rangle$ is
$$\textstyle \bigl\langle
\binom{-1\ \phantom{-}1\;}{\phantom{-}1\ -1\;},
\binom{\phantom{-}1\ -1\;}{-1\ \phantom{-}1\;},
\binom{-1\ \ 1\;}{-1\ \ 1\;}
\bigr\rangle
\rtimes\langle\sigma\rangle.$$
In particular, the index of $G_{r+1}$ in its normalizer in $W_{r+1}$ is~$2$.
\end{enumerate}
\end{Lem}

\begin{Proof}
By Propositions \ref{32GenCaseSquaresInH} and \ref{32SmallCaseSquaresInAndNotInH} we have $(a_sa_r)^2|_{T_r}\in H_r$ whenever $s\ge2$. Thus the inclusion in Lemma \ref{32IndexAtLeast4} (b) is an equality. On the other hand, by Propositions \ref{31GnAll} and \ref{32HAlmostDirect} (a) we have $G_r=W_r = H_r \cdot \langle a_s|_{T_r},a_r|_{T_r}\rangle$. 
Combining this with the remaining information in Lemma \ref{32IndexAtLeast4} implies (a).

{}From (a) we readily obtain the isomorphism in~(b). Next Proposition \ref{32HAlmostDirect} implies that $G_{r+1} = G^1_{r+1}\rtimes\langle\sigma\rangle$ with $G^1_{r+1} = (H_r\times H_r) \cdot \langle (a_s|_{T_r},a_r|_{T_r}), (a_r|_{T_r},a_s|_{T_r}) \rangle$, which directly implies (c). The last part (d) is an easy direct calculation that is left to the reader.
\end{Proof}

\begin{Prop}\label{32MidCaseSquaresNotInH}
If $s=1$, then $(a_sa_r)^2|_{T_n}\not\in H_n$ for all $n\ge r$.
\end{Prop}

\begin{Proof}
If $r=2$, then $H=1$ by definition and the assertion holds trivially. So assume that $s\ge3$.
Then by definition $H$ is the closure of the subgroup generated by all $G$-conjugates of $a_i=(a_{i-1},1)$ for all $i\not=1,s+1,r$ and of $a_{s+1} = (a_s,a_r)$. In other words it is the closure of the subgroup generated by all $G$-conjugates of $(a_j,1)$ for all $j\not=s,r-1,r$ and of $(a_s,a_r) = (a_1,a_r)$. Thus if $H'$ denotes the closure of the subgroup generated by all $G$-conjugates of $a_i$ for all $i\not=r-1$, from Proposition \ref{32G1Prop} we deduce that $H\subset H'\times H'$. This in turn implies that $H_n\subset H'_{n-1}\times H'_{n-1}$.

Now $(a_sa_r)^2 = (a_1a_r)^2 = (a_{r-1},a_{r-1})$ by Lemma \ref{32aiajSquares}, because $r>2=s+1$. Thus $(a_sa_r)^2|_{T_n} = (a_{r-1}|_{T_{n-1}},a_{r-1}|_{T_{n-1}})$, and so to prove the proposition it suffices to show that $a_{r-1}|_{T_{n-1}} \not\in H'_{n-1}$.

For this recall from Proposition \ref{31GenSigns} that $\sgn_{r-1}(a_{r-1})=-1$, while $\sgn_{r-1}(a_i)=1$ for all $i\not=r-1$. Thus $\sgn_{r-1}$ is trivial on all the generators of $H'$ and hence on~$H'$. Moreover, the assumption $n\ge r$ implies that $n-1\ge r-1$; hence $\sgn_{r-1}$ factors through $W_{n-1}$. Thus the resulting function $W_{n-1} \to \{\pm1\}$ is identically $1$ on $H'_{n-1}$ but not on $a_{r-1}|_{T_{n-1}}$; hence $a_{r-1}|_{T_{n-1}} \not\in H'_{n-1}$, as desired.
\end{Proof}

\begin{Prop}\label{32Intersections}
For any $n\ge0$ we have
$$H_n \cap \langle a_s|_{T_n}, a_r|_{T_n} \rangle \ =\ 
\left\{\begin{array}{ll}
1 & \hbox{if $n<r$,} \\[3pt]
1 & \hbox{if $s=1$,} \\[3pt]
1 & \hbox{if $s=2$ and $r=3$ and $n\ge4$,} \\[3pt]
\langle(a_sa_r)^2|_{T_n}\rangle & \hbox{if $s=2$ and $r=3$ and $n=3$,} \\[3pt]
\langle(a_sa_r)^2|_{T_n}\rangle & \hbox{if $s\ge2$ and $r\ge4$ and $n\ge r$.}
\end{array}\right.$$
\end{Prop}

\begin{Proof}
If $n<s$ then $a_s|_{T_n}=a_r|_{T_n}=1$ and the assertion holds trivially. 
If $s\le n<r$ then $a_r|_{T_n}=1$, and so the question concerns $H_n \cap \langle a_s|_{T_n} \rangle$. But then $\sgn_s$ factors through $W_n$ and is trivial on all generators of $H$; hence $\sgn_s$ is trivial on~$H_n$. On the other hand $\sgn_s(a_s|_{T_n})=-1$, and $a_s|_{T_n}$ is an element of order~$2$; hence the intersection is trivial, as desired.

For the rest of the proof assume that $n\ge r$. Then by Lemma \ref{32IndexAtLeast4} (b) we have $H_n \cap \langle a_s|_{T_n}, a_r|_{T_n}\rangle \subset \langle(a_sa_r)^2|_{T_n}\rangle$.

If $s\ge2$ and $r\ge4$, then Proposition \ref{32GenCaseSquaresInH} implies that this inclusion is an equality, as desired. If $s=2$ and $r=3$ and $n=r=3$, the same follows from Proposition \ref{32SmallCaseSquaresInAndNotInH}. 

If $s=2$ and $r=3$ and $n\ge r+1=4$, then $(a_sa_r)^2|_{T_n}\not\in H_n$ by Proposition \ref{32SmallCaseSquaresInAndNotInH}. Since $(a_sa_r)^4=1$ by Proposition \ref{32aiajOrder} in this case, it follows that the intersection is trivial, as desired.

Finally, if $s=1$, then $(a_sa_r)^2|_{T_n}\not\in H_n$ by Proposition \ref{32MidCaseSquaresNotInH}. If in addition $r\ge3$, then $(a_sa_r)^4=1$ by Proposition \ref{32aiajOrder}; hence the intersection is trivial, as desired. If $r=2$, then already $H=1$ by construction and hence $H_n=1$ and the intersection is again trivial. This finishes the proof in all cases.
\end{Proof}


\subsection{Size}
\label{33Size}

\begin{Lem}\label{33Indices1}
For any $n\ge r\ge 3$ we have
$$[G_n:H_n] \ =\ 
\left\{\begin{array}{ll}
8 & \hbox{if $s=1$,} \\[3pt]
8 & \hbox{if $s=2$ and $r=3$ and $n\ge4$,} \\[3pt]
4 & \hbox{if $s=2$ and $r=3$ and $n=3$,} \\[3pt]
4 & \hbox{if $s\ge2$ and $r\ge4$.}
\end{array}\right.$$
\end{Lem}

\begin{Proof}
In all cases, Proposition \ref{32HAlmostDirect} (a) implies that $G_n = H_n \cdot \langle a_s|_{T_n}, a_r|_{T_n} \rangle$ and hence 
$$[G_n:H_n] \ =\ \bigl[ \langle a_s|_{T_n}, a_r|_{T_n} \rangle : H_n \cap \langle a_s|_{T_n}, a_r|_{T_n} \rangle \bigr].$$
If $s=1$, the assumption $r\ge3$ and Lemma \ref{32aiajSquares} show that $(a_sa_r)^2|_{T_n} = (a_{r-1}|_{T_{n-1}},a_{r-1}|_{T_{n-1}})$, where $a_{r-1}|_{T_{n-1}}$ has order $2$ by Lemma \ref{31GenRes}. Thus $\langle a_s|_{T_n}, a_r|_{T_n} \rangle$ has order $8$ in this case, and the desired assertion follows from the case $s=1$ of Proposition \ref{32Intersections}.

If $s=2$ and $r=3$, then Lemma \ref{32aiajSquares} implies that $(a_sa_r)^2|_{T_n} = ((a_1a_2)^2|_{T_{n-1}},1)$ and that $(a_1a_2)^2|_{T_{n-1}} = (a_1|_{T_{n-2}},a_1|_{T_{n-2}})$. Since $n\ge3$, the element $a_1|_{T_{n-2}} = \sigma|_{T_{n-2}}$ has order~$2$; hence so does $(a_sa_r)^2|_{T_n}$. Thus $\langle a_s|_{T_n}, a_r|_{T_n} \rangle$ has order $8$ in this case. Again the desired assertions follows from the corresponding cases of Proposition \ref{32Intersections}.

Finally, from Lemma \ref{32IndexAtLeast4} (c) we already know that $\bigl[ \langle a_s|_{T_n}, a_r|_{T_n} \rangle : \langle (a_sa_r)^2|_{T_n} \rangle \bigr] = \nobreak 4$ in all cases. Thus the last case is a direct consequence of the corresponding case in Proposition \ref{32Intersections}. 
\end{Proof}

\begin{Lem}\label{33Indices2}
For any $n\ge r$ we have
$$[G_{n+1}:H_n\times H_n] \ =\ 2\cdot [G_n:H_n].$$
\end{Lem}

\begin{Proof}
By Proposition \ref{32HAlmostDirect} (c) we have $G^1_{n+1} = (H_n\times H_n) \cdot \langle (a_s|_{T_n},a_r|_{T_n}), (a_r|_{T_n},a_s|_{T_n}) \rangle$. Here the second factor is a dihedral group which maps isomorphically to $\langle a_s|_{T_n}, a_r|_{T_n} \rangle$ under both projections.
{}From Lemma \ref{32IndexAtLeast4} (b) we already know that $H_n \cap \langle a_s|_{T_n}, a_r|_{T_n}\rangle$ is contained in $\langle(a_sa_r)^2|_{T_n}\rangle$.
Since moreover $\bigl((a_s,a_r)\, (a_r,a_s)\bigr)^2 = \bigl((a_sa_r)^2,(a_ra_s)^2\bigr)$, where the second entry is the inverse of the first, we deduce that $(H_n\times H_n) \cap \langle (a_s|_{T_n},a_r|_{T_n}), (a_r|_{T_n},a_s|_{T_n}) \rangle$ maps isomorphically to $H_n\cap \langle a_s|_{T_n}, a_r|_{T_n} \rangle$ under both projections. Thus 
\begin{eqnarray*}
[G_{n+1}^1:H_n\times H_n] 
&\!\!=\!\!& \bigl[ \langle (a_s|_{T_n},a_r|_{T_n}), (a_r|_{T_n},a_s|_{T_n}) \rangle : 
(H_n{\times}H_n) \cap \langle (a_s|_{T_n},a_r|_{T_n}), (a_r|_{T_n},a_s|_{T_n}) \rangle\bigr]\\
&\!\!=\!\!& \bigl[ \langle a_s|_{T_n},a_r|_{T_n} \rangle : 
H_n \cap \langle \langle a_s|_{T_n},a_r|_{T_n} \rangle \bigr]\\
&\!\!=\!\!& [G_n : H_n].
\end{eqnarray*}
Finally, since $n+1>0$, we have $[G_{n+1}:G_{n+1}^1]=2$, and the proposition follows.
\end{Proof}

\begin{Prop}\label{33GnOrder}
For all $n\ge0$ we have 
$$\log_2|G_n|\ =
\left\{\begin{array}{ll}
2^n-1                     & \hbox{if $n\le r$,} \\[3pt]
n+1                       & \hbox{if $s=1$ and $r=2$ and $n\ge2$,} \\[3pt]
2^n - 3\cdot 2^{n-r} + 2\ & \hbox{if $s=1$ and $r\ge3$ and $n\ge r$,} \\[3pt]
2^n - 5\cdot 2^{n-4} + 2  & \hbox{if $s=2$ and $r=3$ and $n\ge4$,} \\[3pt]
2^n - 2^{n-r+1} + 1       & \hbox{if $s\ge2$ and $r\ge4$ and $n\ge r$.}
\end{array}\right.$$
\end{Prop}

\begin{Proof}
The first case $n\le r$ results from Proposition \ref{31GnAll} and the formula (\ref{12WnOrder}). 

In the case $s=1$ and $r=2$ we have $G_n = \langle a_s|_{T_n}, a_r|_{T_n} \rangle$ and will calculate directly. Here $a_1a_2 = \sigma\,(a_1,a_2)$ is non-trivial on level $1$ and hence on $T_n$ for every $n\ge1$. This and the fact that $(a_1a_2)^2 = (a_2a_1,a_1a_2)$ implies that $\ord(a_1a_2|_{T_n}) = 2\cdot\ord(a_1a_2|_{T_{n-1}})$ for all $n\ge1$. By induction on $n$ we deduce that $\ord(a_1a_2|_{T_n}) = 2^n$ for all $n\ge0$. For all $n\ge r=2$ it follows that $G_n \cong D_{2^n}$ and hence $|G_n| = 2^{n+1}$, as desired.

In the other cases, for all $n\ge r$ Lemma \ref{33Indices2} implies that 
$$|G_{n+1}| \ =\  2\cdot [G_n:H_n] \cdot |H_n|^2
\ =\ |G_n|^2 \cdot 2\cdot [G_n:H_n]^{-1}.$$
From this and Lemma \ref{33Indices1} we deduce that
$$\log_2|G_{n+1}| \ =\ 2\cdot\log_2|G_n| - \left\{\begin{array}{ll}
2 & \hbox{if $s=1$,} \\[3pt]
2 & \hbox{if $s=2$ and $r=3$ and $n\ge4$,} \\[3pt]
1 & \hbox{if $s=2$ and $r=3$ and $n=3$,} \\[3pt]
1 & \hbox{if $s\ge2$ and $r\ge4$.}
\end{array}\right.$$
We use this recursion relation to deduce the remaining formulas by induction on~$n$.

In the case $s=1$ and $r\ge3$ observe first that $2^n - 3\cdot 2^{n-r} + 2 = 2^n-1$ for $n=r$; hence $\log_2|G_n| = 2^n - 3\cdot 2^{n-r} + 2$ for $n=r$. If the same formula holds for some $n\ge r$, the above recursion relation implies that $\log_2|G_{n+1}| = 2\cdot(2^n-3\cdot2^{n-r}+2)-2 = 2^{n+1}-3\cdot2^{n+1-r}+2$; hence the same formula holds for $n+1$ in place of~$n$. The formula therefore holds for all $n\ge r$, as desired.

In the case $s=2$ and $r=3$ we already know that $\log_2|G_3|=2^3-1=7$. Thus the recursion relation implies that $\log_2|G_4| = 2\cdot7-1 = 13$. This coincides with the value of $2^n-5\cdot2^{n-4}+2$ for $n=4$. If $\log_2|G_n| = 2^n-5\cdot2^{n-4}+2$ for some $n\ge r$, the above recursion relation implies that $\log_2|G_{n+1}| = 2\cdot(2^n-5\cdot2^{n-4}+2)-2 = 2^{n+1}-5\cdot2^{n+1-4}+2$; hence the same formula holds for $n+1$ in place of~$n$. The formula therefore holds for all $n\ge 4$, as desired.

Finally, in the case $s\ge2$ and $r\ge4$ we already have $\log_2|G_n| = 2^n-1 = 2^n - 2^{n-r+1} + 1$ in the case $n=r$. If $\log_2|G_n| = 2^n-2^{n-r+1}+1$ for some $n\ge r$, the above recursion relation implies that $\log_2|G_{n+1}| = 2\cdot(2^n-2^{n-r+1}+1)-1 = 2^{n+1}-2^{n+1-r+1}+1$; hence the same formula holds for $n+1$ in place of~$n$. The formula therefore holds for all $n\ge r$, as desired.
This finishes the proof in all cases.
\end{Proof}

\begin{Thm}\label{33Hausdorff}
The Hausdorff dimension of $G$ exists and is 
$$\left\{\begin{array}{ll}
0 & \hbox{if $s=1$ and $r=2$,} \\[3pt]
1-3\cdot2^{-r} & \hbox{if $s=1$ and $r\ge3$,} \\[3pt]
1-5\cdot2^{-4} & \hbox{if $s=2$ and $r=3$,} \\[3pt]
1-2^{1-r} & \hbox{if $s\ge2$ and $r\ge4$.}
\end{array}\right.$$
\end{Thm}

\begin{Proof}
By (\ref{12Hausdorff}) the Hausdorff dimension of $G$ is the limit of $\frac{\log_2|G_n|}{2^n-1}$ for $n\to\infty$. Thus the given formulas follow directly from Proposition \ref{33GnOrder}.
\end{Proof}


\subsection{Conjugacy of generators}
\label{34Conjugacy}

The opening remarks of Subsection \ref{24Conjugacy} also apply here to the recursion relations (\ref{3RecRels}). 

\begin{Thm}\label{34SemiRigid}
For any elements $b_1,\ldots,b_r\in W$ the following are equivalent:
\begin{enumerate}
\item[(a)] $\left\{\begin{array}{l}
\hbox{$b_1$ is conjugate to $\sigma$ under $W$,} \\[3pt]
\hbox{$b_{s+1}$ is conjugate to $(b_s,b_r)$ under $W$, and} \\[3pt]
\hbox{$b_i$ is conjugate to $(b_{i-1},1)$ under $W$ for all $i\not=1,s+1$.}\\
\end{array}\right\}$
\item[(b)] $b_i$ is conjugate to $a_i$ under $W$ for any $1\le i\le r$.
\item[(c)] There exist $w\in W$ and $x_i\in G$ such that $b_i=(wx_i)a_i(wx_i)^{-1}$ for any $1\le i\le r$.
\end{enumerate}
Moreover, for any $w$ as in (c) we have $\Langle b_1,\ldots,b_r\Rangle=wGw^{-1}$.
\end{Thm}

\begin{Proof}
The equivalence (a)$\Leftrightarrow$(b) is a special case of Proposition \ref{15RecConj=Conj}. The implication (c)$\Rightarrow$(b) is obvious, and the last statement follows from Lemma \ref{13ConjGen}. It remains to prove the implication (b)$\Rightarrow$(c). This follows by taking the limit over $n$ from the following assertion for all $n\ge0$: 
\begin{enumerate}
\item[($*_n$)] 
For any elements $b_i\in W_n$ conjugate to $a_i|_{T_n}$ under~$W_n$, there exist $w\in W_n$ and $x_i\in G_n$ such that $b_i=(wx_i)(a_i|_{T_n})(wx_i)^{-1}$ for each~$i$. 
\end{enumerate}
This is trivial for $n=0$, so assume that $n>0$ and that ($*_{n-1}$) is true. Take elements $b_i\in W_n$ conjugate to $a_i|_{T_n}$ under~$W_n$.
Since $b_1$ is conjugate to $a_1|_{T_n} = \sigma|_{T_n}$, after conjugating everything by the same element of $W_n$ we may assume that $b_1 = \sigma|_{T_n}$.

For each $i\not=1,s+1$ the element $b_i$ is conjugate to $a_i|_{T_n} = (a_{i-1}|_{T_{n-1}},1)$ under~$W_n$. Thus it is conjugate to $(a_{i-1}|_{T_{n-1}},1)$ or $(1,a_{i-1}|_{T_{n-1}})$ under the subgroup $W_{n-1}\times W_{n-1}$. In the second case we can replace $b_i$ by its conjugate $b_1b_ib_1^{-1}$, which, as in the proof of Theorem \ref{24SemiRigid}, does not change the assertion to be proved. Thus without loss of generality we may assume that $b_i$ is conjugate to $(a_{i-1}|_{T_{n-1}},1)$ under $W_{n-1}\times W_{n-1}$. This means that $b_i = (c_{i-1},1)$ for some $c_{i-1}\in W_{n-1}$ which is conjugate to $a_{i-1}|_{T_{n-1}}$ under $W_{n-1}$.

Likewise the element $b_{s+1}$ is conjugate to $a_{s+1}|_{T_n} = (a_s|_{T_{n-1}},a_r|_{T_{n-1}})$ under~$W_n$. By the same argument as above, after possibly replacing $b_{s+1}$ by its conjugate $b_1b_{s+1}b_1^{-1}$ we may without loss of generality assume that $b_{s+1}$ is conjugate to $(a_s|_{T_{n-1}},a_r|_{T_{n-1}})$ \hbox{under} ${W_{n-1}\times W_{n-1}}$. This means that $b_{s+1} = (c_s,c_r)$, where $c_s\in W_{n-1}$ is conjugate to $a_s|_{T_{n-1}}$ under $W_{n-1}$, and $c_r\in W_{n-1}$ is conjugate to $a_r|_{T_{n-1}}$ under $W_{n-1}$.

Thus for every $1\le i\le r$ we have found an element $c_i \in W_{n-1}$ which is conjugate to $a_i|_{T_{n-1}}$ under $W_{n-1}$. By the induction hypothesis ($*_{n-1}$) there therefore exist $u\in W_{n-1}$ and $x_1,\ldots,x_r\in G_{n-1}$ such that $c_i=(ux_i)(a_i|_{T_{n-1}})(ux_i)^{-1}$ for every~$i$. 

Now observe that the definition (\ref{32KDef}) of $K$ implies that $G = K\cdot\langle a_r\rangle$. Using this and Lemma \ref{32CommK}, we can write the element $x_s^{-1}x_r\in G_n$ in the form $x_s^{-1}x_r = xy^{-1}(a_r^\nu|_{T_{n-1}})$ for some $(x,y)\in G^1_n$ and $\nu\in\{0,1\}$. Set $v := ux_sx \in W_{n-1}$, so that $ux_s = vx^{-1}$ and $ux_r = ux_sxy^{-1}(a_r^\nu|_{T_{n-1}}) = vy^{-1}(a_r^\nu|_{T_{n-1}})$. Then
\begin{eqnarray*}
b_{s+1}
&\!\!=\!\!& \bigl(\ \,(ux_s)(a_s|_{T_{n-1}})(ux_s)^{-1}\ \,,\,\ (ux_r)(a_r|_{T_{n-1}})(ux_r)^{-1} \,\ \bigr) \\
&\!\!=\!\!& \bigl( (vx^{-1})(a_s|_{T_{n-1}})(vx^{-1})^{-1} , (vy^{-1})(a_r|_{T_{n-1}})(vy^{-1})^{-1} \bigr).
\end{eqnarray*}
Thus with $w := (v,v) \in W_n$ and $z_{s+1} := (x^{-1},y^{-1}) \in G^1_n \subset G_n$ we find that 
$$b_{s+1}
\ =\ (wz_{s+1})\,(a_s|_{T_{n-1}},a_r|_{T_{n-1}})(wz_{s+1})^{-1}
\ =\ (wz_{s+1})\,(a_{s+1}|T_n)\,(wz_{s+1})^{-1}.$$
Also, the element $w = (v,v)$ commutes with $\sigma|_{T_n}$; hence with $z_1 := 1\in G_n$ we have
$$b_1 = \sigma|_{T_n} = (wz_1)\,(a_1|T_n)\,(wz_1)^{-1}.$$
Furthermore, for any $i\not=1,s+1$ we have $ux_{i-1} = vx^{-1}x_s^{-1}x_{i-1}$ with $x^{-1}x_s^{-1}x_{i-1} \in \nobreak G_{n-1}$. By Proposition \ref{32G1Prop} we can therefore find an element $y_{i-1}\in G_{n-1}$ such that $z_i := (x^{-1}x_s^{-1}x_{i-1},y_{i-1})$ lies in~$G_n$. Then
\begin{eqnarray*}
b_i
&\!\!=\!\!& \bigl(\qquad\ \ \ (ux_{i-1})(a_{i-1}|_{T_{n-1}})(ux_{i-1})^{-1}\qquad\ \ \ ,\qquad\quad\: 1 \qquad\quad\: \bigr) \\
&\!\!=\!\!& \bigl( (vx^{-1}x_s^{-1}x_{i-1})(a_{i-1}|_{T_{n-1}})(vx^{-1}x_s^{-1}x_{i-1})^{-1} , (vy_{i-1})(vy_{i-1})^{-1} \bigr) \\
&\!\!=\!\!& (wz_i)\,(a_{i-1}|_{T_{n-1}},1)\,(wz_i)^{-1} \\
&\!\!=\!\!& (wz_i)\,(a_i|_{T_n})\,(wz_i)^{-1}.
\end{eqnarray*}
Altogether we have thus found an element $w\in W_n$ and elements $z_i\in G_n$ such that $b_i = (wz_i) (a_i|_{T_n}) (wz_i)^{-1}$ for all $1\le i\le r$. This proves ($*_n$), and so we are done by induction.
\end{Proof}

\begin{Prop}\label{33WConjInGIsGConj}
Consider any $1\le i\le r$. Assume that $i\le s$ or $s=1$. Then any element $b_i\in G$ that is conjugate to $a_i$ under~$W$ is conjugate to $a_i$ under~$G$.
\end{Prop}

The author does not know whether the proposition holds without the extra assumptions on $i$ and $s$.
\bigskip

\begin{Proof}
First we consider the case $i=1$. An arbitrary $W$-conjugate of $a_1=\sigma$ has the form $b_1=(u,v)\,\sigma\,(u,v)^{-1} = (uv^{-1},vu^{-1})\,\sigma$ for $u,v\in W$. Thus it lies in $G$ if and only if $(uv^{-1},vu^{-1})$ lies in~$G$. By Lemmas \ref{32SigmaConj} and Lemma \ref{32CommK} this is equivalent to $uv^{-1}=xy^{-1}$ for some $(x,y)\in G^1$. But then $b_1$ is equal to $(xy^{-1},yx^{-1}) = (x,y)\,\sigma\,(x,y)^{-1}$ and hence conjugate to $\sigma$ under~$G$, as desired.

Next we consider the case $2\le i\le s$. By induction on $i$ we may assume that the proposition holds for $i-1$ in place of~$i$. Take any element $w\in W$ such that $b_i := wa_iw^{-1}$ lies in~$G$. If $w$ acts non-trivially on level~$1$, we replace $b_i$ by its $G$-conjugate $a_1b_ia_1^{-1}$ and $w$ by $a_1w$, which does not change the desired assertion. Afterwards we have $w=(u,v)$ for some elements $u,v\in W$. Then the recursion relation yields
$$b_i = wa_iw^{-1}
= (u,v)\,(a_{i-1},1)\,(u,v)^{-1}
= (ua_{i-1}u^{-1},1).$$
Since $b_i$ lies in~$G$, Proposition \ref{32G1Prop} (b) implies that $b_{i-1} := ua_{i-1}u^{-1}$ lies in~$G$. By the induction hypothesis we therefore have $b_{i-1} = xa_{i-1}x^{-1}$
for some element $x\in G$. By Proposition \ref{32G1Prop} (c)  we can choose an element $y\in G$ for which $z := (x,y)$ lies in~$G$. Then
$$b_i = (b_{i-1},1)
= (xa_{i-1}x^{-1},1)
= (x,y)\,(a_{i-1},1)\,(x,y)^{-1}
= za_iz^{-1},$$
which proves that $b_i$ is conjugate to $a_i$ under~$G$, as desired.

\medskip
The remaining cases require an infinite recursion, which we realize as an induction over finite levels of the tree~$T$. By Lemma \ref{13ConjLimit} it suffices to show that for any $s+1\le i\le r$ and $n\ge0$:
\begin{enumerate}
\item[($*_{i,n}$)] Any element $b_i\in G_n$ that is conjugate to $a_i|_{T_n}$ under~$W_n$ is conjugate to $a_i|_{T_n}$ under~$G_n$. 
\end{enumerate}
This is trivial if $n=0$, so assume that $n>0$. For $s+1<i\le r$, the same argument as in the induction step above shows that ($*_{i-1,n-1}$) implies ($*_{i,n}$), as in the proof of Proposition \ref{23WConj=GConjPower}. So it remains to show that ($*_{r,n-1}$) implies ($*_{s+1,n}$). This is where the extra assumption $s=1$ will be used, but we will see that place more clearly if we keep $s$ arbitary until then. Consider an element $b_{s+1}\in G_n$ and an element $w\in W_n$ such that $b_{s+1} = w\,(a_{s+1}|_{T_n})\,w^{-1}$. 

First, if $w$ acts non-trivially on level~$1$, we replace $w$ by $(a_1|_{T_n})\,w$ and $b_{s+1}$ by its $G_n$-conjugate $(a_1|_{T_n})\,b_{s+1}\,(a_1|_{T_n})^{-1}$, which does not change the desired assertion. Afterwards we have $w=(u,v)$ for some elements $u,v\in W_{n-1}$. Next the recursion relation yields
$$b_{s+1} = w\,(a_{s+1}|_{T_n})\,w^{-1}
= (u,v)\,(a_s|_{T_{n-1}},a_r|_{T_{n-1}})\,(u,v)^{-1}
= \bigl(u(a_s|_{T_{n-1}})u^{-1},v(a_r|_{T_{n-1}})v^{-1}\bigr).$$
Since $b_{s+1}$ lies in~$G_n$, Proposition \ref{32G1Prop} (b) shows that $v(a_r|_{T_{n-1}})v^{-1}$ lies in~$G_{n-1}$. By the hypothesis ($*_{r,n-1}$) we therefore have $v(a_r|_{T_{n-1}})v^{-1} = y(a_r|_{T_{n-1}})y^{-1}$ for some $y\in G_{n-1}$. By Proposition \ref{32G1Prop} (c) we can choose an element $x\in G_{n-1}$ for which $z := (x,y)$ lies in~$G_n$. After replacing $b_{s+1}$ by its $G_n$-conjugate $z^{-1}b_{s+1}z$ and $w$ by $z^{-1}w$, and hence $u$ by $x^{-1}u$, which does not change the desired assertion, we may therefore assume that 
$$b_{s+1} = \bigl(u(a_s|_{T_{n-1}})u^{-1},a_r|_{T_{n-1}}\bigr).$$
Next, since $b_{s+1}$ lies in~$G_n$, Proposition \ref{32G1Prop} (b) shows that $u(a_s|_{T_{n-1}})u^{-1}$ lies in~$G_{n-1}$. By the induction hypothesis we can therefore find a new element $x\in G_{n-1}$ such that $u(a_s|_{T_{n-1}})u^{-1} = x(a_s|_{T_{n-1}})x^{-1}$.
By Proposition \ref{32HAlmostDirect} we can write this in the form $x = h ((a_sa_r)^\lambda a_s^\mu)|_{T_{n-1}}$ for some $h\in H_{n-1}$ and $\lambda\in\BZ$ and $\mu\in\{0,1\}$. Then
\begin{eqnarray*}
b_{s+1} 
&\!\!=\!\!& \bigl( h ((a_sa_r)^\lambda a_s^\mu a_s a_s^{-\mu}(a_sa_r)^{-\lambda})|_{T_{n-1}} h^{-1} ,a_r|_{T_{n-1}}\bigr) \\
&\!\!=\!\!& \bigl( h ((a_sa_r)^{2\lambda} a_s)|_{T_{n-1}} h^{-1} ,a_r|_{T_{n-1}}\bigr) \\
&\!\!\in\!\!& (H_{n-1}\times 1)\cdot
\bigl(((a_sa_r)^{2\lambda})|_{T_{n-1}},1\bigr)
\cdot(a_s,a_r)|_{T_n}.
\end{eqnarray*}
The fact that $b_{s+1}$ lies in~$G_n$ thus implies that $\bigl(((a_sa_r)^{2\lambda})|_{T_{n-1}},1\bigr)$ lies in~$G_n$, and hence in~$G^1_n$. 
Now, and only now,
we use the assumption $s=1$. By it and Proposition \ref{32Intersections}, we then have $H_{n-1} \cap \langle a_s|_{T_{n-1}}, a_r|_{T_{n-1}} \rangle = 1$. By Proposition \ref{32HAlmostDirect} (c) this implies that
$$G^1_n = (H_{n-1}\times H_{n-1}) \rtimes \langle (a_s|_{T_{n-1}},a_r|_{T_{n-1}}), (a_r|_{T_{n-1}},a_s|_{T_{n-1}}) \rangle.$$
The fact that $\bigl(((a_sa_r)^{2\lambda})|_{T_{n-1}},1\bigr)$ lies in $G_n$ is therefore equivalent to $((a_sa_r)^{2\lambda})|_{T_{n-1}}=1$. 
By the above calculation this means that 
$$b_{s+1} = \bigl(h(a_s|_{T_{n-1}})h^{-1},a_r|_{T_{n-1}}\bigr)
= (h,1)\,\bigl(a_s|_{T_{n-1}},a_r|_{T_{n-1}}\bigr)\,(h,1)^{-1}
= (h,1)\,(a_{s+1}|_{T_n})\,(h,1)^{-1}\!.$$
Since $(h,1)\in H_{n-1}\times1$ lies in $G_n$ by Proposition \ref{32HAlmostDirect} (b), this proves that $b_{s+1}$ is conjugate to $a_{s+1}|_{T_n}$ under~$G_n$, as desired.
Thus finishes the induction step and thus proves ($*_{i,n}$) for all $i$ and~$n$.
\end{Proof}


\subsection{Small cases}
\label{37Small}

In the smallest possible case
\UseTheoremCounterForNextEquation
\begin{equation}\label{35dAss}
\hbox{$s=1$\ \ and\ \ $r=2$}
\end{equation}
we have $G=\Langle a_1,\,a_2\Rangle$, where both generators have order two. Thus $G = \Langle a_1a_2\Rangle \rtimes \langle a_1\rangle \cong \BZ_2\rtimes\{\pm1\}$ is infinite dihedral by Proposition \ref{32aiajOrder}. Moreover $a_0 := a_1a_2$ is an odometer by Proposition \ref{38ManyOdos} below. Let $N$ denote the normalizer of $\Langle a_0\Rangle$ in~$W$.
Then by Proposition \ref{16OdoProp} we have $N\cong\BZ_2\rtimes\BZ_2^\times$, which implies that $N$ is also the normalizer of~$G$. In principle this provides all the information that we will need in this case. But we can also describe $N$ explicitly by the same method as in Proposition \ref{16OdoNormExpl}:

\begin{Prop}\label{37OdoNormExpl}
For any $k\in\BZ_2^\times$ consider the elements $v_k, w_k\in W$ defined uniquely~by
$$\biggl\{\begin{array}{rl}
v_k \!\!\!& = (a_0^{\smash{\frac{1-k}{2}}}v_k,v_k), \\[3pt]
w_k \!\!\!& = \ a_0^{\smash{\frac{1-k}{2}}}v_k,
\end{array}\biggr\}$$
where the first one requires Proposition \ref{15RecRelsProp}. Then
\begin{enumerate}
\item[(a)] For any $k\in\BZ_2^\times$ we have $w_ka_0w_k^{-1} = a_0^k$.
\item[(b)] For any $k,k'\in\BZ_2^\times$ we have $w_kw_{k'}=w_{kk'}$.
\item[(c)] There is an isomorphism $\BZ_2\rtimes\BZ_2^\times \stackrel{\sim}{\longto} N$, $(i,k)\mapsto a_0^iw_k$.
\end{enumerate}
\end{Prop}

\begin{Proof}
For (a) we abbreviate $\ell := \frac{k-1}{2} \in \BZ_2$ and $v := v_k = (a_0^{-\ell}v,v)$. Since 
$$a_0^2\ =\ a_1a_2a_1a_2 \ =\ \sigma\,(a_1,a_2)\,\sigma\,(a_1,a_2) 
\ =\ (a_2a_1,a_1a_2) \ =\ (a_0^{-1},a_0),$$
we have $a_0^{-k} = a_0^{-1-2\ell} = a_0^{-1}\,(a_0^\ell,a_0^{-\ell})$, as well as
\UseTheoremCounterForNextEquation
\begin{eqnarray}
\label{37OdoNormExpl1}
va_1v^{-1} &\!\!\!=\!\!\!& (a_0^{-\ell}v,v)\,\sigma\,(v^{-1}a_0^\ell,v^{-1})
\,=\, (a_0^{-\ell},a_0^\ell)\,\sigma
\,=\, a_0^{2\ell}a_1, \quad\hbox{and} \\
\nonumber
va_2v^{-1} &\!\!=\!\!& va_1v^{-1} va_1a_2v^{-1}\ =\ a_0^{2\ell}a_1 va_0v^{-1}.
\end{eqnarray}
Therefore
\begin{eqnarray*}
t \ :=\ v a_0 v^{-1} a_0^{-k}
&\!\!=\!\!& (a_0^{-\ell}v,v)\, \sigma\,(a_1,a_2)\, (v^{-1}a_0^\ell,v^{-1})\, 
(a_1^{-1},a_2^{-1})\,\sigma^{-1}\, (a_0^\ell,a_0^{-\ell}) \\
&\!\!=\!\!& \bigl(a_0^{-\ell}va_2v^{-1}a_2^{-1}a_0^\ell \,,\, 
va_1v^{-1}a_0^\ell a_1^{-1}a_0^{-\ell} \bigr) \\
&\!\!=\!\!& \bigl(a_0^{-\ell}a_0^{2\ell}a_1 va_0v^{-1}a_2^{-1}a_0^\ell \,,\, 
a_0^{2\ell}a_1a_0^\ell a_1^{-1}a_0^{-\ell} \bigr) \\
&\!\!=\!\!& \bigl(a_0^\ell a_1 ta_0^ka_2^{-1}a_0^\ell \,,\, 
a_0^{2\ell}a_0^{-\ell}a_0^{-\ell} \bigr)\\
&\!\!=\!\!& \bigl(a_0^\ell a_1 ta_1a_2a_0^{2\ell}a_2^{-1}a_0^\ell \,,\, 1 \bigr) \\
&\!\!=\!\!& \bigl(a_0^\ell a_1 ta_1a_0^{-2\ell}a_0^\ell \,,\, 1 \bigr) \\
&\!\!=\!\!& (a_0^\ell a_1 ta_1^{-1}a_0^{-\ell} \,,\, 1 ).
\end{eqnarray*}
By Proposition \ref{15RecTriv} this implies that $t=1$. Thus $v a_0 v^{-1} = a_0^k$, and therefore $w_ka_0w_k^{-1} = a_0^{-\ell}va_0v^{-1}a_0^\ell = a_0^{-\ell}a_0^ka_0^\ell = a_0^k$, proving (a). For (b), using (\ref{37OdoNormExpl1}) we first show that
$$w_ka_1w_k^{-1} \,=\, a_0^{-\ell}va_1v^{-1}a_0^\ell 
\,=\, a_0^{-\ell} a_0^{2\ell}a_1 a_0^\ell 
\,=\, a_0^{-\ell} a_0^{2\ell} a_0^{-\ell} a_1 
\,=\, a_1.$$
Thus $w_k$ commutes with~$a_1$. In (b) it follows that $u := w_k w_{k'} w_{kk'}^{-1}$ commutes with~$a_1$. But this element also commutes with $a_0$ and is therefore equal to $a_0^\lambda$ for some $\lambda\in\BZ_2$. Since this commutes with $a_1$, we have $a_0^\lambda = a_1a_0^\lambda a_1^{-1} = a_0^{-\lambda}$ and hence $\lambda=0$. Thus $u=1$, proving (b). 
Finally, (c) is a direct consequence of (a) and (b).
\end{Proof}

\medskip
The case $s=1$ and $r=2$ is also the only one among the strictly pre-periodic cases where the strong rigidity in the sense of Theorem \ref{25SmallRigid} holds. 
In fact, the rigidity is easy to verify directly in this case using the structure of $G$ as a dihedral group.

\medskip
The next smallest case $s=1$ and $r=3$ is the closure of the Grigorchuk group, first introduced by Grigorchuk \cite{Grigorchuk1980}. Grigorchuk \cite[\S\S12-15]{Grigorchuk2000b} and Grigorchuk-Sidki \cite{GrigorchukSidki2004} obtain results for the Grigorchuk group similar to the ones we obtain for its closure, for instance a description of the normalizer. Also, the results of Nekrashevych \cite{Nekrashevych-2007} imply some of the results of the preceding subsection in this case.



\subsection{Normalizer, subcase (a)}
\label{35Normalizer}

Next we will determine the normalizer
\UseTheoremCounterForNextEquation
\begin{equation}\label{35NormDef}
N := \Norm_W(G).
\end{equation}
Throughout this subsection we assume:
\UseTheoremCounterForNextEquation
\begin{equation}\label{35Ass}
\hbox{$s\ge2$\ \ and\ \ $r\ge4$.}
\end{equation}

\begin{Lem}\label{35Lem0}
In this case we have:
\begin{enumerate}
\item[(a)] $H\cap\Langle a_s,a_r\Rangle = \Langle(a_sa_r)^2\Rangle$
\item[(b)] $G/H$ is of order $4$ and generated by the images of $a_s$ and~$a_r$.
\item[(c)] $H = \Ker(\sgn_r|_G)\cap\Ker(\sgn_s|_G)$.
\end{enumerate}
\end{Lem}

\begin{Proof}
Part (a) follows from Proposition \ref{32GenCaseSquaresInH} and Lemma \ref{32IndexAtLeast4}, and it implies (b). For (c) observe that the proof of Lemma \ref{32IndexAtLeast4}, without restriction to~$T_n$, shows that $H$ is contained in the subgroup $\Ker(\sgn_r|_G)\cap\Ker(\sgn_s|_G)$ of index $4$ of~$G$. Since $H$ already has index $4$ in the present case, it is therefore equal to that subgroup. 
\end{Proof}

\begin{Lem}\label{35Lem1}
The group $N$ normalizes $H$ and acts trivially on $G/H$.
\end{Lem}

\begin{Proof}
Direct consequence of Lemma \ref{35Lem0} (c) and the fact that $\sgn_s$ and $\sgn_r$ are defined on all of~$W$.
\end{Proof}

\medskip
As before, let $\diag\colon W\to W\times W$, $w\mapsto (w,w)$ denote the diagonal embedding.

\begin{Lem}\label{35Lem2}
We have $N=G\cdot\diag(N)$.
\end{Lem}

\begin{Proof}
Consider any element $w\in W$ and write it in the form $\sigma^\lambda\,(u,v)$ with $\lambda\in\{0,1\}$ and $u,v\in W$. Since $\sigma\in G$, we have $w\in N$ if and only if $(u,v)\in N$. If that is so, then $(u,v)$ normalizes $G^1 = G\cap(W\times W)$; so by Proposition \ref{32G1Prop} both $u$ and $v$ normalize $\proj_1(G^1)=\proj_2(G^1)= \nobreak G$. Thus a necessary condition is that $u,v\in N$.

Conversely assume that $u,v\in N$. Then $(u,v)$ already normalizes $H\times H$ by Lemma \ref{35Lem1}. It also acts trivially on $(G\times G)/(H\times H)$, and so the conjugate $(u,v)\,(a_s,a_r)\,(u,v)^{-1}$ lies in $(H\times H)\cdot(a_s,a_r)$. From Proposition \ref{32HAlmostDirect} (b) it thus follows that $(u,v)$ normalizes $G$ if and only if the conjugate $(u,v)\,\sigma\,(u,v)^{-1}$ lies in~$G$. Since $\sigma\in G$, this is so if and only if the commutator $(u,v)\,\sigma\,(u,v)^{-1}\sigma^{-1} = (uv^{-1},vu^{-1})$ lies in~$G$.
By Lemma \ref{32SigmaConj} this is equivalent to $uv^{-1}\in K$, and hence by Lemma \ref{32CommK} to the existence of $(x,y)\in G^1$ such that $uv^{-1}=xy^{-1}$. But that last equality is equivalent to 
$$w\ =\ \sigma^\lambda\,(u,v) 
\ =\ \sigma^\lambda\,(xy^{-1}v,v)
\ =\ \sigma^\lambda\,(x,y)\,(y^{-1}v,y^{-1}v),$$
or equivalently to $w\in z\diag(N)$ for $z := \sigma^\lambda\,(x,y) \in G$.
\end{Proof}

\medskip
Now we recursively define elements 
\UseTheoremCounterForNextEquation
\begin{equation}\label{35WiDef}
\left\{\!\begin{array}{lll}
w_1     &\!\!\!:=\, (a_s,a_s), & \\[3pt]
w_{i+1} &\!\!\!:=\, (w_i,w_i) & \hbox{for all $i\ge1$.}\\
\end{array}\!\right\}
\end{equation}
Since $a_s\in G\subset N$, by induction Lemma \ref{35Lem2} implies that $w_i\in N$ for all $i\ge1$. Also, since $a_s$ has order two, by induction the same follows for all~$w_i$. Moreover, by induction we find that the restriction $w_i|_{T_n}$ is trivial for all $i\ge n$; hence the sequence $w_1,w_2,\ldots$ converges to $1$ within~$N$. Thus the following map is well-defined:
\UseTheoremCounterForNextEquation
\begin{equation}\label{35MapDef}
\textstyle \phi\colon
\prod\limits_{i=1}^\infty\BF_2 \longto N, \quad (k_1,k_2,\ldots) \mapsto w_1^{k_1}w_2^{k_2}\cdots.
\end{equation}
By construction it is continuous, but not a homomorphism. Note that it satisfies the basic formula
\UseTheoremCounterForNextEquation
\begin{equation}\label{35MapForm}
\phi(k_1,k_2,\ldots) \ =\ w_1^{k_1}\cdot\diag(\phi(k_2,k_3,\ldots)).
\end{equation}

\begin{Lem}\label{35Lem3}
The map $\phi$ induces a homomorphism $\bar\phi\colon \prod_{i=1}^\infty\BF_2 \to N/G$.
\end{Lem}

\begin{Proof}
It suffices to show that for all $1\le i<j$ the images of $w_i$ and $w_j$ in $N/G$ commute with each other, or in other words that the commutator $[w_i,w_j]$ lies in~$G$. We will prove this by induction on~$i$.

If $i=1$, we have $[w_1,w_j] = \diag([a_s,w_{j-1}])$. By Lemma \ref{35Lem1} the element $w_{j-1}$ normalizes $H$ and acts trivially on $G/H$; hence $[a_s,w_{j-1}]$ lies in~$H$. Therefore $[w_1,w_j]$ lies in $H\times H \subset G$, as desired.

For $i>1$ we have  $[w_i,w_j] = \diag([w_{i-1},w_{j-1}])$, and by the induction hypothesis we already know that $[w_{i-1},w_{j-1}]\in G$. On the other hand, as a commutator in~$W$, the element $[w_{i-1},w_{j-1}]$ lies in the kernel of both $\sgn_s$ and $\sgn_r$. By Lemma \ref{35Lem0} (c) it therefore lies in~$H$, and so $[w_i,w_j]$ lies in $H\times H \subset G$, as desired.
\end{Proof}

\begin{Lem}\label{35Lem4}
The homomorphism $\bar\phi$ is injective.
\end{Lem}

\begin{Proof}
Let $(k_1,k_2,\ldots)$ be an element of the kernel of~$\bar\phi$. Then $w := \phi(k_1,k_2,\ldots)$ lies in~$G$. We first prove that $k_1=0$.

Suppose to the contrary that $k_1=1$. Then the definition of $\phi$ implies that $w = w_1\cdot\diag(u)$ for some element $u\in\diag(N)$. Thus $w=(a_su,a_su)$ actually lies in~$G^1$. By Proposition \ref{32HAlmostDirect} (c) and Lemma \ref{35Lem0} there therefore exist $\lambda,\mu\in\{0,1\}$ such that
\begin{eqnarray*}
(a_su,a_su) &\in& (H\times H)\cdot(a_s,a_r)^\lambda\cdot(a_r,a_s)^\mu \\[3pt]
\Longleftrightarrow\qquad\quad
a_su &\in& H a_s^\lambda a_r^\mu \cap H a_r^\lambda a_s^\mu.
\end{eqnarray*}
By the structure of $G/H$ the last condition requires that $\lambda=\mu$ and is therefore equivalent to 
$$a_su\ \in\ H a_s^\lambda a_r^\lambda.\qquad\qquad$$ 
Now the fact that $u=(v,v)$ for some $v\in W$ implies that $\sgn_s(u)=\sgn_r(u)=1$. Since $\sgn_s$ is also trivial on $H$ and on~$a_r$, but not on~$a_s$, the condition implies that $\lambda=1$. On the other hand, since $\sgn_r$ is also trivial on $H$ and on~$a_s$, but not on~$a_r$, the condition implies that $\lambda=0$. We have thus reached a contradiction and thereby proved that $k_1=0$.

Since $k_1=0$, the formula (\ref{35MapForm}) shows that $w= \diag(w')$ for $w' := \phi(k_2,k_3,\ldots)$. Thus $w=(w',w')$ actually lies in~$G^1$, and so $w'$ lies in $G$ by Proposition \ref{32G1Prop} (c). This means that $(k_2,k_3,\ldots)$ also lies in the kernel of~$\bar\phi$.

For every element $(k_1,k_2,\ldots)$ of the kernel of~$\bar\phi$ we have thus proved that $k_1=0$ and that $(k_2,k_3,\ldots)$ again lies in the kernel of~$\bar\phi$. By an induction on $i$ we can deduce from this that for every $i\ge1$ and every element $(k_1,k_2,\ldots)$ of the kernel of~$\bar\phi$ we have $k_i=0$. This means that the kernel of $\bar\phi$ is trivial, and so the homomorphism $\bar\phi$ is injective, as desired.
\end{Proof}

\begin{Lem}\label{35Lem5}
The homomorphism $\bar\phi$ is surjective.
\end{Lem}

\begin{Proof}
We first claim that $G\cdot\diag(G) = G\cdot\langle w_1\rangle$. To see this, recall that $G=H\cdot\Langle a_s,a_r\Rangle$ with $H\triangleleft G$ of index~$4$. This implies that 
\begin{eqnarray*}
(H\times H)\cdot\diag(G) 
&=& (H\times H)\cdot \Langle (a_s,a_s),(a_r,a_r)\Rangle \\
&=& (H\times H)\cdot \Langle (a_sa_r,a_sa_r)\Rangle \cdot\langle w_1\rangle.
\end{eqnarray*}
Since $(a_sa_r)^2\in H$ in the present case, we have $a_sa_r \in H a_ra_s$ and hence 
\begin{eqnarray*}
(H\times H)\cdot\diag(G) 
&=& (H\times H)\cdot \Langle (a_sa_r,a_ra_s)\Rangle \cdot\langle w_1\rangle.
\end{eqnarray*}
As $(a_sa_r,a_ra_s) = (a_s,a_r)\,(a_r,a_s)\in G$, it follows that
$$G\cdot\diag(G) 
\ =\ G\cdot \Langle (a_sa_r,a_ra_s)\Rangle \cdot\langle w_1\rangle
\ =\ G\cdot \langle w_1\rangle,$$
as claimed.

For the surjectivity we must prove that $N = G\cdot\phi\bigl(\prod_{i=1}^\infty\BF_2\bigr)$. For this it suffices to show that $N_n = G_n\cdot\phi\bigl(\prod_{i=1}^\infty\BF_2\bigr)|_{T_n}$ for all~$n\ge0$. This is trivial for $n=0$, so assume that $n>0$ and that the equality holds for~$n-1$. Using, in turn, Lemma \ref{35Lem2}, the induction hypothesis, the above claim, and the formula (\ref{35MapForm}) we deduce that
\begin{eqnarray*}
N_n &=& G_n\cdot\diag(N_{n-1}) \\
&=& \textstyle G_n\cdot\diag\bigl( G_{n-1} \cdot 
\phi\bigl(\prod_{i=1}^\infty\BF_2\bigr)|_{T_{n-1}} \bigr) \\
&=& \textstyle G_n\cdot \langle w_1|_{T_n}\rangle \cdot 
\diag\bigl(\phi\bigl(\prod_{i=1}^\infty\BF_2\bigr)|_{T_{n-1}}\bigr) \\
&=& \textstyle G_n\cdot \langle w_1|_{T_n}\rangle \cdot 
\phi\bigl(\{0\}\times\prod_{i=2}^\infty\BF_2\bigr)|_{T_n} \\
&=& \textstyle G_n\cdot \phi\bigl(\prod_{i=1}^\infty\BF_2\bigr)|_{T_n},
\end{eqnarray*}
so the equality holds for~$n$. Thus it follows for all $n\ge0$ by induction, and we are done.
\end{Proof}

\medskip
Combining Lemmas \ref{35Lem3} through \ref{35Lem5} now implies:

\begin{Thm}\label{35NormThm}
In the case (\ref{35Ass}) the map $\phi$ induces an isomorphism 
$$\textstyle \bar\phi\colon \prod\limits_{i=1}^\infty\BF_2 \stackrel{\sim}{\longto} N/G.$$
\end{Thm}


\subsection{Normalizer, subcase (b)}
\label{36Normalizer}

Throughout this subsection we assume:
\UseTheoremCounterForNextEquation
\begin{equation}\label{36Ass}
\hbox{$s=1$\ \ and\ \ $r\ge3$.}
\end{equation}

\begin{Lem}\label{36Lem0}
In this case we have:
\begin{enumerate}
\item[(a)] $\langle a_s,a_r\rangle$ is a finite dihedral group of order~$8$ with center $\langle(a_sa_r)^2\rangle$.
\item[(b)] $G = H\rtimes\langle a_s,a_r\rangle$.
\item[(c)] $G = (H\times H) \rtimes \langle \sigma, (a_s,a_r) \rangle$.
\item[(d)] $H^+ := H\rtimes\langle(a_sa_r)^2\rangle = \Ker(\sgn_s|_G)\cap\Ker(\sgn_r|_G)$.
\end{enumerate}
\end{Lem}

\begin{Proof}
Part (a) follows from Proposition \ref{32aiajOrder}, part (b) from Proposition \ref{32Intersections}, and part (c) from Proposition \ref{32HAlmostDirect} (b). Part (d) follows by the same arguments as in the proof of Lemma \ref{35Lem0}.
\end{Proof}

\begin{Lem}\label{36Lem1}
The group $N$ normalizes $H$ and induces only inner automorphisms on $G/H$.
\end{Lem}

\begin{Proof}
Recall from (\ref{32HDef}) that $H$ is the closure of the subgroup of $G$ generated by all $G$-conjugates of $a_i$ for all $i\not=s,r$. By Proposition \ref{33WConjInGIsGConj} these $G$-conjugates are precisely the $W$-conjugates of $a_i$ which lie in~$G$. The set of these is therefore preserved under conjugation by~$N$; hence $H$ is normalized by~$N$, as desired.

Next, the same argument as in the proof of Lemmas  \ref{35Lem0} (c) and \ref{35Lem1} shows that $N$ normalizes $H^+$ and acts trivially on $G/H^+$. This means that the action of $N$ on the dihedral group $G/H\cong D_4$ is trivial modulo the center of~$D_4$. But every automorphism of $D_4$ which is trivial modulo its center is inner.
\end{Proof}

\medskip
Let $M$ denote the subgroup of all elements of $N$ which induce the identity on~$G/H$. 

\begin{Lem}\label{36Lem2}
We have $\diag(M) \subset N$.
\end{Lem}

\begin{Proof}
Consider any element $u\in M$. Then $\diag(u) = (u,u)$ normalizes $H\times H$ by Lemma \ref{36Lem1}, and it evidently centralizes~$\sigma$. Moreover, since $u$ acts trivially on $G/H$, the conjugate $(u,u)\,(a_s,a_r)\,(u,u)^{-1} = (ua_su^{-1},ua_ru^{-1})$ lies in $Ha_s\times Ha_r = (H\times H)\cdot(a_s,a_r)$ and hence in~$G$. By Lemma \ref{36Lem0} (c) the element $(u,u)$ therefore conjugates $G$ into itself and thus lies in~$N$, as desired.
\end{Proof}

%

\begin{Lem}\label{36Lem4}
We have $N\cap(1\times W) = 1\times H^+$.
\end{Lem}

\begin{Proof}
We take any element $v\in W$ and ask under which conditions $(1,v)$ normalizes~$G$. 

First the conjugate $(1,v)\,\sigma\,(1,v)^{-1} = \sigma\,(v,v^{-1})$ must lie in~$G$. This is equivalent to $(v,v^{-1}) \in G$, and hence by Lemma \ref{32SigmaConj} to $v\in K$. By the definition (\ref{32KDef}) of $K$ this means that $v \in H\cdot(a_sa_r)^\lambda$ for some $\lambda\in\BZ_2$. Granting this, we in particular have $v\in G$; hence $(1,v)$ normalizes $H\times H$. 
We also have
$$va_rv^{-1} \ \in\ H\cdot (a_sa_r)^\lambda a_r (a_sa_r)^{-\lambda} 
\ =\ H\cdot a_r(a_sa_r)^{-2\lambda}$$
and hence 
$$(1,v)\,(a_s,a_r)\,(1,v)^{-1}
\ =\ (a_s,va_rv^{-1})
\ \in\ (1\times H)\cdot(a_s,a_r)\cdot \bigl(1,(a_sa_r)^{-2\lambda}\bigr).$$
By Lemma \ref{36Lem0} (c) the last condition is that this conjugate lies in~$G$, which therefore means that
$\bigl(1,(a_sa_r)^{-2\lambda}\bigr)$ lies in~$G$.
By Lemma \ref{36Lem0} this is equivalent to $(a_sa_r)^{-2\lambda}=1$ and hence to $2\lambda\equiv0\mod4$. Thus $\lambda$ must be even, and so $(1,v)$ normalizes $G$ if and only if $v\in H\cdot\langle(a_sa_r)^2\rangle = H^+$, as desired.
\end{Proof}

\begin{Lem}\label{36Lem5}
We have $N = G\cdot\langle(1,(a_sa_r)^2)\rangle\cdot\diag(M)$.
\end{Lem}

\begin{Proof}
By Lemma \ref{36Lem2} we have $\diag(M) \subset N$, and by Lemma \ref{36Lem4} we in particular have $\langle(1,(a_sa_r)^2)\rangle \subset N$. Together this implies the inclusion ``$\supset$''. 

To prove the inclusion ``$\subset$'' consider an arbitrary element $w\in N$. Write it in the form $w=\sigma^\lambda\,(u,v)$ for some $\lambda\in\{0,1\}$ and $u,v\in W$. Then $(u,v)$ normalizes $G$ and hence $G^1 = G\cap(W\times W)$, and so Proposition \ref{32G1Prop} shows that $u$ normalizes $\proj_1(G^1)=\nobreak G$. Thus $u$ lies in~$N$.

By Lemma \ref{36Lem1} the element $u$ therefore induces an inner automorphism on the dihedral group $G/H$. Choose $\mu,\nu\in\{0,1\}$ such that $a_s^\mu a_r^\nu$ induces the same inner automorphism on $G/H$. Then $(a_s,a_r)^\mu\,(a_r,a_s)^\nu$ lies in~$G$, and so we can write $w = \sigma^\lambda\,(a_s,a_r)^\mu\,(a_r,a_s)^\nu\,w'$ for another element $w'=(u',v')\in N$, where now $u'$ induces the trivial automorphism on $G/H$. By the definition of $M$ this means that $u'\in M$. 

By Lemma \ref{36Lem2} we therefore have $(u',u')\in N$, and so $w' = w''\,(u',u')$ for some element $w'' \in N\cap(1\times W)$. By Lemma \ref{36Lem4} this means that $w'' \in 1\times H^+ = 1\times (H\cdot\langle(a_sa_r)^2\rangle)$. Altogether we therefore have
\begin{eqnarray*}
w &\in& G\cdot \bigl(1\times (H\cdot\langle(a_sa_r)^2\rangle)\bigr) \cdot \diag(M) \\[3pt]
&=& G\cdot (1\times H)\cdot \langle(1,(a_sa_r)^2)\rangle \cdot \diag(M) \\[3pt]
&=& G\cdot \langle(1,(a_sa_r)^2)\rangle \cdot \diag(M),
\end{eqnarray*}
as desired.
\end{Proof}

\begin{Lem}\label{36Lem9}
We have $\diag(H^+) \subset H^+$.
\end{Lem}

\begin{Proof}
Since $(a_sa_r)^2 = (a_ra_s)^2$, we have
\begin{eqnarray*}
\diag(H^+)\ =\ \diag\bigl(H\cdot\langle(a_sa_r)^2\rangle\bigr)
&\!\!=\!\!&  \diag(H)\cdot\langle((a_sa_r)^2,(a_ra_s)^2)\rangle \\
&\!\!=\!\!& \diag(H)\cdot\langle((a_s,a_r)\,(a_r,a_s))^2\rangle 
\ \subset\ G.
\end{eqnarray*}
On the other hand, both $\sgn_s$ and $\sgn_r$ are trivial on all elements of the form $(u,u)$. Thus Lemma \ref{36Lem0} (d) implies that $\diag(H^+) \subset H^+$, as desired.
\end{Proof}

\begin{Lem}\label{36Lem6}
We have $\langle(1,(a_sa_r)^2)\rangle \subset M$.
\end{Lem}

\begin{Proof}
{}From Lemma \ref{36Lem5} we already know that $(1,(a_sa_r)^2)$ normalizes $G$. By Lemma \ref{36Lem1} it therefore normalizes~$H$, and it remains to show that conjugation by it induces the identity on $G/H$. For this it suffices to show that the commutators $[(1,(a_sa_r)^2),a_s]$ and $[(1,(a_sa_r)^2),a_r]$ lie in~$H$. 
Since $r\ge3>2=s+1$ and therefore $a_r=(a_{r-1},1)$, the second commutator is trivial. Since $a_s=a_1=\sigma$, the first commutator is 
$$\bigl((a_1a_r)^{-2},(a_1a_r)^2\bigr)
\ =\ ((a_ra_1)^2,(a_1a_r)^2)
\ =\ \bigl((a_r,a_1)\,(a_1,a_r)\bigr)^2
\ =\ \bigl(\sigma\,a_2\,\sigma^{-1}\,a_2\bigr)^2.$$
As $a_2$ and $\sigma\,a_2\,\sigma^{-1}$ lie in $H$ by the definition of~$H$, it follows that this commutator lies in~$H$, as desired.
\end{Proof}

\begin{Lem}\label{36Lem7}
We have $\diag(M) \subset M$.
\end{Lem}

\begin{Proof}
Consider any element $u\in M$. From Lemma \ref{36Lem2} we already know that $(u,u)$ normalizes $G$. By Lemma \ref{36Lem1} it therefore normalizes~$H$, and it remains to show that conjugation by it induces the identity on $G/H$. For this we must show that the commutators $[(u,u),a_s]$ and $[(u,u),a_r]$ lie in~$H$. 

Since $s=1$, we have $a_s=a_1=\sigma$, and the first commutator is trivial. Since $r>s+1$, the second commutator is equal to $[(u,u),(a_{r-1},1)] = ([u,a_{r-1}],1)$. {}From Lemma \ref{36Lem1} we know that $(u,u)$ induces an inner automorphism on $G/H$. Thus the commutator already lies in $H\cdot\langle(a_sa_r)^2\rangle$. As $(a_sa_r)^2 = (\sigma\,(a_{r-1},1))^2 = (a_{r-1},a_{r-1})$, which is of order~$2$, we therefore have 
\UseTheoremCounterForNextEquation
\begin{equation}\label{36Lem7Disp}
([u,a_{r-1}],1) \in H\cdot(a_{r-1},a_{r-1})^\lambda
\qquad\hbox{for some $\lambda\in\{0,1\}$.}
\end{equation}
Now recall that $H$ is the closure of the subgroup generated by all $G$-conjugates of $a_i$ for all $i\not=s,r$, in other words for all $2\le i\le r-1$. Here $a_i=(a_{i-1},1)$ for $i>2$ and $a_2=(a_1,a_r)$. Thus all generators have the form $(x,y) \in W\times W$ with $\sgn_{r-1}(x)=\sgn_{r-1}(y)=1$, and so all elements of $H$ have the same property. Since the element $([u,a_{r-1}],1)$ does so, too, the relation (\ref{36Lem7Disp}) implies that $(a_{r-1},a_{r-1})^\lambda$ does as well, and hence that $\lambda=0$. Thus the commutator lies in~$H$, as desired.
\end{Proof}

\medskip
Now we recursively define elements 
\UseTheoremCounterForNextEquation
\begin{equation}\label{36WiDef}
\left\{\!\begin{array}{lll}
w_1     &\!\!\!:=\, (1,(a_sa_r)^2), & \\[3pt]
w_{i+1} &\!\!\!:=\, (w_i,w_i) & \hbox{for all $i\ge1$.}\\
\end{array}\!\right\}
\end{equation}
Since $w_1\in M$ by Lemma \ref{36Lem6}, by induction Lemma \ref{36Lem7} implies that $w_i\in M$ for all $i\ge1$. Also, since $(a_sa_r)^2$ has order two, so does $w_1$, and by induction the same follows for all~$w_i$. Moreover, by induction we find that the restriction $w_i|_{T_n}$ is trivial for all $i\ge n$; hence the sequence $w_1,w_2,\ldots$ converges to $1$ within~$N$. Thus the following map is well-defined:
\UseTheoremCounterForNextEquation
\begin{equation}\label{36MapDef}
\textstyle \phi\colon
\prod\limits_{i=1}^\infty\BF_2 \longto M, \quad (k_1,k_2,\ldots) \mapsto w_1^{k_1}w_2^{k_2}\cdots.
\end{equation}
By construction it is continuous, but not a homomorphism. Note that it satisfies the basic formula
\UseTheoremCounterForNextEquation
\begin{equation}\label{36MapForm}
\phi(k_1,k_2,\ldots) \ =\ w_1^{k_1}\cdot\diag(\phi(k_2,k_3,\ldots)).
\end{equation}

\begin{Lem}\label{36Lem13}
The map $\phi$ induces a homomorphism $\bar\phi\colon \prod_{i=1}^\infty\BF_2 \to N/G$.
\end{Lem}

\begin{Proof}
It suffices to show that for all $1\le i<j$ the images of $w_i$ and $w_j$ in $N/G$ commute with each other; in other words that the commutator $[w_i,w_j]$ lies in~$G$. We will prove this by induction on~$i$.
If $i=1$, we have 
$$[w_1,w_j] = \bigl[(1,(a_sa_r)^2),(w_{j-1},w_{j-1})\bigr] = \bigl(1,[(a_sa_r)^2,w_{j-1}]\bigr).$$
Since $(a_sa_r)^2 \in G$ and $w_{j-1}\in M$, the commutator $[(a_sa_r)^2,w_{j-1}]$ lies in $H$, and so $[w_1,w_j]$ lies in $1\times H \subset G$, as desired.
For $i>1$ we have  $[w_i,w_j] = \diag([w_{i-1},w_{j-1}])$, and by the induction hypothesis we already know that $[w_{i-1},w_{j-1}]\in G$. But being a commutator, it also lies in the kernel of both $\sgn_s$ and $\sgn_r$; hence we have $[w_{i-1},w_{j-1}]\in H^+$ by Lemma \ref{36Lem0} (d). Thus Lemma \ref{36Lem9} implies that $[w_i,w_j] \in H^+ \subset G$, as desired.
\end{Proof}

\begin{Lem}\label{36Lem14}
The homomorphism $\bar\phi$ is injective.
\end{Lem}

\begin{Proof}
Let $(k_1,k_2,\ldots)$ be an element of the kernel of~$\bar\phi$. Then $w := \phi(k_1,k_2,\ldots)$ lies in~$G$. We first prove that $k_1=0$.

Suppose to the contrary that $k_1=1$. Then the definition of $\phi$ implies that $w = w_1\cdot\diag(u)$ for some element $u\in M$. Thus $w=(u,(a_sa_r)^2u)$ actually lies in~$G^1$. By Proposition \ref{32HAlmostDirect} (c) and Lemma \ref{36Lem0} there therefore exist $\lambda\in\BZ$ and $\mu\in\{0,1\}$ such that 
\begin{eqnarray*}
(u,(a_sa_r)^2u) &\in& (H\times H)\cdot(a_sa_r,a_ra_s)^\lambda\cdot(a_s,a_r)^\mu \\[5pt]
\Longleftrightarrow\quad\ \quad
\rlap{$\biggl\{\begin{array}{lll}
u &\!\!\!\in\ \;H (a_sa_r)^\lambda a_s^\mu & \hbox{and}\\[3pt]
(a_sa_r)^2u &\!\!\!\in\ \;H (a_ra_s)^\lambda a_r^\mu &
\end{array}\biggr\}$} \qquad\qquad\  && \\[5pt]
\Longleftrightarrow\qquad\qquad\quad\ \;
u \ \,&\in& H (a_sa_r)^\lambda a_s^\mu \cap (a_sa_r)^{-2} H (a_ra_s)^\lambda a_r^\mu.
\end{eqnarray*}
In particular $u\in G$ and therefore $u\in G\cap M = H^+ = H\rtimes\langle(a_sa_r)^2\rangle$. Therefore $\mu=0$ and $\lambda$ is even, and $H (a_sa_r)^\lambda = (a_sa_r)^{-2} H (a_ra_s)^\lambda$. But this contradicts the fact that $(a_sa_r)^2 = (a_ra_s)^2 \not\in H$. We have thus reached a contradiction and thereby proved that $k_1=0$.

Since $k_1=0$, the formula (\ref{36MapForm}) shows that $w= \diag(w')$ for $w' := \phi(k_2,k_3,\ldots) \in \nobreak M$. Thus $w=(w',w')\in G$ actually lies in~$G^1$, and so $w'$ lies in $G$ by Proposition \ref{32G1Prop} (c). This means that $(k_2,k_3,\ldots)$ also lies in the kernel of~$\bar\phi$.

For every element $(k_1,k_2,\ldots)$ of the kernel of~$\bar\phi$ we have thus proved that $k_1=0$ and that $(k_2,k_3,\ldots)$ again lies in the kernel of~$\bar\phi$. By an induction on $i$ we can deduce from this that for every $i\ge1$ and every element $(k_1,k_2,\ldots)$ of the kernel of~$\bar\phi$ we have $k_i=0$. This means that the kernel of $\bar\phi$ is trivial, and so the homomorphism $\bar\phi$ is injective, as desired.
\end{Proof}

\begin{Lem}\label{36Lem15}
The homomorphism $\bar\phi$ is surjective.
\end{Lem}

\begin{Proof}
We must prove that $N = G\cdot\phi\bigl(\prod_{i=1}^\infty\BF_2\bigr)$. For this it suffices to show the equality $N_n = G_n\cdot\phi\bigl(\prod_{i=1}^\infty\BF_2\bigr)|_{T_n}$ for all~$n\ge0$, which we will achieve by induction on~$n$. If $n\le r$, Proposition \ref{31GnAll} implies that $G_n=N_n$, from which the equality follows. Thus we may assume that $n>r$ and that $N_{n-1} = G_{n-1}\cdot\phi\bigl(\prod_{i=1}^\infty\BF_2\bigr)|_{T_{n-1}}$.

Since $\phi\bigl(\prod_{i=1}^\infty\BF_2\bigr) \subset M$ by construction, the induction hypothesis implies that $M_{n-1} = (G_{n-1}\cap M_{n-1})\cdot\phi\bigl(\prod_{i=1}^\infty\BF_2\bigr)|_{T_{n-1}}$. We claim that $G_{n-1}\cap M_{n-1} = H_{n-1}^+$ and consequently 
\UseTheoremCounterForNextEquation
\begin{equation}\label{36Lem15Diag}
\textstyle M_{n-1} = H_{n-1}^+\cdot\phi\bigl(\prod\limits_{i=1}^\infty\BF_2\bigr)|_{T_{n-1}}.
\end{equation}
To show the claim observe that under the given assumptions on $s,r,n$, Proposition \ref{32MidCaseSquaresNotInH} implies that $G_{n-1} = H_{n-1}\rtimes\langle a_s|_{T_{n-1}},a_r|_{T_{n-1}}\rangle$ and that restriction to $T_{n-1}$ induces an isomorphism $\langle a_s,a_r\rangle \stackrel{\sim}{\to} \langle a_s|_{T_{n-1}},a_r|_{T_{n-1}}\rangle$. Thus the restriction induces an isomorphism $G/H \stackrel{\sim}{\to} G_{n-1}/H_{n-1}$. It follows that $M$ is the subgroup of all elements of $N$ which induce the identity on $G_{n-1}/H_{n-1}$. Consequently $G_{n-1}\cap M_{n-1}$ is the subgroup of all elements of $G_{n-1}$ which induce the identity on $G_{n-1}/H_{n-1}$, which by the same semidirect product decomposition is $H_{n-1}^+ = H_{n-1}\rtimes\langle (a_sa_r)^2|_{T_{n-1}}\rangle$. Thus $G_{n-1}\cap M_{n-1} = H_{n-1}^+$, as claimed.

Now we can perform the induction step. Using, in turn, Lemma \ref{36Lem5}, the equality (\ref{36Lem15Diag}), Lemma \ref{36Lem9}, and the formula (\ref{36MapForm}) we deduce that
\begin{eqnarray*}
N_n &=& G_n\cdot\langle w_1|_{T_n}\rangle\cdot\diag(M_{n-1}) \\
&=& \textstyle G_n\cdot\langle w_1|_{T_n}\rangle\cdot
\diag\bigl( H_{n-1}^+\cdot\phi\bigl(\prod_{i=1}^\infty\BF_2\bigr)|_{T_{n-1}} \bigr) \\
&=& \textstyle G_n\cdot\langle w_1|_{T_n}\rangle\cdot
\diag\bigl( \phi\bigl(\prod_{i=1}^\infty\BF_2\bigr) \bigr)|_{T_n} \\
&=& \textstyle G_n\cdot\phi\bigl(\prod_{i=1}^\infty\BF_2\bigr)|_{T_n},
\end{eqnarray*}
so the desired equality holds for~$n$, and we are done.
\end{Proof}

\medskip
Combining Lemmas \ref{36Lem13} through \ref{36Lem15} now implies:

\begin{Thm}\label{36NormThm}
In the case (\ref{36Ass}) the map $\phi$ induces an isomorphism 
$$\textstyle \bar\phi\colon \prod\limits_{i=1}^\infty\BF_2 \stackrel{\sim}{\longto} N/G.$$
\end{Thm}


\subsection{Normalizer, subcase (c)}
\label{37Normalizer}

Throughout this subsection we assume:
\UseTheoremCounterForNextEquation
\begin{equation}\label{37Ass}
\hbox{$s=2$\ \ and\ \ $r=3$.}
\end{equation}
Thus we have $G=\Langle a_1,a_2,a_3\Rangle$ with $a_1=\sigma$ and $a_2=(a_1,1)$ and $a_3=(a_2,a_3)$, and the definitions (\ref{32HDef}) and (\ref{32LDef}) of $H$ and $L$ imply that $H=L$. Surprisingly, the determination of the normalizer requires more work in this case than in any of the other cases. Perhaps it is a kind of mixture of the cases in the two preceding subsections.

\begin{Lem}\label{37Lem0}
In this case we have:
\begin{enumerate}
\item[(a)] $\langle a_s,a_r\rangle$ is a finite dihedral group of order~$8$ with center $\langle(a_sa_r)^2\rangle$.
\item[(b)] $G = H\rtimes\langle a_s,a_r\rangle$.
\item[(c)] $G = (H\times H) \rtimes \langle \sigma, (a_s,a_r) \rangle$.
\item[(d)] $H^+ := H\rtimes\langle(a_sa_r)^2\rangle = \Ker(\sgn_s|_G)\cap\Ker(\sgn_r|_G)$.
\end{enumerate}
\end{Lem}

\begin{Proof}
Same as that of Lemma \ref{36Lem0}.
\end{Proof}

\begin{Lem}\label{37Lem1}
The group $N$ normalizes $H$ and induces only inner automorphisms on $G/H$.
\end{Lem}

\begin{Proof}
By definition $H$ is the closure of the subgroup of $G$ generated by all $G$-conjugates of~$a_1$. By Proposition \ref{33WConjInGIsGConj} these generators are precisely the $W$-conjugates of $a_1$ which lie in~$G$. The set of these is therefore preserved under conjugation by~$N$, and hence $H$ is normalized by~$N$. The rest of the proof proceeds as that of Lemma \ref{36Lem1}.
\end{Proof}

\medskip
Let $M$ denote the subgroup of all elements of $N$ which induce the identity on~$G/H$. The same proofs as those of Lemmas \ref{36Lem2} through \ref{36Lem9} yield the next four lemmas:

\begin{Lem}\label{37Lem2}
We have $\diag(M) \subset N$.
\end{Lem}


\begin{Lem}\label{37Lem4}
We have $N\cap(1\times W) = 1\times H^+$.
\end{Lem}


\begin{Lem}\label{37Lem5}
We have $N = G\cdot\langle(1,(a_sa_r)^2)\rangle\cdot\diag(M)$.
\end{Lem}


\begin{Lem}\label{37Lem9}
We have $\diag(H^+) \subset H^+$.
\end{Lem}


\medskip
At this point things begin to differ with respect to Subsection \ref{36Normalizer}. First we recall the relation between the subgroups $H_n \subset H_n^+ \subset G_n$ which, as usual, denote the images of $H\subset H^+\subset G$ in~$W_n$. 

\begin{Lem}\label{37Lem41}
For all $n\ge4$ we have:
\begin{enumerate}
\item[(a)] $G_n = H_n\rtimes\langle a_s|_{T_n},a_r|_{T_n}\rangle$ with $\langle a_s|_{T_n},a_r|_{T_n}\rangle \cong \langle a_s,a_r\rangle \cong G/H$ dihedral of order~$8$.
\item[(b)] $H^+_n = H_n \rtimes \langle (a_sa_r)^2|_{T_n} \rangle$
with $\langle (a_sa_r)^2|_{T_n} \rangle \cong \langle (a_sa_r)^2 \rangle$ cyclic of order~$2$.
\item[(c)] But $(a_sa_r)^2|_{T_3}\in H_3$.
\end{enumerate}
\end{Lem}

\begin{Proof}
Since $(a_sa_r)^2$ has order~$2$, so does its restriction $(a_sa_r)^2|_{T_n}$ for all $n\ge4$ by Proposition \ref{32SmallCaseSquaresInAndNotInH}. Using this, 
Lemma \ref{37Lem0} implies both (a) and (b). Assertion (c) is part of Proposition \ref{32SmallCaseSquaresInAndNotInH}.
\end{Proof}

\begin{Lem}\label{37Lem6}
We have $\langle(1,(a_sa_r)^2)\rangle \subset M$.
\end{Lem}

\begin{Proof}
{}From Lemma \ref{37Lem5} we already know that $(1,(a_sa_r)^2)$ normalizes $G$. By Lemma \ref{37Lem1} it therefore normalizes~$H$, and it remains to show that conjugation by it induces the identity on $G/H$. For this it suffices to show that the commutators $[(1,(a_sa_r)^2),a_s]$ and $[(1,(a_sa_r)^2),a_r]$ lie in~$H$. 
Since $a_s=a_2=(a_1,1)$ in this case, the first commutator is trivial. Since $a_r = a_3 = (a_s,a_r)$, and $(a_sa_r)^2$ commutes with $a_r$ by Lemma \ref{37Lem0} (a), the second commutator is trivial as well, and we are done.
\end{Proof}

\begin{Lem}\label{37Lem7}
For any $u\in M$ with $u|_{T_3}=1$ we have $\diag(u) \in M$.
\end{Lem}

\begin{Proof}
{}From Lemma \ref{37Lem2} we already know that $(u,u)$ normalizes~$G$. By Lemma \ref{37Lem1} it therefore normalizes~$H$. On the other hand the assumption implies that $(u,u)|_{T_4}=1$; hence conjugation by it induces the identity on $G_4/H_4$. Thus Lemma \ref{37Lem41} (a) implies that conjugation by $(u,u)$ induces the identity on $G/H$; hence $(u,u)\in M$, as desired.
\end{Proof}

\medskip
As in Subsection \ref{36Normalizer} we recursively define elements 
\UseTheoremCounterForNextEquation
\begin{equation}\label{37WiDef}
\left\{\!\begin{array}{lll}
w_1     &\!\!\!:=\, (1,(a_sa_r)^2), & \\[3pt]
w_{i+1} &\!\!\!:=\, (w_i,w_i) & \hbox{for all $i\ge1$.}\\
\end{array}\!\right\}
\end{equation}

\begin{Lem}\label{37Lem11}
\begin{enumerate}
\item[(a)] For all $i\ge1$ we have $w_i\in M$.
\item[(b)] For all $i\ge1$ we have $w_i|_{T_4} \in G_4$.
\end{enumerate}
\end{Lem}

\begin{Proof}
First observe that $(a_1a_2)^2|_{T_1} = 1$ because it is a square. Therefore $(a_2a_3)^2|_{T_2} = ((a_1,1)(a_2,a_3))^2|_{T_2} = ((a_1a_2)^2,1)|_{T_2} = 1$, and hence $w_1|_{T_3} = 1$. By induction on $i$ it follows that $w_i|_{T_3} = 1$ for all $i\ge1$ and that $w_i|_{T_4} = 1$ for all $i\ge2$. 

Since $w_1\in M$ by Lemma \ref{37Lem6}, by induction Lemma \ref{37Lem7} now implies that $w_i\in M$ for all $i\ge1$, proving (a).
Also, from Lemma \ref{37Lem41} (c) we deduce that $w_1|_{T_4} = (1,(a_sa_r)^2|_{T_3})$ lies in $1\times H_3$, and hence in $G_4$ by Lemma \ref{37Lem0} (c). For $i\ge2$ the element $w_i|_{T_4} = 1$ trivially lies in~$G_4$. This proves (b).
\end{Proof}

\medskip
To describe $N$ we will need one more element $w_0\in W$\!, defined by the recursion relation
\UseTheoremCounterForNextEquation
\begin{equation}\label{37W0Def}
w_0\ =\ (a_2w_0,a_3w_0)\ =\ a_3\,(w_0,w_0).
\end{equation}

\begin{Lem}\label{37Lem01}
We have $[w_0,a_2] = [[a_3,a_1],a_2] \in H$.
\end{Lem}

\begin{Proof}
The fact that $(w_0,w_0)$ commutes with $a_1=\sigma$ implies that
$$[a_2w_0,a_1] = [a_2a_3(w_0,w_0),\sigma] = [a_2a_3,a_1].$$
{}From this we deduce that
$$[w_0,a_2] = [(a_2w_0,a_3w_0),(a_1,1)] = ([a_2w_0,a_1],1) = ([a_2a_3,a_1],1) = [(a_2a_3,a_3a_2),(a_1,1)].$$
Since 
$$(a_2a_3,a_3a_2) = (a_2,a_3)\,(a_3,a_2)^{-1} = [(a_2,a_3),\sigma] = [a_3,a_1],$$
it follows that $[w_0,a_2] = [[a_3,a_1],a_2]$. But by definition $a_1$ lies in~$H$, and $H$ is a normal subgroup of~$G$; hence $[[a_3,a_1],a_2] \in H$, as desired.
\end{Proof}

\begin{Lem}\label{37Lem02}
We have $[w_0,a_3] \in H$.
\end{Lem}

\begin{Proof}
Set $x := [w_0,a_3]$. Then 
$$x = [(a_2w_0,a_3w_0),(a_2,a_3)] 
    = ([a_2w_0,a_2],[a_3w_0,a_3]) 
    = (a_2[w_0,a_2]a_2^{-1},a_3[w_0,a_3]a_3^{-1}),$$
which using Lemma \ref{37Lem01} implies that
\UseTheoremCounterForNextEquation
\begin{equation}\label{37XDef}
x = \bigl(a_2[[a_3,a_1],a_2]a_2^{-1},a_3xa_3^{-1}\bigr).
\end{equation}
This recursion relation determines $x$ uniquely and independently of~$w_0$. Thus we must show that the element $x\in W$ characterized by (\ref{37XDef}) lies in~$H$.
We will show by induction that $x|_{T_n}$ lies in $H_n$ for all $n\ge0$. 

As a preparation for this we calculate
$$[a_1a_3a_1^{-1},a_2] = \bigl[\sigma\,(a_2,a_3)\,\sigma^{-1},(a_1,1)\bigr]
 = [(a_3,a_2),(a_1,1)] = ([a_3,a_1],1),$$
whence
\begin{eqnarray*}
a_3[[a_1a_3a_1^{-1},a_2],a_3]a_3^{-1} 
&\!\!=\!\!& (a_2,a_3)\,\bigl[([a_3,a_1],1),(a_2,a_3)\bigr]\,(a_2,a_3)^{-1} \\
&\!\!=\!\!& (a_2,a_3)\,\bigl( [[a_3,a_1],a_2],1\bigr)\,(a_2,a_3)^{-1} \\
&\!\!=\!\!& \bigl( a_2[[a_3,a_1],a_2]a_2^{-1},1\bigr).
\end{eqnarray*}
Now observe that, since $a_3|_{T_2}=1$, we also have $(a_2[[a_3,a_1],a_2]a_2^{-1})|_{T_2}=1$, and hence $x|_{T_3}=1$. This and the above calculation together imply that 
$$x|_{T_4} = \bigl(a_2[[a_3,a_1],a_2]a_2^{-1}),1\bigr)|_{T_4}
= \bigl(a_3[[a_1a_3a_1^{-1},a_2],a_3]a_3^{-1}\bigr)|_{T_4}.$$
But $[a_1a_3a_1^{-1},a_2]$ lies in the commutator group of $G$ and hence in~$H^+$. Since $H^+/H$ is the center of $G/H$, we have $[H^+,G]\subset H$, and therefore $a_3[[a_1a_3a_1^{-1},a_2],a_3]a_3^{-1} \in H$. Together this shows that $x|_{T_4}$ lies in~$H_4$.

To finish the proof by induction, assume that $n>4$ and that $x|_{T_{n-1}}$ lies in~$H_{n-1}$. Then
$$x|_{T_n} = \bigl((a_2[[a_3,a_1],a_2]a_2^{-1})|_{T_{n-1}},(a_3xa_3^{-1})|_{T_{n-1}}\bigr),$$
where the first entry lies in $H_{n-1}$ by Lemma \ref{37Lem01}, and the second entry lies in $H_{n-1}$ by the induction hypothesis. Thus $x|_{T_n}$ lies in $H_{n-1}\times H_{n-1}$ and hence in~$G_n$. On the other hand, Lemma \ref{37Lem41} (a) implies that the restriction to $T_4$ induces an isomorphism $G_n/H_n \stackrel{\sim}{\to} G_4/H_4$, and we already know that $x|_{T_4}$ lies in~$H_4$. Together this implies that $x|_{T_n}$ lies in~$H_n$, as desired.
\end{Proof}

\begin{Lem}\label{37Lem03}
We have $w_0\in M$.
\end{Lem}

\begin{Proof}
Lemmas \ref{37Lem01} and \ref{37Lem02} already imply that $w_0a_2w_0^{-1}$ and $w_0a_3w_0^{-1}$ lie in~$G$. Moreover, the definition (\ref{37W0Def}) of $w_0$ shows that 
$$w_0a_1w_0^{-1} = 
(a_2w_0,a_3w_0)\sigma(a_2w_0,a_3w_0)^{-1} 
= (a_2,a_3)\sigma(a_2,a_3)^{-1} 
= a_3a_1a_3^{-1} 
\in G.$$
Since $G=\Langle a_1,a_2,a_3\Rangle$, it follows that $w_0Gw_0^{-1}\subset G$ and hence $w_0Gw_0^{-1}=G$. Thus $w_0$ lies in~$N$. From Lemmas \ref{37Lem01} and \ref{37Lem02} it now follows that $w_0$ lies in~$M$.
\end{Proof}

\begin{Lem}\label{37Lem04}
We have $w_0^2\in G$.
\end{Lem}

\begin{Proof}
It suffices to show that $w_0^2|_{T_n}$ lies in $G_n$ for all $n\ge0$. 
For $n\le3$ this follows from the fact that $G_n=W_n$. So assume that $n>3$ and that $w_0^2|_{T_{n-1}}$ lies in~$G_{n-1}$. Since $w_0^2|_{T_{n-1}}$ is a square, it lies in the kernel of $\sgn_2$ and $\sgn_3$, and therefore already in $H_{n-1}^+$. By Lemma \ref{37Lem9} this implies that $(w_0^2,w_0^2)|_{T_n}$ lies in $H^+_n \subset G_n$. 
But by the construction (\ref{37W0Def}) of $w_0$ we have
$$w_0^2 \ =\ (a_2w_0a_2w_0,a_3w_0a_3w_0)
\ =\ \bigl(a_2[w_0,a_2]a_2^{-1}w_0^2,a_3[w_0,a_3]a_3^{-1}w_0^2\bigr).$$
Since $[w_0,a_2]$ and $[w_0,a_3]$ lie in $H$ by Lemmas \ref{37Lem01} and \ref{37Lem02}, it follows that $w_0^2 \in (H\times H)\cdot(w_0^2,w_0^2) \subset G\cdot(w_0^2,w_0^2)$. Therefore $w_0^2|_{T_n} \in G_n\cdot(w_0^2,w_0^2)|_{T_n} = G_n$, as desired.
\end{Proof}

\begin{Lem}\label{37Lem05}
We have $[w_0,w_i] \in G$ for all $i\ge1$.
\end{Lem}

\begin{Proof}
The constructions (\ref{37WiDef}) and (\ref{37W0Def}) of $w_0$ and $w_1$ imply that
$$[w_0,w_1] = \bigl[(a_2w_0,a_3w_0),(1,(a_2a_3)^2)\bigr]
 = \bigl(1,[a_3w_0,(a_2a_3)^2]\bigr).$$
Since $w_0$ lies in~$M$, conjugation by it fixes the image of $(a_2a_3)^2$ in $G/H$. Since $a_3$ commutes with $(a_2a_3)^2$, it follows that conjugation by $a_3w_0$ also fixes the image of $(a_2a_3)^2$ in $G/H$. Therefore $[a_3w_0,(a_2a_3)^2]$ lies in~$H$, and hence $[w_0,w_1]$ lies in $1\times H\subset G$. This proves the desired assertion for $i=1$.

Assume that $i>1$ and that $[w_0,w_{i-1}]$ lies in~$G$. Being a commutator, this element lies in the kernel of $\sgn_2$ and $\sgn_3$, and therefore already in~$H^+$. By Lemma \ref{37Lem9} this implies that $[(w_0,w_0),w_i] = ([w_0,w_{i-1}],[w_0,w_{i-1}])$ lies in $H^+ \subset G$. Since $a_3$ lies in~$G$, and $G$ is normalized by $w_0$ and $w_i$, it follows that $[w_0,w_i] = [a_3(w_0,w_0),w_i]$ lies in $G$ as well.
Thus the lemma follows for all $i\ge1$ by induction on~$i$.
\end{Proof}

\begin{Lem}\label{37Lem06}
We have $N_4=G_4\cdot\langle w_0|_{T_4}\rangle$ and $[N_4:G_4]=2$.
\end{Lem}

\begin{Proof}
The defining relation $w_0=(a_2w_0,a_3w_0)$ recursively implies that
\begin{eqnarray*}
\sgn_1(w_0) &\!\!=\!\!& 1, \\
\sgn_2(w_0) &\!\!=\!\!& \sgn_1(a_2w_0)\cdot\sgn_1(a_3w_0)=1, \quad\hbox{and}\\
\sgn_3(w_0) &\!\!=\!\!& \sgn_2(a_2w_0)\cdot\sgn_2(a_3w_0)=-1.
\end{eqnarray*}
Thus the image of the residue class of $w_0|_{T_4}$ under the isomorphism from Lemma \ref{32Grplusone} (b) is the matrix
$$ \binom{\;\sgn_2(a_2w_0)\ \sgn_2(a_3w_0)\;}{\sgn_3(a_2w_0)\ \sgn_3(a_3w_0)} 
= \binom{-1\ \phantom{-}1}{-1\ \phantom{-}1}.$$
By Lemma \ref{32Grplusone} (c) this does not lie in the image of~$G_4$; hence $w_0|_{T_4} \not\in G_4$.
Also, since $N$ normalizes~$G$, its restriction $N_4$ normalizes~$G_4$; and hence Lemma \ref{32Grplusone} (d) implies that $[N_4:G_4]\le 2$. 
Since we already have $w_0|_{T_4} \in\nobreak N_4$ by Lemma \ref{37Lem03}, this leaves only the possibility $[N_4:G_4]=2$ with $N_4=G_4\cdot\langle w_0|_{T_4}\rangle$, as desired.
\end{Proof}

\medskip
Now we can put everything together. Lemmas \ref{37Lem11} (a) and \ref{37Lem03} imply that $w_i\in M$ for all $i\ge0$. By induction on $i$ we find that the restriction $w_i|_{T_n}$ is trivial for all $i\ge n$; hence the sequence $w_0,w_1,w_2,\ldots$ converges to $1$ within~$M$. Moreover, although in fact $w_0^2\not=1$, let us by abuse of notation for all $i\ge0$ and $k\in\BF_2$ define
$$w_i^k := 
\biggl\{\begin{array}{ll}
1   & \hbox{if $k=0$ in $\BF_2$,}\\[3pt]
w_i & \hbox{if $k=1$ in $\BF_2$.}
\end{array}$$
Then the following map is well-defined:
\UseTheoremCounterForNextEquation
\begin{equation}\label{37MapDef}
\textstyle \phi\colon
\prod\limits_{i=0}^\infty\BF_2 \longto M, \quad 
(k_0,k_1,k_2,\ldots) \mapsto w_0^{k_0}w_1^{k_1}w_2^{k_2}\cdots.
\end{equation}
By construction it is continuous, but not a homomorphism. By the construction (\ref{37WiDef}) of $w_1,w_2,\ldots$ it satisfies the basic formula
\UseTheoremCounterForNextEquation
\begin{equation}\label{37MapForm}
\phi(k_0,k_1,k_2,\ldots) 
\ =\ w_0^{k_0}\cdot\phi(0,k_1,k_2,\ldots) 
\ =\ w_0^{k_0}\cdot w_1^{k_1}\cdot\diag(\phi(0,k_2,k_3,\ldots)).
\end{equation}

\begin{Lem}\label{37Lem13}
The map $\phi$ induces a homomorphism $\bar\phi\colon \prod_{i=0}^\infty\BF_2 \to N/G$.
\end{Lem}

\begin{Proof}
It suffices to show that the images of all $w_i$ in $N/G$ have order dividing $2$ and commute with each other; in other words that $w_i^2$ and $[w_i,w_j]$ lie in~$G$ for all $i,j\ge0$.
Since $(a_sa_r)^2$ has order~$2$, the construction (\ref{37WiDef}) implies that $w_1^2=1$ and hence $w_i^2=1$ for all $i\ge1$. Also, the same proof as that of Lemma \ref{36Lem13} shows that $[w_i,w_j] \in G$ for all $i,j\ge1$. On the other hand, by Lemmas \ref{37Lem04} and \ref{37Lem05} the elements $w_0^2$ and $[w_0,w_i]$ lie in~$G$ for all $i\ge1$. Thus all cases are covered.
\end{Proof}

\begin{Lem}\label{37Lem14}
The homomorphism $\bar\phi$ is injective.
\end{Lem}

\begin{Proof}
Let $(k_0,k_1,k_2,\ldots)$ be an element of the kernel of~$\bar\phi$. Then $w := \phi(k_0,k_1,k_2,\ldots)$ lies in~$G$. 
In particular $w|_{T_4}$ lies in~$G_4$. But by Lemma \ref{37Lem11} (b), the definition of $\phi$ implies that $w|_{T_4} = w_0^{k_0}|_{T_4}$. Since $w_0|_{T_4} \not\in G_4$ by Lemma \ref{37Lem06}, it follows that $k_0=0$.
The same arguments as in the proof of Lemma \ref{36Lem15} now show that $k_1=0$ and that $(0,k_2,k_3,\ldots)$ also lies in the kernel of~$\bar\phi$. By an induction on $i$ we can deduce from this that for every $i\ge0$ and every element $(k_0,k_1,k_2,\ldots)$ of the kernel of~$\bar\phi$ we have $k_i=0$. This means that the kernel of $\bar\phi$ is trivial, and so the homomorphism $\bar\phi$ is injective, as desired.
\end{Proof}

\begin{Lem}\label{37Lem15}
The homomorphism $\bar\phi$ is surjective.
\end{Lem}

\begin{Proof}
We must prove that $N = G\cdot\phi\bigl(\prod_{i=0}^\infty\BF_2\bigr)$. For this it suffices to show the equality $N_n = G_n\cdot\phi\bigl(\prod_{i=0}^\infty\BF_2\bigr)|_{T_n}$ for all~$n\ge0$. For $n\le 4$ this follows from Lemma \ref{37Lem06}. So assume that $n>4$ and that $N_{n-1} = G_{n-1}\cdot\phi\bigl(\prod_{i=0}^\infty\BF_2\bigr)|_{T_{n-1}}$.

Since $n-1\ge4$, the same argument as in the proof of Lemma \ref{36Lem14} with Lemma \ref{37Lem41} in place of Proposition \ref{32MidCaseSquaresNotInH} shows that $G_{n-1}\cap M_{n-1} = H_{n-1}^+$. As $\phi\bigl(\prod_{i=0}^\infty\BF_2\bigr) \subset M$ by Lemmas \ref{37Lem11} (a) and \ref{37Lem03}, the induction hypothesis implies that 
\UseTheoremCounterForNextEquation
\begin{equation}\label{37Lem15Diag}
\textstyle M_{n-1} \ =\  
(G_{n-1}\cap M_{n-1})\cdot\phi\bigl(\prod\limits_{i=0}^\infty\BF_2\bigr)|_{T_{n-1}}
\ =\ H_{n-1}^+\cdot\phi\bigl(\prod\limits_{i=0}^\infty\BF_2\bigr)|_{T_{n-1}}.
\end{equation}
Using Lemma \ref{37Lem5}, the equality (\ref{37Lem15Diag}), the formula (\ref{37MapForm}), Lemma \ref{37Lem9}, the fact that $w_0=a_3(w_0,w_0)$, and rearranging factors because $N_n/G_n$ is commutative, we deduce that
\begin{eqnarray*}
N_n &=& G_n\cdot\langle w_1|_{T_n}\rangle\cdot\diag(M_{n-1}) \\
&=& \textstyle G_n\cdot\langle w_1|_{T_n}\rangle\cdot
\diag\bigl( H_{n-1}^+\cdot\phi\bigl(\prod_{i=0}^\infty\BF_2\bigr)|_{T_{n-1}} \bigr) \\
&=& \textstyle G_n\cdot\langle w_1|_{T_n}\rangle\cdot
\diag(H_{n-1}^+)\cdot\{1,(w_0,w_0)|_{T_n}\}
\cdot\diag\bigl(\phi\bigl(\{0\}\times\prod_{i=1}^\infty\BF_2\bigr)|_{T_{n-1}}\bigr) \\
&=& \textstyle G_n\cdot\diag(H_{n-1}^+)\cdot\{1,(w_0,w_0)|_{T_n}\}
\cdot\langle w_1|_{T_n}\rangle\cdot
\phi\bigl(\{0\}^2\times\prod_{i=2}^\infty\BF_2\bigr)|_{T_n} \\
&=& \textstyle G_n\cdot\{1,w_0|_{T_n}\}
\cdot\phi\bigl(\{0\}\times\prod_{i=1}^\infty\BF_2\bigr)|_{T_n} \\
&=& \textstyle G_n\cdot\phi\bigl(\prod_{i=0}^\infty\BF_2\bigr)|_{T_n},
\end{eqnarray*}
so the desired equality holds for~$n$. By induction it therefore follows for all $n\ge0$, as desired.
\end{Proof}

\medskip
Combining Lemmas \ref{37Lem13} through \ref{37Lem15} now implies:

\begin{Thm}\label{37NormThm}
In the case (\ref{37Ass}) the map $\phi$ induces an isomorphism 
$$\textstyle \bar\phi\colon \prod\limits_{i=0}^\infty\BF_2 \stackrel{\sim}{\longto} N/G.$$
\end{Thm}


\subsection{Odometers}
\label{38aOdometers}

In this subsection we study the set of odometers in~$G$, in particular concerning conjugacy. First we show that odometers are abundant in~$G$. By the \emph{proportion} of a measurable subset $X\subset G$ we mean the ratio $\mu(X)/\mu(G)$ for any Haar measure $\mu$ on~$G$. 

\begin{Prop}\label{38ManyOdos}
\begin{enumerate}
\item[(a)] For any $W$\!-conjugates $b_i$ of $a_i$ and any permutation $\rho$ of $\{1,\ldots,r\}$ the product $b_{\rho1}\cdots b_{\rho r}$ is an odometer.
\item[(b)] The set of odometers in $G$ is open and closed in $G$ and has proportion $2^{-r}$.
\end{enumerate}
\end{Prop}

\begin{Proof}
Same as for Proposition \ref{28ManyOdos}, using Proposition \ref{31GenSigns} instead of \ref{21GenSigns}.
\end{Proof}


\begin{Lem}\label{38OdosConjLem1}
If $r\ge3$, then any odometer in $G$ is conjugate under $N$ to an element of ${(Ha_s\times\{a_r\})\,\sigma}$.
\end{Lem}

\begin{Proof}
Any odometer $c$ in $G$ acts non-trivially on $T_1$ and is therefore equal to $(x,y)\,\sigma$ for some $(x,y)\in G^1$. Recall from Proposition \ref{32HAlmostDirect} (c) that $G^1 = (H\times H) \cdot \Langle (a_s,a_r), (a_r,a_s) \Rangle$. Since $\sgn_s(xy) = \sgn_{s+1}(c)= -1$ and $\sgn_s|H=1$, it follows that $(x,y)$ is an element of $H\times H$ times an odd number of factors $(a_s,a_r)$ and/or $(a_r,a_s)$. 
After possibly replacing $c$ by its conjugate under $\sigma=a_1\in G$, we may without loss of generality assume that the factor $(a_s,a_r)$ occurs an odd number of times and the factor $(a_r,a_s)$ an even number of times. Then $(x,y)\,(a_s,a_r)^{-1}$ is a product of an even number of each kind of factors. Since both $a_s$ and $a_r$ have order two, we deduce that
$$(x,y)\,(a_s,a_r)^{-1}
\ =\ \bigl((a_s,a_r)\,(a_r,a_s)\bigr)^{2\lambda}
\ =\ \bigl((a_sa_r)^{2\lambda},(a_ra_s)^{2\lambda}\bigr)
\ =\ (v^{\lambda},v^{-\lambda})$$
for some $\lambda\in\BZ_2$, where for brevity we write $v := (a_sa_r)^2$. Thus
$$c\ \in\ (H\times H)\cdot (v^{\lambda}a_s,v^{-\lambda}a_r)\cdot\sigma.$$
Next, if $s\ge2$ and $r\ge4$, then $v=(a_sa_r)^2\in H$ by Proposition \ref{32GenCaseSquaresInH}, and so $c$ already lies in $(H\times H)\cdot (a_s,a_r)\cdot\sigma$. Otherwise, the element $(1,v)=(1,(a_sa_r)^2)$ lies in $N$ by Lemma \ref{36Lem5} or \ref{37Lem5}, and so we may replace $c$ by its conjugate
\begin{eqnarray*}
(1,v^{\lambda})\,c\,(1,v^{-\lambda})
&\in& (1,v^{\lambda})\cdot(H\times H)\cdot (v^{\lambda}a_s,v^{-\lambda}a_r)\cdot\sigma\cdot (1,v^{-\lambda}) \\
&=& (H\times H)\cdot (v^{\lambda}a_sv^{-\lambda},a_r)\cdot\sigma \\
&=& (H\times H)\cdot (v^{2\lambda}a_s,a_r)\cdot\sigma.
\end{eqnarray*}
Since in this case $v^2 = (a_sa_r)^4=1$ by Lemma \ref{36Lem0} (a) or \ref{37Lem0} (a), we can then drop the term $v^{2\lambda}$. In either case have may thus without loss of generality assume that
$$c\ \in\ (H\times H)\cdot (a_s,a_r)\cdot\sigma.$$
Now write $c = (h'a_s,ha_r)\,\sigma$ with $h',h\in H$. Then $(1,h)\in G$; hence we may replace $c$ by its conjugate
$$(1,h)^{-1}\,(h'a_s,ha_r)\,\sigma\,(1,h) \ =\ (h'a_sh,a_r)\,\sigma.$$
Since $h'a_sh \in Ha_s$, the element $c$ then has the required form.
\end{Proof}

\begin{Lem}\label{38OdosConjLem2}
If $r\ge3$, for any $u\in N$ we have $(1,[a_r,u])\in N$.
\end{Lem}

\begin{Proof}
Since $u$ normalizes~$G$, the commutator $[a_r,u]$ lies in~$G$. Being a commutator, all its signs are~$+1$. If $s\ge2$ and $r\ge4$, by Lemma \ref{35Lem0} (c) this implies that $[a_r,u]\in H$; hence $(1,[a_r,u])$ lies in $G\subset N$, as desired. Otherwise Lemma \ref{36Lem0} (d) or \ref{37Lem0} (d) implies that $[a_r,u] \in H\cdot\Langle (a_sa_r)^2\Rangle$, and hence
$$(1,[a_r,u]) \ \in\ (1\times H)\cdot \Langle (1,(a_sa_r)^2)\Rangle.$$
But in this case $(1,(a_sa_r)^2) \in N$ by Lemma \ref{36Lem5} or \ref{37Lem5}, and since $1\times H\subset G\subset N$, we again deduce that $(1,[a_r,u])\in N$, as desired.
\end{Proof}

\begin{Thm}\label{38OdosConj}
All odometers in $G$ are conjugate under~$N$.
\end{Thm}

\begin{Proof}
In the case $s=1$ and $r=2$ recall from Subsection \ref{37Small} that 
$G = \Langle a_1a_2\Rangle \rtimes \langle a_1\rangle \cong \BZ_2\rtimes\{\pm1\}$ where $a_1a_2$ is an odometer. Thus the odometers in $G$ are precisely the elements $(a_1a_2)^k$ for all $k\in\BZ_2^\times$, and they are all conjugate to $a_1a_2$ by an element of $N\cong \BZ_2\ltimes\BZ_2^\times$. For the rest of the proof we therefore assume that $r\ge3$.
By Lemma \ref{13ConjLimit} it suffices to show for all $n\ge0$:
\begin{enumerate}
\item[($*_n$)] For any odometers $c,c'\in G$, there exists $w\in N$ with $wcw^{-1}|_{T_n} = c'|_{T_n}$. 
 \end{enumerate}
For $n\le r$ this is a direct consequence of the fact that $G_r=W_r$. So assume that $n>r$ and that ($*_{n-1}$) is true. 
Consider any odometers $c,c'\in G$. By Lemma \ref{38OdosConjLem1} we may without loss of generality assume that $c=(ha_s,a_r)\,\sigma$ and $c'=(h'a_s,a_r)\,\sigma$ for some $h,h'\in H$. Then 
$$c^2 = (ha_s,a_r)\,\sigma\,(ha_s,a_r)\,\sigma = (ha_sa_r,a_rha_s) \in G,$$ 
and so $d := ha_sa_r$ is again an odometer, which lies in $G$ by Proposition \ref{32G1Prop}. Likewise the element $d' := h'a_sa_r$ is an odometer in~$G$. Thus by the induction hypothesis ($*_{n-1}$) we can choose an element $u\in N$ with $udu^{-1}|_{T_{n-1}} = d'|_{T_{n-1}}$. For later use observe that we can replace $u$ by $ud^\lambda$ for any $\lambda\in\BZ_2$ without losing the equality. But for the moment we fix $u$ and set $w := (u,a_rua_r^{-1}) \in W$. Since $c=(da_r^{-1},a_r)\,\sigma$ and $c'=(d'a_r^{-1},a_r)\,\sigma$, we find that
$$wcw^{-1} 
\ =\ (u,a_rua_r^{-1})\,(da_r^{-1},a_r)\,\sigma\,(u^{-1},a_ru^{-1}a_r^{-1})
\ =\ (udu^{-1}a_r^{-1},a_r)\,\sigma$$
and deduce that
$$wcw^{-1}|_{T_n}
\ =\ \bigl((udu^{-1}a_r^{-1})|_{T_{n-1}},a_r|_{T_{n-1}}\bigr)\,\sigma
\ =\ \bigl((d'a_r^{-1})|_{T_{n-1}},a_r|_{T_{n-1}}\bigr)\,\sigma
\ =\ c'|_{T_n}.$$
To prove ($*_n$) it remains to show that we can choose $u\in N$ such that $w$ lies in~$N$. For this observe that 
$$w\ =\ (u,a_rua_r^{-1}) \ =\ (1,[a_r,u])\cdot(u,u).$$
By Lemma \ref{38OdosConjLem2} the first factor always lies in~$N$, so we must choose $u$ such that $(u,u)\in N$. If $s\ge2$ and $r\ge4$, this already follows from Lemma \ref{35Lem2}, and we are done. For the rest of the proof we therefore assume that $s=1$ or $r=3$. 

Since $n>r$, we have $n-1\ge r$. Let us postpone the case $(s,r,n) = (2,3,4)$ to the end of the proof and exclude it for now. Then Lemma \ref{36Lem0} and Proposition \ref{32MidCaseSquaresNotInH} if $s=1$, respectively Lemma \ref{37Lem0} and Proposition \ref{32SmallCaseSquaresInAndNotInH} if $s=2$, imply that that $G_{n-1} = H_{n-1}\rtimes\langle a_s|_{T_{n-1}},a_r|_{T_{n-1}}\rangle$ and that restriction to $T_{n-1}$ induces an isomorphism $\langle a_s,a_r\rangle \stackrel{\sim}{\to} \langle a_s|_{T_{n-1}},a_r|_{T_{n-1}}\rangle$. Thus the restriction induces an isomorphism $G/H \stackrel{\sim}{\to} G_{n-1}/H_{n-1}$. 

Recall that $M$ was defined as the subgroup of all elements of $N$ which induce the identity on~$G/H$, and that $N=MG$ by Lemma \ref{36Lem1} or \ref{37Lem1}. Recall also that $uha_sa_ru^{-1}|_{T_{n-1}} = udu^{-1}|_{T_{n-1}} = d'|_{T_{n-1}} = h'a_sa_r|_{T_{n-1}}$. It follows that $H_{n-1}(ua_sa_ru^{-1})|_{T_{n-1}} = H_{n-1}(a_sa_r)|_{T_{n-1}}$, in other words, the action of $u$ on $G_{n-1}/H_{n-1}$ induced by conjugation fixes the image of $a_sa_r$. The same is therefore true for the action of $u$ on $G/H$.
As $G/H$ is dihedral with cyclic part $\langle a_sa_r\rangle$ of order~$4$, this action is also induced by an element of $\langle a_sa_r\rangle$. Since $u\in N=MG$, it follows that $u\in M\langle a_sa_r\rangle$. But the fact that $d\in Ha_sa_r$ implies that $M\langle a_sa_r\rangle = M\langle d\rangle$. After replacing $u$ by $ud^\lambda$ for a suitable $\lambda\in\BZ_2$ we can therefore assume that $u\in M$. Then $(u,u)\in N$ by Lemma \ref{36Lem2} or \ref{37Lem2}, and we have established ($*_n$).

It remains to finish the induction step in the case $(s,r,n) = (2,3,4)$. For this we use the description of $G_{r+1}/(H_r\times H_r) \subset W_{r+1}/(H_r\times H_r) \cong \Mat_{2\times2}(\{\pm1\})\rtimes\langle\sigma\rangle$ from Lemma \ref{32Grplusone}. Since $c=(ha_s,a_r)\,\sigma$ and $c=(h'a_s,a_r)\,\sigma$ with $h,h'\in H$, the image of both elements under the map from \ref{32Grplusone} (b) is
$$\textstyle \binom{-1\ \phantom{-}1\;}{\phantom{-}1\ -1\;}\cdot\sigma.$$
Since any two odometers are conjugate under~$W$, we can choose an element $v\in W$ with $vcv^{-1}=c'$. Then the image of $v$ in $\Mat_{2\times2}(\{\pm1\})\rtimes\langle\sigma\rangle$ commutes with 
$\binom{-1\ \phantom{-}1\;}{\phantom{-}1\ -1\;}\cdot\sigma$.
By direct calculation the centralizer of this element is 
$$\textstyle \bigl\langle
\binom{-1\ -1\;}{\phantom{-}1\ \phantom{-}1\;},
\binom{\phantom{-}1\ \phantom{-}1\;}{-1\ -1\;},
\binom{-1\ \phantom{-}1\;}{\phantom{-}1\ -1\;}\cdot\sigma
\bigr\rangle.$$
{}From \ref{32Grplusone} (d) we see that this centralizer is contained in the normalizer of the image of~$G_{r+1}$. Moreover, Lemma \ref{37Lem06} implies that this normalizer is equal to the image of $N_{r+1}$ in $\Mat_{2\times2}(\{\pm1\})\rtimes\langle\sigma\rangle$. Thus the image of $N_{r+1}$ contains the image of~$v$. 
Since $r+1=4=n$, this is equivalent to $v|_{T_n} \in N_n$. But this in turn means that there exists a new element $w\in N$ with $w|_{T_n} = v|_{T_n}$. The equation $vcv^{-1}=c'$ then implies that $wcw^{-1}|_{T_n} = vcv^{-1}|_{T_n} = c'|_{T_n}$, as desired.
This proves ($*_n$) in the remaining case, finishing the proof of Theorem \ref{38OdosConj}.
\end{Proof}

\medskip
The author does not know which analogues of Theorem \ref{28OdosConj} (b) and Theorem \ref{28OdoSemiRigid}, if any, might hold in the strictly pre-periodic case. 
However, we can describe how the normalizer of any odometer in $G$ sits inside~$N$. 

\begin{Prop}\label{38OdoNorm1}
Consider any odometer $c\in G$ and any $k\in\BZ_2^\times$.
\begin{enumerate}
\item[(a)] There exists an element $w\in N$ with $wcw^{-1}=c^k$. 
\item[(b)] Every element $w\in W$ with $wcw^{-1}=c^k$ lies in~$N$.
\item[(c)] For any such $w$ the coset $wG$ depends only on~$k$.
\end{enumerate}
\end{Prop}

\begin{Proof}
As $c^k$ is another odometer in~$G$, part (a) follows from Theorem \ref{38OdosConj}. Next consider two elements $w,w'\in W$ with $wcw^{-1}=w'cw^{\prime-1}=c^k$. Then $w^{-1}w'\in W$ commutes with~$c$. By Proposition \ref{16OdoProp} (b) it therefore lies in $\Langle c\Rangle$ and hence in~$G$. Thus $wG=w'G$. Since $w\in N$ for some choice of~$w$, it follows that $w\in N$ for every choice, proving (b).  It also follows that the coset $wG$ is independent of~$w$. Moreover, as all odometers are conjugate under $N$ by Theorem \ref{38OdosConj}, the conjugacy class of $wG$ in $N/G$ is independent of~$c$. Since $N/G$ is abelian, it follows that the coset $wG$ itself is independent of~$c$, proving (c). 
\end{Proof}

\medskip
For the rest of this subsection we fix elements $c$, $k$, $w$ as in Proposition \ref{38OdoNorm1}. Like the description of $N$ itself, the determination of the coset $wG$ depends on the respective subcases. Until Proposition \ref{MonoRep} we assume that $r\ge3$. Recall that Theorem \ref{35NormThm}, \ref{36NormThm}, or \ref{37NormThm} then yields an isomorphism $\bar\phi$ between $N/G$ and $\prod_{i=1}^\infty\BF_2$ or $\prod_{i=0}^\infty\BF_2$. Let
\UseTheoremCounterForNextEquation
\begin{equation}\label{38CharDecomp}
\biggl\{\begin{array}{rl}
(k_1,k_2,\ldots) \!\!\!& \in \prod_{i=1}^\infty\BF_2, 
\quad\hbox{resp.} \\[3pt]
(k_0,k_1,k_2,\ldots) \!\!\!& \in \prod_{i=0}^\infty\BF_2,
\end{array}\biggr\}
\end{equation}
be the tuple corresponding to the coset~$wG$.

\begin{Lem}\label{38LemDiag}
The elements $k_i$ for $i\ge1$ are all equal.
\end{Lem}

\begin{Proof}
After possibly replacing $w$ by $wc$ we may without loss of generality assume that $w=(u,u')\in W\times W$. Write $c^2=(d,d')$. Then $d$ is another odometer in~$G$, and the equation $wcw^{-1}=c^k$ implies that $udu^{-1}=d^k$. By Proposition \ref{38OdoNorm1} we therefore have $u\in N$ and $uG=wG$.

Suppose first  that $(s,r)\not=(2,3)$. Since $w\in W\times W$, the choice (\ref{38CharDecomp}) then means that $w\in G^1\cdot\phi(k_1,k_2,\ldots)$. Note that by the definition (\ref{35WiDef}) or (\ref{36WiDef}) of $w_1$ we have $w_1\in W\times W$ and
$$\proj_1(w_1) \ =\ \biggl\{\begin{array}{ll}
\proj_1((a_s,a_s)) = a_s    & \hbox{if $s\ge2$ and $r\ge4$,} \\[3pt]
\proj_1((1,(a_sa_r)^2)) = 1 & \hbox{if $s = 1$ and $r\ge3$.}
\end{array}\biggr\}\ \in\ G.$$
Using the functional equation (\ref{35MapForm}) or (\ref{36MapForm}) and the fact that $\proj_1(G^1)=G$ by Proposition \ref{32G1Prop} (d), we deduce that
\begin{eqnarray*}
G\cdot\phi(k_1,k_2,\ldots)
&\!\!=\!\!& Gw \;=\; Gu \;=\; G\cdot\proj_1(w) \\
&\!\!=\!\!& 
G\cdot\proj_1\bigl(G^1\cdot w_1^{k_1}\cdot\diag(\phi(k_2,k_3,\ldots))\bigr) \\
&\!\!=\!\!& G\cdot\phi(k_2,k_3,\ldots).
\end{eqnarray*}
Therefore $k_1=k_2=k_3=\ldots$, as desired.

Suppose now that $(s,r)=(2,3)$. Then (\ref{38CharDecomp}) means that $w\in G^1\cdot\phi(k_0,k_1,k_2,\ldots)$. 
Also, by the definitions (\ref{37WiDef}) and (\ref{37W0Def}) of $w_1$ and $w_0$ we have $w_1,w_0\in W\times W$ and 
$$\biggl\{\begin{array}{l}
\proj_1(w_1) = \proj_1((1,(a_sa_r)^2)) = 1,\ \ \hbox{and} \\[3pt]
\proj_1(w_0) = \proj_1((a_2w_0,a_3w_0)) = a_2w_0 \in Gw_0.
\end{array}$$
Using the functional equation (\ref{37MapForm}) and the fact that $\proj_1(G^1)=G$, we deduce that
\begin{eqnarray*}
G\cdot\phi(k_0,k_1,k_2,\ldots)
&\!\!=\!\!& Gw \;=\; Gu \;=\; G\cdot\proj_1(w) \\
&\!\!=\!\!& 
G\cdot\proj_1\bigl(G^1\cdot w_0^{k_0}\cdot w_1^{k_1}\cdot\diag(\phi(0,k_2,k_3,\ldots))\bigr) \\
&\!\!=\!\!& G\cdot w_0^{k_0}\cdot\phi(0,k_2,k_3,\ldots) \\
&\!\!=\!\!& G\cdot\phi(k_0,k_2,k_3,\ldots).
\end{eqnarray*}
Therefore again $k_1=k_2=k_3=\ldots$, as desired.
\end{Proof}

\medskip
By Lemma \ref{38LemDiag} it remains to determine~$k_1$, respectively $k_1$ and~$k_0$. The first tool for this is a closer look at the subgroups $G_{r+1}\subset N_{r+1} \subset W_{r+1}$ in the next two lemmas:

\begin{Lem}\label{38LemWiInGrPlusOne}
\begin{enumerate}
\item[(a)] If $s=2$ and $r=3$, we have $N_4 = G_4\cdot\langle w_0|_{T_4}\rangle$
and $[N_4:G_4]=2$ and $w_i|_{T_4}\in G_4$ for all $i\ge1$.
\item[(b)] If $s\ge2$ and $r\ge4$, we have $N_{r+1} = G_{r+1}\cdot\langle w_1|_{T_{r+1}}\rangle$ and $[N_{r+1}:G_{r+1}]=2$ and $w_i|_{T_{r+1}}\in G_{r+1}$ for all $i\ge2$.
\end{enumerate}
\end{Lem}

\begin{Proof}
Part (a) is just a repetition of Lemmas \ref{37Lem11} (b) and \ref{37Lem06}. 
For (b) recall from the construction (\ref{35WiDef}) that $w_0=(a_s,a_s)$; hence the image of the residue class of $w_0|_{T_{r+1}}$ under the isomorphism from Lemma \ref{32Grplusone} (b) is the matrix
$$ \binom{\;\sgn_s(a_s)\ \ \sgn_s(a_s)\;}{\sgn_r(a_s)\ \ \sgn_r(a_s)} 
= \binom{{-1}\ \ {-1}\,}{\phantom{-}1\ \ \phantom{-}1\,}.$$
By Lemma \ref{32Grplusone} (c) this matrix does not lie in the image of $G_{r+1}$, which proves that $w_0|_{T_{r+1}} \not\in\nobreak G_{r+1}$.
Also, since $N$ normalizes~$G$, its restriction $N_{r+1}$ normalizes~$G_{r+1}$; and hence Lemma \ref{32Grplusone} (d) implies that $[N_{r+1}:G_{r+1}]\le 2$. Since $w_1|_{T_{r+1}}$ already lies in~$N_{r+1}$, but not in~$G_{r+1}$, this leaves only the possibility $[N_{r+1}:G_{r+1}]=2$ with $N_{r+1}=G_{r+1}\cdot\langle w_1|_{T_{r+1}}\rangle$.
Finally, for all $i\ge2$ the construction (\ref{35WiDef}) implies that $w_i = (w_{i-1},w_{i-1})$ with $w_{i-1} \in \diag(W)$ and hence $\sgn_s(w_{i-1})=\sgn_r(w_{i-1})=1$. Thus the image of the residue class of $w_i|_{T_{r+1}}$ under the isomorphism from Lemma \ref{32Grplusone} (b) is the trivial matrix, which shows that $w_i|_{T_{r+1}}$ lies in~$G_{r+1}$.
All parts of (b) are now proved.
\end{Proof}

\begin{Lem}\label{38LemOdo}
If $s\ge2$, then $w|_{T_{r+1}} \in G_{r+1}$ if and only if $\frac{k-1}{2}$ is even.
\end{Lem}

\begin{Proof}
Abbreviate $\ell := \frac{k-1}{2}\in\BZ_2$. As a preparation let $a\in W$ be the standard odometer from (\ref{16OdoRec}) and $z\in W$ the element satisfying ${zaz^{-1}=a^k}$ from Remark \ref{16OdoNormExpl}. Recall that $z$ was characterized by the recursion relation $z=(z,a^\ell z)$. Thus we have $\sgn_1(z)=1$, and for all $n\ge2$ we can calculate
$$\sgn_n(z) = \sgn_{n-1}(z)\cdot\sgn_{n-1}(a^\ell z) = \sgn_{n-1}(a)^\ell = (-1)^\ell.$$
Since $r>s\ge2$, the image of $z|_{T_{r+1}}$ in $\Mat_{2\times2}(\{\pm1\})\rtimes\langle\sigma\rangle$ under the map from Lemma \ref{32Grplusone} (b) is therefore the matrix
$$ \binom{\;\sgn_s(z)\ \ \sgn_s(a^\ell z)\;}{\sgn_r(z)\ \ \sgn_r(a^\ell z)} 
= \binom{(-1)^\ell\ \ 1}{(-1)^\ell\ \ 1}.$$
The conjugacy class of this element in $\Mat_{2\times2}(\{\pm1\})\rtimes\langle\sigma\rangle$ is the subset
$$X\ :=\ \biggl\{
\binom{(-1)^\ell\ \ 1}{(-1)^\ell\ \ 1},\ 
\binom{1\ \ (-1)^\ell}{1\ \ (-1)^\ell}
\biggr\}.$$
Now recall that $a$ is conjugate to $c$ under~$W$. Proposition \ref{16OdoProp} therefore implies that $z$ is conjugate to $wc^m$ under $W$ for some $m\in\BZ_2$. It follows that the image of $wc^m|_{T_{r+1}}$ in $\Mat_{2\times2}(\{\pm1\})\rtimes\langle\sigma\rangle$ is contained in~$X$. Using the description of the image of $G_{r+1}$ from Lemma \ref{32Grplusone} (c), we find that $X$ is a subset of the image of $G_{r+1}$ if $\ell$ is even, and disjoint from it if $\ell$ is odd. Thus the image of $wc^m|_{T_{r+1}}$ is contained in the image of $G_{r+1}$ if and only if $\ell$ is even. This means that $wc^m|_{T_{r+1}} \in G_{r+1}$ if and only if $\ell$ is even. As $c|_{T_{r+1}}$ itself lies in $G_{r+1}$, it follows that $w|_{T_{r+1}} \in G_{r+1}$ if and only if $\ell$ is even, as desired.
\end{Proof}

\medskip
The second tool is a closer study of $G/(H\times H)$ in the following two lemmas:

\begin{Lem}\label{38LemWiInGModHTimesH}
Assume that either $s=1$ and $r\ge3$, or $s=2$ and $r=3$. Then:
\begin{enumerate}
\item[(a)] 
$G = (H\times H)\rtimes\langle\sigma,(a_s,a_r)\rangle$
where $\langle\sigma,(a_s,a_r)\rangle$ is a dihedral group with the maximal cyclic subgroup $\langle\sigma(a_s,a_r)\rangle$ of order~$8$. 
\item[(b)] The subgroup $H\times H$ is normalized by~$N$.
\item[(c)] Conjugation by $w_1$ induces the automorphism $x\mapsto x^5$ on the cyclic subgroup of order~$8$ of $G/(H\times H)$.
\item[(d)] For every $i\ge2$, conjugation by $w_i$ induces the trivial automorphism on $G/(H\times H)$.
\item[(e)] If $s=2$ and $r=3$, conjugation by $w_0$ induces an inner automorphism on $G/(H\times H)$.
\end{enumerate}
\end{Lem}

\begin{Proof}
The semidirect product decomposition in (a) is a reminder from Lemmas \ref{36Lem0} (c) and \ref{37Lem0} (c), and the element $\sigma\,(a_s,a_r) = a_1a_{s+1}$ has order $8$ by Proposition \ref{32aiajOrder} in both cases, proving (a). 

Next, by (\ref{36WiDef}) and (\ref{37WiDef}) we have $w_1 = (1,(a_sa_r)^2)$ in both cases. Since $(a_sa_r)^2 = (a_ra_s)^2$ generates the center of $\langle a_s,a_r\rangle$ by Lemma \ref{36Lem0} (a) or \ref{37Lem0} (a), it follows that $w_1$ commutes with $(a_s,a_r)$. On the other hand we calculate
$$(\sigma\,(a_s,a_r))^4
= (a_ra_s,a_sa_r)^2 
= ((a_ra_s)^2,(a_sa_r)^2) 
= ((a_sa_r)^2,(a_sa_r)^2) 
= [w_1,\sigma].$$
Together this implies that 
$$w_1\cdot\sigma\,(a_s,a_r)\cdot w_1^{-1} 
= w_1\sigma w_1^{-1}\sigma^{-1} \cdot \sigma\,(a_s,a_r) 
= (\sigma\,(a_s,a_r))^4\,\sigma\,(a_s,a_r)
= (\sigma\,(a_s,a_r))^5.$$
Also, since $(a_sa_r)^2$ normalizes~$H$, the element $w_1$ normalizes $H\times H$, and together this shows (c).

By (\ref{36WiDef}) and (\ref{37WiDef}) in both cases we also have $w_i = (w_{i-1},w_{i-1})$ for all $i\ge2$. This element trivially commutes with~$\sigma$, and since conjugation by $w_{i-1} \in M$ normalizes $H$ and acts trivially on~$G/H$, it follows that conjugation by $w_i$ normalizes $H\times H$ and fixes the coset $(a_s,a_r)(H\times H)$. This implies (d).

In the case $s=2$ and $r=3$ we have $w_0 = a_3\,(w_0,w_0)$ by the construction (\ref{37W0Def}). Also since $w_0\in M$ by Lemma \ref{37Lem03}, conjugation by $w_0$ normalizes $H$ and acts trivially on~$G/H$. Thus conjugation by $(w_0,w_0)$ normalizes $H\times H$ and fixes the coset ${(a_s,a_r)}{(H\times H)}$. Since $(w_0,w_0)$ also commutes with~$\sigma$, it follows that conjugation by it induces the trivial automorphism of $G/(H\times H)$. This implies that $w_0 = a_3\,(w_0,w_0)$ normalizes $H\times H$ and induces the same inner automorphism of $G/(H\times H)$ as~$a_3$, proving (e).

Finally, since $N$ is generated by $G$ and all the elements $w_i$ considered, and all these elements normalize $H\times H$, this implies (b).
\end{Proof}

\begin{Lem}\label{38LemUff}
In the situation of Lemma \ref{38LemWiInGModHTimesH}, the image of $c$ in the dihedral group $G/(H\times H)$ generates the maximal cyclic subgroup of order~$8$.
\end{Lem}

\begin{Proof}
Since all odometers in $G$ are conjugate under~$N$ by Theorem \ref{38OdosConj}, and $N$ normalizes $H\times H$ by Lemma \ref{38LemWiInGModHTimesH} (b), it suffices to prove this for any one choice of~$c$. We will prove it for $c=a_1\cdots a_r$, which is an odometer in $G$ by Proposition \ref{38ManyOdos} (a). Recall that $a_i=(a_{i-1},1)\in H\times 1$ for all $i\not=1,s+1$. Thus the image of $a_1\cdots a_r$ in $G/(H\times H)$ is equal to the image of $a_1a_{s+1}=\sigma\,(a_s,a_r)$. By Lemma \ref{38LemWiInGModHTimesH} (a) this image generates the maximal cyclic subgroup of $G/(H\times H)$, as desired.
\end{Proof}

\medskip
Now we can reap the fruits of the preceding efforts:

\begin{Lem}\label{38LemMod48Char}
\strut
\begin{enumerate}
\item[(a)] If $s\ge2$ and $r\ge4$, we have $k_1\equiv\frac{k-1}{2}\mod 2$.
\item[(b)] If $s=1$ and $r\ge3$, we have $k_1\equiv\frac{k^2-1}{8}\mod 2$.
\item[(c)] If $s=2$ and $r=3$, we have $k_0\equiv\frac{k-1}{2}\mod 2$.
\item[(d)] If $s=2$ and $r=3$, we have $k_1\equiv\frac{k^2-1}{8}\mod 2$.
\end{enumerate}
\end{Lem}

\begin{Proof}
If $s\ge2$ and $r\ge4$, the definition (\ref{35MapDef}) of $\phi$ implies that 
$$w\in G\cdot\phi(k_1,k_2,\ldots) = G\, w_1^{k_1}w_2^{k_2}\cdots.$$
By Lemma \ref{38LemWiInGrPlusOne} (b) it follows that $w|_{T_{r+1}} \in G_{r+1}\cdot (w_1|_{T_{r+1}})^{k_1}$ with $w_1|_{T_{r+1}}\not\in G_{r+1}$. But by Lemma \ref{38LemOdo} we have  $w|_{T_{r+1}} \in G_{r+1}$ if and only if $\frac{k-1}{2}\equiv0\mod 2$. Thus $k_1\equiv\frac{k-1}{2}\mod 2$, proving (a).

\medskip
If $s=1$ and $r\ge3$, the definition (\ref{36MapDef}) of $\phi$ again implies that 
$$w\in G\cdot\phi(k_1,k_2,\ldots) = G\, w_1^{k_1}w_2^{k_2}\cdots.$$
{}From Lemma \ref{38LemWiInGModHTimesH} (c--d) we find that conjugation by $w_1^{k_1}w_2^{k_2}\cdots$ induces the automorphism $x\mapsto x^{1+4k_1}$ on the maximal cyclic subgroup of order $8$ of ${G/(H\times H)}$.
Moreover, as $G/(H\times H)$ is a dihedral group, the automorphism of the maximal cyclic subgroup induced by any element of $G$ has the form $x\mapsto x^{\pm1}$.
On the other hand, by Lemma \ref{38LemUff} this subgroup is generated by the image of~$c$, and so conjugation by $w$ induces the automorphism $x\mapsto x^k$ on it.
Together this implies that $x^{\pm(1+4k_1)}=x^k$, or equivalently $1+4k_1 \equiv \pm k \mod 8$. Therefore $k_1\equiv\frac{k^2-1}{8}\mod 2$, proving (b).

\medskip
If $s=2$ and $r=3$, the definition (\ref{37MapDef}) of $\phi$ implies that 
$$w\in G\cdot\phi(k_0,k_1,k_2,\ldots) = G\,w_0^{k_0}w_1^{k_1}w_2^{k_2}\cdots.$$
By Lemma \ref{38LemWiInGrPlusOne} (a) it follows that $w|_{T_4} \in G_4\cdot (w_0|_{T_4})^{k_0}$ with $w_0|_{T_4}\not\in G_4$. But by Lemma \ref{38LemOdo} we have  $w|_{T_4} \in G_4$ if and only if $\frac{k-1}{2}\equiv0\mod 2$. Thus $k_0\equiv\frac{k-1}{2}\mod 2$, proving~(c).

Moreover, from Lemma \ref{38LemWiInGModHTimesH} (c--d) we find that conjugation by $w_1^{k_1}w_2^{k_2}\cdots$ induces the automorphism $x\mapsto x^{1+4k_1}$ on the maximal cyclic subgroup of order $8$ of ${G/(H\times H)}$.
With Lemma \ref{38LemWiInGModHTimesH} (e) we also deduce that the automorphism induced by any element of $Gw_0^{k_0}$ has the form $x\mapsto x^{\pm1}$.
Continuing exactly as in the proof of (b) now shows (d).
\end{Proof}

\medskip
To combine everything consider the two homomorphisms
\UseTheoremCounterForNextEquation
\begin{equation}\label{MonoRep}
\biggl\{\begin{array}{l}
\theta_1\colon\;\BZ_2^\times \onto\BF_2,\ k\mapsto \frac{k-1}{2}\mod 2, \\[3pt]
\theta_2\colon\;\BZ_2^\times \onto\BF_2,\ k\mapsto \frac{k^2-1}{8}\mod 2.
\end{array}\biggr\}
\end{equation}
Recall that Theorem \ref{35NormThm}, \ref{36NormThm}, or \ref{37NormThm}, respectively Proposition \ref{37OdoNormExpl} yields an isomorphism between $N/G$ and $\prod_{i=1}^\infty\BF_2$ or $\prod_{i=0}^\infty\BF_2$, respectively $\BZ_2^\times / \{\pm1\}$. 

\begin{Prop}\label{38OdoNorm}
Consider any odometer $c\in G$ and any element $w\in N$ with $wcw^{-1}=c^k$ for some $k\in\BZ_2^\times$. Then the coset $wG \in N/G$ corresponds to the element
$$\left\{\begin{array}{rll}
(\theta_1,\theta_1,\ldots)(k) \!\! & 
\in\, \prod_{i=1}^\infty\BF_2 \,\cong\, N/G &
\hbox{if $s\ge2$ and $r\ge4$,} \\[5pt]
(\theta_2,\theta_2,\ldots)(k) \!\! & 
\in\, \prod_{i=1}^\infty\BF_2 \,\cong\, N/G &
\hbox{if $s=1$ and $r\ge3$,} \\[5pt]
(\theta_1,\theta_2,\theta_2,\ldots)(k) \!\! &
\in\, \prod_{i=0}^\infty\BF_2 \,\cong\, N/G &
\hbox{if $s=2$ and $r=3$,} \\[5pt]
k\;\mod\;\{\pm1\} \!\! &
\in\, \BZ_2^\times / \{\pm1\} \,\cong\, N/G &
\hbox{if $s=1$ and $r=2$.}
\end{array}\right\}$$
\end{Prop}

\begin{Proof}
The first three cases follow by combining the formula (\ref{38CharDecomp}) with Lemmas \ref{38LemDiag} and \ref{38LemMod48Char}. The last case is a direct consequence of Proposition \ref{37OdoNormExpl}.
\end{Proof}


\subsection{Iterated monodromy groups}
\label{38Monodromy}

Now let $G^\geom\subset G^\arith\subset W$ be the geometric and arithmetic fundamental groups associated to a quadratic polynomial over~$k$, as described in Subsection \ref{17Monodromy}, and assume that we are in the strictly pre-periodic case \ref{17Cases} (c). Then by (\ref{17GgeomGensPol}) we have $G^\geom = \Langle b_{p_1},\ldots,b_{p_r}\Rangle$, where by Proposition \ref{17GenRec} the generators satisfy the recursion relations
\UseTheoremCounterForNextEquation
\begin{equation}\label{38bRecRels}
\left\{\begin{array}{l}
\hbox{$b_{p_1}$ is conjugate to $\sigma$ under~$W$,}\\[3pt]
\hbox{$b_{p_{s+1}}$ is conjugate to $(b_{p_s},b_{p_r})$ under~$W$, and}\\[3pt]
\hbox{$b_{p_i}$ is conjugate to $(b_{p_{i-1}},1)$ under~$W$ for all $2\le i\le r$ with $i\not=s+1$.}
\end{array}\right\}
\end{equation}
By Theorem \ref{34SemiRigid} we deduce:

\begin{Thm}\label{38GgeomThm}
There exists $w\in W$ such that $G^\geom=wGw^{-1}$ and $b_{p_i}$ is conjugate to $wa_iw^{-1}$ under $G^\geom$ for every $1\le i\le r$.
\end{Thm}

For the following we change the identification of trees used in Subsection \ref{17Monodromy} by the automorphism~$w$, after which we have $G^\geom=G$, and $b_{p_i}$ is conjugate to $a_i$ under $G$ for every $1\le i\le r$.
Then $G^\arith$ is contained in the normalizer~$N$, and to describe it it suffices to describe the factor group $G^\arith/G \subset N/G$. More precisely we will determine the composite homomorphism
\UseTheoremCounterForNextEquation
\begin{equation}\label{38GalRep}
\bar\rho\colon \Gal(\bar k/k) \onto G^\arith/G \into N/G
\end{equation}
obtained from (\ref{17FactorMonoRep}). Recall from (\ref{17GgeomGensRelPol}) that 
$b_\infty = (b_{p_1}\!\cdots b_{p_r})^{-1}$, and consider the cyclotomic character
$$\cycl\colon \Gal(\bar k/k) \longto \BZ_2^\times,\ \tau\mapsto\cycl(\tau).$$

\begin{Lem}\label{38LemCyclEnough}
The element $b_\infty$ is an odometer in~$G$, and for any $\tau\in\Gal(\bar k/k)$ we have $\bar\rho(\tau) = wG$ for some element $w\in N$ which satisfies $wb_\infty w^{-1} = b_\infty^{\cycl(\tau)}$.
\end{Lem}

\begin{Proof}
By the construction in Subsection \ref{17Monodromy} the element $b_\infty$ is a generator of the image of an inertia group $I_\infty \subset\pi_1^\et(\BP^1_{\bar k} \setminus P,x_0)$ above $\infty\in\BP^1_k$. Let $D_\infty$ be the associated decomposition group. Since $\infty$ is a $k$-rational point, we have a commutative diagram with exact rows
$$\xymatrix{
1 \ar[r] & I_\infty \ar@{^{ (}->}[d]
\ar[r] & D_\infty \ar@{^{ (}->}[d]
\ar[r] & \Gal(\bar k/k) \ar[r] \ar@{=}[d] &  1 \\
1 \ar[r] & \pi_1^\et(\BP^1_{\bar k} \setminus P,x_0) \ar@{->>}[d]_\rho
\ar[r] & \pi_1^\et(\BP^1_k \setminus P,x_0) \ar@{->>}[d]_\rho
\ar[r] & \Gal(\bar k/k) \ar[r] \ar@{->>}[d]_{\bar\rho} & 1 \\
1 \ar[r] & G
\ar[r] & G^\arith
\ar[r] & G^\arith/G \ar[r] &  1\rlap{.} \\}$$
Since the target of the homomorphism $\rho$ is a pro-$2$ group, the restriction $\rho|_{I_\infty}$ factors through the maximal pro-$2$ quotient of~$I_\infty$. As the characteristic of $k$ is not~$2$, this quotient is the pro-$2$ tame inertia group above $\infty$ and hence isomorphic~$\BZ_2$ and generated by~$b_\infty$. Moreover, $D_\infty/I_\infty \cong \Gal(\bar k/k)$ acts on it through the cyclotomic character. Thus for any $\delta\in D_\infty$ with image $\tau\in\Gal(\bar k/k)$ the element $w:=\rho(\delta) \in G^\arith\subset N$ satisfies $wb_\infty w^{-1} = b_\infty^{\cycl(\tau)}$. 
Since $b_\infty$ is an odometer by Proposition \ref{38ManyOdos}, the lemma follows.
\end{Proof}

\medskip
Finally recall that Theorem \ref{35NormThm}, \ref{36NormThm}, or \ref{37NormThm}, respectively Proposition \ref{37OdoNormExpl} yields an isomorphism between $N/G$ and $\prod_{i=1}^\infty\BF_2$ or $\prod_{i=0}^\infty\BF_2$, respectively $\BZ_2^\times / \{\pm1\}$. By combining Lemma \ref{38LemCyclEnough} with Proposition \ref{38OdoNorm} we deduce:

\begin{Thm}\label{38GarithThm1}
The homomorphism $\bar\rho\colon \Gal(\bar k/k) \to N/G$ has the form
$$\left\{\begin{array}{rll}
(\theta_1,\theta_1,\ldots)\circ\cycl :\!\! &
\Gal(\bar k/k) \stackrel{\bar\rho}{\longto} N/G \,\cong\, \prod_{i=1}^\infty\BF_2 & 
\hbox{if $s\ge2$ and $r\ge4$,} \\[5pt]
(\theta_2,\theta_2,\ldots)\circ\cycl :\!\! &
\Gal(\bar k/k) \stackrel{\bar\rho}{\longto} N/G \,\cong\, \prod_{i=1}^\infty\BF_2 & 
\hbox{if $s=1$ and $r\ge3$,} \\[5pt]
(\theta_1,\theta_2,\theta_2,\ldots)\circ\cycl :\!\! &
\Gal(\bar k/k) \stackrel{\bar\rho}{\longto} N/G \,\cong\, \prod_{i=0}^\infty\BF_2 & 
\hbox{if $s=2$ and $r=3$,} \\[5pt]
\cycl\;\mod\{\pm1\} :\!\! &
\Gal(\bar k/k) \stackrel{\bar\rho}{\longto} N/G \,\cong\, \BZ_2^\times / \{\pm1\} & 
\hbox{if $s=1$ and $r=2$.}
\end{array}\right\}$$
\end{Thm}

\begin{samepage}
\begin{Cor}\label{38GarithThm2}
\strut
\begin{enumerate}
\item[(a)] If $s\ge2$ and $r\ge4$, the order of $G^\arith/G$ divides $2$.
\item[(b)] If $s\ge2$ and $r\ge4$ and $k$ contains the fourth roots of unity, then $G^\arith=G$.
\item[(c)] If $s=1$ and $r\ge3$, the order of $G^\arith/G$ divides $2$.
\item[(d)] If $s=1$ and $r\ge3$ and $k$ contains a square root of~$2$, then $G^\arith=G$.
\item[(e)] If $s=2$ and $r=3$, the order of $G^\arith/G$ divides $4$.
\item[(f)] If $s=2$ and $r=3$ and $k$ contains the eighth roots of unity, then $G^\arith=G$.
\item[(g)] If $s=1$ and $r=2$ and $k$ is finitely generated over its prime field, then $G^\arith/G$ a subgroup of finite index of $\BZ_2^\times / \{\pm1\}$.
\end{enumerate}
\end{Cor}
\end{samepage}

\begin{Proof}
Parts (a--f) are direct consequences of Theorem \ref{38GarithThm1}; part (g) follows by the same argument as in the proof of Theorem \ref{28GarithThm} (b).
\end{Proof}


%
%
%
%
%
%
%

\newpage

\addcontentsline{toc}{section}{References}



\end{document}